\documentclass[aap]{imsart}

\RequirePackage{amsthm,amsmath,amsfonts,amssymb,color}
\RequirePackage[numbers,sort&compress]{natbib}
\RequirePackage[colorlinks,citecolor=blue,urlcolor=blue]{hyperref}
\RequirePackage{graphicx}

\startlocaldefs
\theoremstyle{plain}

\newtheorem{theorem}{Theorem}[section]
\newtheorem{lemma}[theorem]{Lemma}
\newtheorem{proposition}[theorem]{Proposition}
\newtheorem{corollary}[theorem]{Corollary}
\newtheorem{conjecture}[theorem]{Conjecture}
\theoremstyle{definition}
\newtheorem{definition}[theorem]{Definition}
\newtheorem{remark}[theorem]{Remark}

\usepackage{etoolbox}
\usepackage{enumerate}
\usepackage{verbatim}
\usepackage{mathrsfs}
\usepackage{centernot}

\allowdisplaybreaks[1]
\date{\today}

\def\E#1{\mathbb{E}}
\def\supp{\operatorname{supp}}

\def\Ent{\operatorname{Ent}}

\def\d{{\rm d}}

\def\conc{{\rm conc\,}}

\makeatother

\endlocaldefs

\begin{document}

\begin{frontmatter}
\title{Local Optimality of Dictator Functions  with Applications to Courtade--Kumar
and Li--M\'edard Conjectures\thanks{\textbf{A\MakeLowercase{ccepted for publication in the} A\MakeLowercase{nnals of} A\MakeLowercase{pplied} P\MakeLowercase{robability}. T\MakeLowercase{his extended version includes additional proofs, auxiliary formulas, and comprehensive Matlab scripts not present in the published version.}}}
}
\runtitle{Local Optimality of Dictator Functions}

\begin{aug}
\author[A]{\fnms{Lei}~\snm{Yu} \ead[label=e1]{leiyu@nankai.edu.cn}}
\address[A]{School of Statistics and Data Science, LPMC, KLMDASR, and LEBPS, 
Nankai University, Tianjin 300071, China\printead[presep={ ,\ }]{e1}}
\end{aug}

\begin{abstract}
Given a convex function $\Phi:[0,1]\to\mathbb{R}$, the $\Phi$-stability
of a Boolean function $f$ is defined as    $\mathbb{E}[\Phi(T_{\rho}f(\mathbf{X}))]$,
where $\mathbf{X}$ is a random vector uniformly distributed on the
discrete cube $\{\pm1\}^{n}$ and $T_{\rho}$ is the Bonami-Beckner
operator. In this paper, we prove that dictator functions are locally
optimal in maximizing the $\Phi$-stability of $f$ over all balanced
Boolean functions. When focusing on the symmetric $q$-stability,
combining this result with our previous bound, we use computer-assisted
methods to prove that dictator functions maximize the symmetric $q$-stability
for $q=1$ and $\rho\in[0,0.914]$ or for $q\in[1.36,2)$ and all
$\rho\in[0,1]$. In other words, we confirm  the (balanced) Courtade--Kumar
conjecture with the correlation coefficient $\rho\in[0,0.914]$ and
the (symmetrized) Li--M\'edard conjecture with $q\in[1.36,2)$. We
conjecture that dictator functions maximize both the symmetric and
asymmetric $\frac{1}{2}$-stability over all balanced Boolean functions.
Our proofs are based on majorization of noise operators and hypercontractivity
inequalities.
\end{abstract}



\end{frontmatter}

\section{\label{sec:Introduction}Introduction}

Let $\mathbf{X}$ be a random vector uniformly distributed on the
discrete cube $\{\pm1\}^{n}$. Let $\mathbf{Y}\in\{\pm1\}^{n}$ be
the random vector obtained by independently changing the sign of each
component of $\mathbf{X}$ with the same probability $\frac{1-\rho}{2}$.
Here, $\rho\in[0,1]$ corresponds to the correlation coefficient between
each component of $\mathbf{X}$ and the corresponding one of $\mathbf{Y}$.
Let $T_{\rho}=T_{\rho}^{(n)}$ be the noise operator (also known as
the \emph{Bonami-Beckner operator}) which acts on Boolean functions
$f:\{\pm1\}^{n}\to\{0,1\}$ such that $T_{\rho}f(\mathbf{x})=\mathbb{E}[f(\mathbf{Y})|\mathbf{X}=\mathbf{x}]$.
Let $\Phi:[0,1]\to\mathbb{R}$ be a continuous and strictly convex
function. 
\begin{definition}
For a Boolean function $f:\{\pm1\}^{n}\to\{0,1\}$, the \emph{$\Phi$-stability}
of $f$ with respect to (w.r.t.) correlation coefficient $\rho$ is
defined as 
\begin{align*}
\mathbf{Stab}_{\Phi}[f] & =\mathbb{E}[\Phi(T_{\rho}f(\mathbf{X}))].
\end{align*}
\end{definition}

The noise stability problem, in a general sense, concerns which Boolean
functions (or measurable sets) are the ``most stable'' under the
action of the noise operator. In terms of $\Phi$-stability, the noise
stability problem is formulated as 
\[
\max_{f:\{\pm1\}^{n}\to\{0,1\},\mathbb{E}f=\alpha}\mathbf{Stab}_{\Phi}[f].
\]
It is known (e.g., in \cite{yu2023phi}) that 
\[
I_{\Phi}(f(\mathbf{Y});\mathbf{X})=\mathbf{Stab}_{\Phi}[f]-\Phi(\mathbb{E}f),
\]
where $I_{\Phi}(f(\mathbf{Y});\mathbf{X}):=\Ent_{\Phi}(f(\mathbf{Y}))-\Ent_{\Phi}(f(\mathbf{Y})|\mathbf{X})$
is the $\Phi$-mutual information between $f(\mathbf{Y})$ and $\mathbf{X}$.
Here, $\Ent_{\Phi}(f(\mathbf{Y})):=\mathbb{E}[\Phi(f(\mathbf{Y}))]-\Phi(\mathbb{E}f(\mathbf{Y}))$
is the \emph{$\Phi$-entropy} of $f(\mathbf{Y})$, and $\Ent_{\Phi}(f(\mathbf{Y})|\mathbf{X}):=\mathbb{E}[\Phi(f(\mathbf{Y}))]-\mathbb{E}\Phi(\mathbb{E}[f(\mathbf{Y})|\mathbf{X}])$
is the \emph{conditional $\Phi$-entropy} of $f(\mathbf{Y})$ given
$\mathbf{X}$. So, maximizing the $\Phi$-stability w.r.t. volume
$\alpha\in[0,1]$ is equivalent to maximizing the $\Phi$-mutual information
between $f(\mathbf{Y})$ and $\mathbf{X}$ over all Boolean functions
with $\mathbb{E}f=\alpha$; see details in \cite{yu2023phi}.

In the learning languages, the noise stability problem admits a natural interpretation. Let $\mathbf{X}$ denote the data we are truly interested in, and let $\mathbf{Y}$ denote its noisy counterpart---the data actually accessible to us. We aim to learn a single bit of information about $\mathbf{X}$ but only using $\mathbf{Y}$. A key question then arises: what is the most informative bit about $\mathbf{X}$ that we can extract from $\mathbf{Y}$? For specific functions $\Phi$, it was conjectured in the literature that the dictator function is optimal for this purpose. This conjecture is known as the ``most informative bit'' conjecture, and it will be introduced in detail in the following subsection.

\subsection{Related Conjectures}

Consider two instances of $\Phi$. For $q>0$, define $\Phi_{q},\Phi_{q}^{\mathrm{sym}}:[0,1]\to\mathbb{R}$
as 
\begin{align*}
\Phi_{q}(t) & :=t\ln_{q}(t)\quad\textrm{ and \quad}\Phi_{q}^{\mathrm{sym}}(t):=t\ln_{q}(t)+(1-t)\ln_{q}(1-t),
\end{align*}
where $0\ln_{q}(0):=0$, and the function $\ln_{q}:(0,+\infty)\to\mathbb{R}$
is defined as 
\[
\ln_{q}(t):=\begin{cases}
\frac{t^{q-1}-1}{q-1}, & q\neq1\\
\ln(t), & q=1
\end{cases}
\]
and is known as the \emph{$q$-logarithm} introduced by Tsallis \cite{tsallis1994numbers},
but with a slight reparameterization.  In particular, for $q=1$,
we also denote 
\[
h(t):=\Phi_{1}^{\mathrm{sym}}(t)=t\ln t+(1-t)\ln(1-t).
\]
The function $-h$ is known as  the  binary entropy function.

\begin{definition}
The $\Phi_{q}$-stability and the $\Phi_{q}^{\mathrm{sym}}$-stability
are respectively called the\emph{ asymmetric and symmetric $q$-stabilities},
denoted as $\mathbf{Stab}_{q}$ and $\mathbf{Stab}_{q}^{\mathrm{sym}}$. 
\end{definition}
For a Boolean function $f$, by definition, $I_{\Phi_{q}^{\mathrm{sym}}}(f(\mathbf{Y});\mathbf{X})$
reduces to the \emph{Tsallis mutual information} between $f(\mathbf{Y})$
and $\mathbf{X}$. In particular, when $q=1$, $I_{\Phi_{1}^{\mathrm{sym}}}(f(\mathbf{Y});\mathbf{X})$
reduces to the \emph{Shannon mutual information} (denoted as $I(f(\mathbf{Y});\mathbf{X})$
for short). 

There are several fundamental conjectures on the noise stability.
A Boolean function $f_{\mathrm{d}}$ is called a dictator function
if $f_{\mathrm{d}}=1\{x_{k}=1\}$ or $1\{x_{k}=-1\}$ for some $1\le k\le n$.
\begin{conjecture}[Asymmetric $q$-Stability Conjecture]
\label{conj:AsymmetricStability} For $0\le\rho\le1$ and $q\in[1,9]$,
it holds that $\mathbf{Stab}_{q}(f)\le\mathbf{Stab}_{q}(f_{\mathrm{d}})$
for all balanced Boolean functions $f$, where $f_{\mathrm{d}}$ is
a dictator function.
\end{conjecture}
\begin{conjecture}[Symmetric $q$-Stability Conjecture]
\label{conj:SymmetricStability} For $0\le\rho\le1$ and $q\in[1,9]$,
it holds that $\mathbf{Stab}_{q}^{\mathrm{sym}}(f)\le\mathbf{Stab}_{q}^{\mathrm{sym}}(f_{\mathrm{d}})$
for all balanced Boolean functions $f$, where $f_{\mathrm{d}}$ is
a dictator function. 
\end{conjecture}
These two conjectures are unifications of  several existing conjectures \cite{yu2023phi}.
These two conjectures are known as the Mossel--O'Donnell
conjecture \cite{mossel2005coin} when $q\in(2,9]$, known as the Li--M\'edard conjecture \cite{li2020boolean} when $q\in(1,2)$, and known as the Courtade--Kumar conjecture when $q=1$.
In particular, for  $q=1$, the Courtade--Kumar
conjecture can be interpreted as the problem of maximizing the Shannon
mutual information $I(f(\mathbf{Y});\mathbf{X})$ over all Boolean
functions $f$ with a given mean. This conjecture attracts considerable
interest from different fields \cite{anantharam2013on,kindler2015remarks,ordentlich2016improved,samorodnitsky2016entropy,pichler2018dictator,li2020boolean},
and it is regarded as one of the most fundamental conjectures at the
interface of information theory and the analysis of Boolean functions.
Furthermore, it is known that the Courtade--Kumar conjecture is equivalent
to the Li--M\'edard conjecture \cite{barnes2020courtade}.  It is
also interesting to ask a stronger question: What is the interval
$[a,b]\subset[1,\infty)$ for which dictator functions maximize asymmetric
(or symmetric) $q$-stability if and only if $q\in[a,b]$ \cite{yu2023phi}?
It is known that such an interval exists \cite{barnes2020courtade,yu2022common}.

Besides the motivation from learning theory given above, 
other motivations for studying the noise stability problem are as follows.
Maximizing noise stability for Boolean functions with given mean is
equivalent to minimizing the size of a ``soft-boundary'' for subsets
with given size. In fact, the classic edge-isoperimetric inequality
proven by Harper can be recovered from the noise stability problem
(by taking the limit $\rho\to1$) \cite{ledoux1994semigroup,courtade2014boolean,li2020boolean}.
So, the noise stability problem can be seen as a probabilistic generalization
of the classic isoperimetric problem. Furthermore, the noise stability
problem can recover the classic hypercontractivity inequality, which
means that the noise stability problem is a strengthening of the hypercontractivity
inequality \cite{nair2014equivalent,ledoux2014remarks}. The noise
stability problem has also many applications in Boolean function analysis.
By Fourier expansion, the noise stability of a Boolean function is
exactly the generating function of Fourier weights of this Boolean
function, where the correlation coefficient $\rho$ is the indeterminate
\cite{O'Donnell14analysisof}. On the other hand, estimating Fourier
weights is a fundamental question in Boolean function analysis. 

\subsection{Related Works}

The study of the noise stability problem, or more precisely, a two-function
(or two-set) version of the noise stability problem called the \emph{non-interactive
correlation distillation (NICD) problem}, dates back to G\'acs and K\"orner's
and Witsenhausen's seminal papers \cite{gacs1973common,witsenhausen1975sequences}.
By utilizing the tensorization property of the maximal correlation,
Witsenhausen \cite{witsenhausen1975sequences} showed that for $q=2$,
the asymmetric and symmetric $2$-stability w.r.t. $\alpha=1/2$ are
attained by dictator functions. The symmetric $q$-stability problem
with $q\in\left\{ 3,4,5,...\right\} $ was studied by Mossel and O'Donnell
\cite{mossel2005coin}, and  the case $q=3$ was solved by them.
Specifically, they showed that for $q=3$, the symmetric $q$-stability
w.r.t. $\alpha=1/2$ is maximized by dictator functions. This result
was improved by the present author in \cite{yu2023phi} by showing
that the asymmetric and symmetric $q$-stability w.r.t. $\alpha=1/2$
are maximized by dictator functions for $q\in[2,3]$ and $q\in[2,5]$
respectively. As for the other direction, Mossel and O'Donnell observed 
that the asymmetric and symmetric $q$-stability are not maximized
by dictator functions for $q=10$.  

In 2016, Samorodnitsky \cite{samorodnitsky2016entropy} made a significant
breakthrough on the Courtade--Kumar conjecture. Specifically, he
proved the existence of a {\em dimension-independent} threshold
$\rho_{0}$ for which the Courtade--Kumar conjecture holds for all
$\rho\in[0,\rho_{0}]$. However, the value of $\rho_{0}$ was not
explicitly given in his paper, and required to be ``sufficiently
small''. In 2023, the present author \cite{yu2023phi} used a different
method to prove an explicit threshold $\rho_{1}\approx0.46$ for which
the Courtade--Kumar conjecture holds for all $\rho\in[0,\rho_{1}]$.
In addition, the present author also provided a bound for the Courtade--Kumar
conjecture in \cite{yu2023phi} which is expressed in the form of
finite-dimensional program. By evaluating this bound, numerical results
indicate that the threshold $\rho_{1}$ here can be improved to $0.83$.
In addition, Anantharam et al. proposed a strengthened version of  the Courtade--Kumar conjecture
by considering another specific $\Phi$-stability \cite{anantharam2017conjecture},
and Chen and Nair proved bounds for this kind of $\Phi$-stability
by using ideas from isoperimetric inequalities in \cite{chen2024optimality}.
By considering a variant of noise model, Eldan, Mikulincer, and Raghavendra
recently proved a variant version of the Courtade--Kumar conjecture
\cite{eldan2022noise}, whose proof is based on the so-called renormalized
Brownian motion. A weaker version (i.e., the two-function version) of Courtade--Kumar
conjecture was solved by Pichler, Piantanida, and Matz \cite{pichler2018dictator}
by using Fourier analysis. Specifically, they showed that $I(f(\mathbf{Y});g(\mathbf{X}))$
is maximized by a pair of identical dictator functions over all Boolean
functions $(f,g)$, not limited to balanced ones.

The $\Phi$-stability problem in the Gaussian setting with $\Phi$
restricted to be convex and increasing was fully resolved by Borell
\cite{borell1985geometric} in 1985. In particular, he showed that
the $\Phi$-stability is maximized by the indicators of half-spaces
over all measurable Boolean functions $f:\mathbb{R}^{n}\to\{0,1\}$
of the same measure. Such a result is known as \emph{Borell's Isoperimetric
Theorem}. The Gaussian analogues of Conjectures \ref{conj:AsymmetricStability}
and \ref{conj:SymmetricStability} with $q\in[1,\infty)$ are implied
by Borell's Isoperimetric Theorem, which were also proved respectively
by Kindler, O'Donnell, and Witmer \cite{kindler2015remarks} and
by Eldan \cite{eldan2015two} using alternative approaches. 

\subsection{Our Contributions}

\begin{figure}
\centering \includegraphics[scale=0.6]{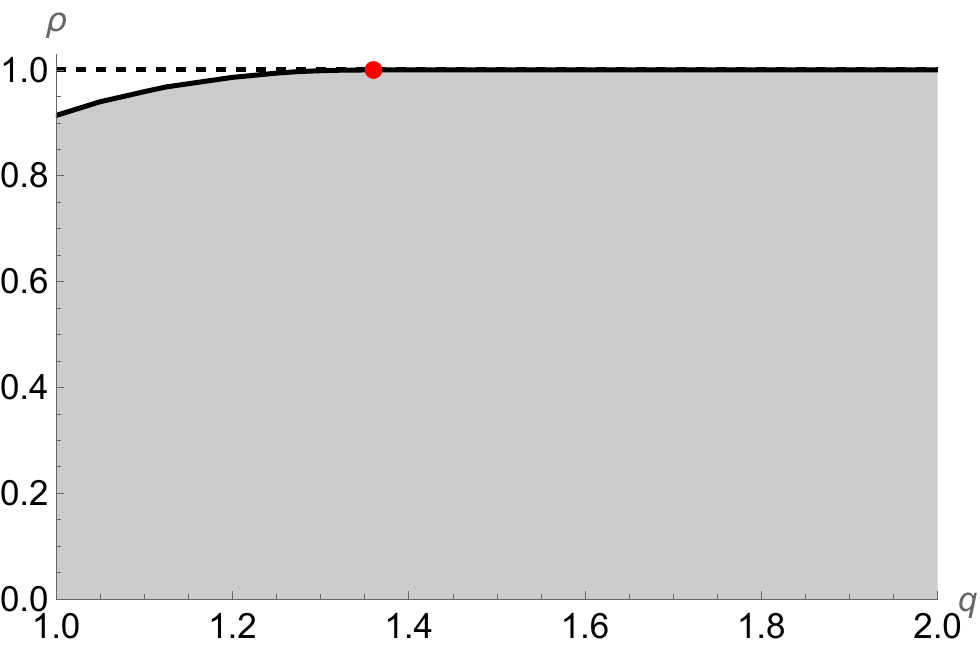}

\caption{\label{fig:rho_q}The region of $(q,\rho)$ (the gray region here)
for which dictator functions maximize the $q$-stability (i.e., for
which the Li--M\'edard conjecture is true).  For $q=1$, the maximum
$\rho$ on the boundary of this region is $\rho=0.914$. For $q\ge1.36$,
the corresponding  $\rho$ on the boundary of this region is $\rho=1$.}
\end{figure}

Rearrangement and majorization techniques are efficient tools in proving
geometric inequalities. Such techniques were also used in the study
of the noise stability problem \cite{borell1985geometric,kindler2015remarks,chen2024optimality}.
In this paper, we investigate the majorization of $T_{\rho}f$ for
a Boolean function $f$, and connect it to the noise stability. We
then use the majorization and hypercontractivity inequalities to
provide upper bounds on the $\Phi$-stability. By these bounds, we
prove that dictator functions locally maximize the $\Phi$-stability
among all balanced Boolean functions. When focusing on the symmetric
$q$-stability, combining this result with our previous bound, we
use computer-assisted methods to prove that dictator functions maximize
the symmetric $q$-stability for $q=1$ and $\rho\in[0,0.914]$ or
for $q\in[1.36,2)$ and all $\rho\in[0,1]$. 
In other words, our main results are the following two theorems. 
\begin{theorem}
\label{thm:CK2-2}The Courtade--Kumar
conjecture is true for all $\rho\in[0,0.914]$.
\end{theorem}
\begin{theorem}
\label{thm:Li-Medard-2}The (symmetrized) Li--M\'edard conjecture is
true for all $q\in[1.36,2)$ and $\rho\in[0,1]$, and also true for all $q\in(1,1.36)$ and $\rho\in[0,0.914]$.
\end{theorem}
For $q\in(1,1.36)$, we numerically
evaluate our bound and find a region of $(q,\rho)$ for which dictator
functions maximize the $q$-stability (i.e., for which the Li--M\'edard
conjecture is true). We plot this region in Fig. \ref{fig:rho_q}.
In addition, it is also interesting to consider the case $q\in(0,1)$.
For this case, we conjecture that dictator functions maximize both
the symmetric and asymmetric $\frac{1}{2}$-stability over all balanced
Boolean functions.

\emph{Note Added:} After the submission of this work, a related study \cite{chen2025differential} by Chen, Gohari, and Nair was posted on arXiv in Feb. 2025. In their work, they reduced the Courtade--Kumar conjecture to a set of conjectured inequalities involving up to $4$ variables.

\subsection{Organization}

This paper is organized as follows. In Section 2, we introduce majorization
theory. In Section 3, we investigate the majorization for noise operators.
In Section 4, we investigate the local optimality of dictator functions
in maximizing the $\Phi$-stability. In Section 5, we focus on the
$q$-stability and investigate the global optimality of dictator functions
in maximizing the $q$-stability. We improve our previous bounds on
the Courtade--Kumar conjecture and the Li--M\'edard conjecture. In
Section 6, we extend our bounds on the $\Phi$-stability to arbitrary
Boolean functions and also to the Gaussian setting.  The proofs of
our results are provided in Appendices A-H. 

\section{Preliminaries}

Since proofs of our results in this paper rely on majorization, we
introduce basic concepts and properties in this field. Here, we adopt the notations
in \cite{marshall2011inequalities}.  Let $(\mathcal{X},\Sigma_{\mathcal{X}},\mu)$
be a measure space and let $f$ be a nonnegative $\mu$-integrable
(real-valued) function defined on $\mathcal{X}$. Define the distribution
function of $f$ with respect to $\mu$ as 
\begin{equation}
\mathfrak{m}_{f}(t):=\mu\{x:f(x)>t\},\;t\ge0.\label{eq:-1}
\end{equation}
Define the decreasing rearrangement of $f$ as 
\begin{equation}
f_{\downarrow}(v):=\sup\{t:\mathfrak{m}_{f}(t)>v\},\;v\ge0.\label{eq:}
\end{equation}
It is easy to see that the distribution function $\mathfrak{m}_{f}$ is nonincreasing and right-continuous on $[0,\infty)$, and finite on $(0,\infty)$. Moreover,
\[
f_{\downarrow}(v)=\lambda\{t:\mathfrak{m}_{f}(t)>v\}=\mathfrak{m}_{\mathfrak{m}_{f}}(v),\;\forall v\ge0,
\]
where $\lambda$ is the Lebesgue measure on $\mathbb{R}_{+}:=[0,\infty)$.
Thus, the decreasing rearrangement of $f$ is the distribution function
(w.r.t. the Lebesgue measure $\lambda$) of the distribution function (w.r.t. the original  measure $\mu$)  of $f$. Furthermore, it should be
noting that the domain of $f_{\downarrow}$ is $[0,\mu(\mathcal{X})]$
which usually differs from the domain $\mathcal{X}$ of $f$. 
\begin{definition}
Let $f$ and $g$ be nonnegative integrable functions on $(\mathcal{X},\Sigma_{\mathcal{X}},\mu)$.
The function $f$ is majorized by $g$, denoted $f\prec g$, if 
\begin{enumerate}
\item $\int_{0}^{t}f_{\downarrow}(v)\d v\le\int_{0}^{t}g_{\downarrow}(v)\d v,\;\forall t\ge0$;
\item $\int f\d\mu=\int g\d\mu$. 
\end{enumerate}
\end{definition}
Since $\int f\d\mu=\int_{0}^{\infty}\mathfrak{m}_{f}(t)\d t=\int_{0}^{\infty}f_{\downarrow}(v)\d v$,
the second condition can be rewritten as $\int_{0}^{\infty}\mathfrak{m}_{f}(t)\d t=\int_{0}^{\infty}\mathfrak{m}_{g}(t)\d t$,
or equivalently, $\int_{0}^{\infty}f_{\downarrow}(v)\d v=\int_{0}^{\infty}g_{\downarrow}(v)\d v$.

The quantity $\int_{0}^{t}f_{\downarrow}(v)\d v$ is in fact closely
related to the optimal tradeoff between two types of errors in hypothesis
testing. Note that 
$A\in\Sigma_{\mathcal{X}}\mapsto\int_{A}f\d\mu$ defines a measure,
which is denoted by $\mu_{f}$. For a nonnegative
function $f$, define the concentration spectrum 
\begin{align}
\mathfrak{D}_{f}(t) & :=\sup_{A:\mu(A)\le t}\mu_{f}(A).\label{eq:concentration-2}
\end{align}
and its smooth version 
\begin{align}
\mathfrak{C}_{f}(t) & :=\sup_{\phi:\mu(\phi)\le t}\mu_{f}(\phi)=\sup_{\phi:\int\phi\d\mu\le t}\int\phi f\d\mu,\label{eq:concentration}
\end{align}
where the suprema in \eqref{eq:concentration} are over all measurable
$\phi:\mathcal{X}\to[0,1]$.  By the Neyman--Pearson lemma, this
concentration spectrum completely characterizes the optimal tradeoff
between two types of errors in the hypothesis testing between two
probability measures $\mu_{f}$ and $\mu$. It is known that (see,
e.g., \cite{yu2022theentro}) 
\begin{align}
\mathfrak{C}_{f}(t) & =\conc\mathfrak{D}_{f}(t) =\int_{0}^{t}f_{\downarrow}(v)\d v, \label{eq:CD}
\end{align}
where $\conc g$ denotes the upper concave envelope of a function
$g$. 

Another related quantity is called the $E_{\gamma}$-divergence \cite{polyanskiy2010channel}.
For two probability measures $P,Q$, 
\begin{align*}
E_{\gamma}(P\|Q) & :=\max_{A}\{P(A)-\gamma Q(A)\}\\
 & =P\{\d P/\d Q>\gamma\}-\gamma Q\{\d P/\d Q>\gamma\}\\
 & =\int[\d P/\d Q-\gamma]^{+}\mathrm{d}Q,
\end{align*}
where $[x]^{+}:=\max\{0,x\}$.
\begin{proposition}
\label{prop:majorization-1} \cite[ Proposition H.1.a]{marshall2011inequalities}
Let $f$ and $g$ be nonnegative integrable functions on $(\mathcal{X},\Sigma_{\mathcal{X}},\mu)$.
Then, the following conditions are equivalent. 
\begin{enumerate}
\item $f\prec g$;
\item $\mathfrak{C}_{f}(t)\le\mathfrak{C}_{g}(t),\;\forall t\ge0$;
\item $\int\Phi(f)\mathrm{d}\mu\le\int\Phi(g)\mathrm{d}\mu$ for all continuous
convex functions $\Phi:\mathbb{R}_{+}\to\mathbb{R}$ for which the
integrals exist;
\item $\int_{t}^{\infty}\mathfrak{m}_{f}(s)\d s\le\int_{t}^{\infty}\mathfrak{m}_{g}(s)\d s,\;\forall t\ge0$; 
\item $E_{\gamma}(\mu_{f}\|\mu)\le E_{\gamma}(\mu_{g}\|\mu)$, i.e., $\int[f-\gamma]^{+}\mathrm{d}\mu\le\int[g-\gamma]^{+}\mathrm{d}\mu,\;\forall\gamma\ge0$. 
\end{enumerate}
\end{proposition}

The concept of majorization can be extended to a pair of functions
defined on two totally different measure spaces. 
\begin{definition}
\label{def:majorization}Let $f$ and $g$ be nonnegative integrable
functions respectively defined on $(\mathcal{X},\Sigma_{\mathcal{X}},\mu)$
and $(\mathcal{Y},\Sigma_{\mathcal{Y}},\nu)$. The function $f$ is
majorized by $g$, denoted $f\prec g$, if 
\begin{enumerate}
\item $\int_{0}^{t}f_{\downarrow}(v)\d v\le\int_{0}^{t}g_{\downarrow}(v)\d v,\;\forall t\ge0$,
where $f_{\downarrow}$ is defined in \eqref{eq:} with respect to
$\mu$ and $g_{\downarrow}$ is defined similarly but with respect
to $\nu$;
\item $\int_{\mathcal{X}}f\d\mu=\int_{\mathcal{Y}}g\d\nu$. 
\end{enumerate}
\end{definition}
Although in this general definition, $f$ and $g$ are defined on
two different measure spaces,  their decreasing rearrangements $f_{\downarrow}$
and $g_{\downarrow}$ are defined on the same one, i.e., $(\mathbb{R}_{+},\Sigma_{\mathbb{R}_{+}},\lambda)$. Clearly, the decreasing
rearrangements of $f_{\downarrow}$ and $g_{\downarrow}$ with respect
to the Lebesgue measure $\lambda$ are $f_{\downarrow}$ and $g_{\downarrow}$
themselves. Hence, according to Definition \ref{def:majorization}, we
arrive at the following observation. 
\begin{proposition}
\label{prop:majorization} Let $f$ and $g$ be nonnegative integrable
functions respectively defined on $(\mathcal{X},\Sigma_{\mathcal{X}},\mu)$
and $(\mathcal{Y},\Sigma_{\mathcal{Y}},\nu)$. Then, $f\prec g$ if
and only if $f_{\downarrow}\prec g_{\downarrow}$. 
\end{proposition}
In this paper, we focus on the integral $\int_{\mathcal{X}}\Phi(f)\mathrm{d}\mu$
for a given measurable function $\Phi$. The following well known result shows
that this integral can be expressed as the integral of $\Phi(f_{\downarrow})$.
\begin{proposition}
\label{prop:integral}Let $f$ be a nonnegative  measurable  function
defined on $(\mathcal{X},\Sigma_{\mathcal{X}},\mu)$.
 such that its distribution function is finite on $(0,\infty)$. Then, 
\begin{equation}
\int_{\mathcal{X}}\Phi(f)\mathrm{d}\mu=\int_{\mathbb{R}_{+}}\Phi(f_{\downarrow})\mathrm{d}\lambda, \label{eq:integration}
\end{equation}
for any measurable (not necessarily convex) function $\Phi:\mathbb{R}_{+}\to\mathbb{R}$,
for which at least one of (equivalently, both of) the integrals exists.
\end{proposition}
\begin{proof}
By definition, the distribution function of $f$ w.r.t. $\mu$ is identical to the distribution function of $f_{\downarrow}$ w.r.t. $\lambda$. Moreover,  a  distribution function that takes finite values  on $(0,\infty)$  determines a unique Borel measure on $(0,\infty)$ (see e.g., \cite[Theorem 1.16]{folland_real_1999}), and the measure of $\{0\}$ is also determined by the distribution function. So, the measures $\mu \circ f^{-1}$ and $\lambda \circ f_{\downarrow}^{-1}$ are identical. 
On the other hand, for any measurable function $\Phi:\mathbb{R}_{+}\to\mathbb{R}$,
\begin{align*}
\int_{\mathcal{X}}\Phi(f)\mathrm{d}\mu & =\int_{\mathbb{R}_{+}}\Phi\mathrm{d}(\mu \circ f^{-1})=\int_{\mathbb{R}_{+}}\Phi\mathrm{d}(\lambda \circ f_{\downarrow}^{-1})=\int_{\mathbb{R}_{+}}\Phi(f_{\downarrow})\mathrm{d}\lambda.
\end{align*}
\end{proof}
By Propositions \ref{prop:majorization-1}, \ref{prop:majorization},
and \ref{prop:integral}, in order to compare two integrals $\int_{\mathcal{X}}\Phi(f)\mathrm{d}\mu$
and $\int_{\mathcal{Y}}\Phi(g)\mathrm{d}\nu$ for a continuous convex
function $\Phi$, it suffices to restrict our attention to the majorization
relation between $f_{\downarrow}$ and $g_{\downarrow}$. 

\begin{proposition}
\label{prop:majorization-1-1} Let $f$ and $g$ be nonnegative integrable
functions respectively defined on $(\mathcal{X},\Sigma_{\mathcal{X}},\mu)$
and $(\mathcal{Y},\Sigma_{\mathcal{Y}},\nu)$. Then, the following
conditions are equivalent. 
\begin{enumerate}
\item $f\prec g$;
\item $f_{\downarrow}\prec g_{\downarrow}$; 
\item $\mathfrak{C}_{f}(t)\le\mathfrak{C}_{g}(t),\;\forall t\ge0$;
\item $\int_{\mathcal{X}}\Phi(f)\mathrm{d}\mu\le\int_{\mathcal{Y}}\Phi(g)\mathrm{d}\nu$
for all continuous convex functions $\Phi:\mathbb{R}_{+}\to\mathbb{R}$
for which the integrals exist;
\item $\int_{t}^{\infty}\mathfrak{m}_{f}(s)\d s\le\int_{t}^{\infty}\mathfrak{m}_{g}(s)\d s,\;\forall t\ge0$; 
\item $E_{\gamma}(\mu_{f}\|\mu)\le E_{\gamma}(\nu_{g}\|\nu)$, i.e., $\int_{\mathcal{X}}[f-\gamma]^{+}\mathrm{d}\mu\le\int_{\mathcal{Y}}[g-\gamma]^{+}\mathrm{d}\nu,\;\forall\gamma\ge0$. 
\end{enumerate}
\end{proposition}
The concept of majorization can be extended to the relation between
a set of functions and a function. 
\begin{definition}
Let $\mathcal{F}$ be a set of nonnegative functions on $(\mathcal{X},\Sigma_{\mathcal{X}},\mu)$
and $g$ be a nonnegative function on $(\mathcal{Y},\Sigma_{\mathcal{Y}},\nu)$.
We say $\mathcal{F}$ is \textit{majorized} by $g$ (and write $\mathcal{F}\prec g$)
if $f\prec g$ for all $f\in\mathcal{F}$. 
\end{definition}
\begin{proposition}
\label{prop:maj_to_convex-1-1} Let $\mathcal{F}$ be a set of nonnegative
functions on $(\mathcal{X},\Sigma_{\mathcal{X}},\mu)$ and $g$ be
a nonnegative function on $(\mathcal{Y},\Sigma_{\mathcal{Y}},\nu)$
such that $\mathcal{F}\prec g$. Then, for any continuous convex function
$\Phi:\mathbb{R}_{+}\to\mathbb{R}$, 
\[
\int_{\mathcal{X}}\Phi(f)\mathrm{d}\mu\le\int_{\mathcal{Y}}\Phi(g)\mathrm{d}\nu,\;\forall f\in\mathcal{F}.
\]
 
\end{proposition}

\section{Majorization for Noise Operators }

We now focus on the hypercube $\{0,1\}^{n}$ and investigate the majorization for $T_{\rho}f$ with Boolean functions $f$. Let $\mu$ be
the uniform distribution on $\{0,1\}^{n}$. We say $f$ is a Boolean
function if it is a function $f:\{0,1\}^{n}\to\{0,1\}$ for some $n$.
Denote $\mu(f):=\mathbb{E}_{\mu}[f]$.
\begin{definition}
Given $\alpha,\rho \in[0,1]$, we denote $\mathcal{F}_{\alpha,\rho}$ (or
briefly, $\mathcal{F}_{\alpha}$) as the set of $T_{\rho}f$ for
Boolean functions $f$ with mean $\mu(f)=\alpha$. 
\end{definition}

The concentration spectrum of $T_{\rho}f$ is 
\begin{align}
\mathfrak{D}_{T_{\rho}f}(\beta) & =\max_{\textrm{Bool }\phi:\,\mu(\phi)\le \beta}\int\phi T_{\rho}f\d\mu,\label{eq:concentration-2-1}
\end{align}
and its smooth version is 
\begin{align}
\mathfrak{C}_{T_{\rho}f}(\beta) & =\max_{\phi:\mu(\phi)\le \beta}\int\phi T_{\rho}f\d\mu,\label{eq:concentration-1}
\end{align}
where the maximization in \eqref{eq:concentration-1} is over all  
$\phi:\{0,1\}^{n}\to[0,1]$ such that $\mu(\phi) \le \beta$. By the fundamental theorem of linear programming, these two functions coincide when $\beta\in2^{-n}[2^{n}]$.

A  concept  related to these two quantities  is 
the maximal noise stability. 
For $\alpha,\beta\in2^{-n}[2^{n}]$,
define the \textit{maximal
noise stability} at $(\alpha,\beta)$  as 
\[
\mathbf{S}_{n}(\alpha,\beta):=\max_{\textrm{Bool }f,\phi:\,\mu(f) = \alpha,\mu(\phi)= \beta}\int\phi T_{\rho}f\d\mu,
\]
where the subscript  $n$   refers to the dimension.
In the following, we show that $\mathfrak{C}_{T_{\rho}f}(\beta)\le \lim_{m\to \infty} \mathbf{S}_{m}(\alpha,\beta_m)$ holds  for any $\beta\in [0,1]$ and any sequence $\beta_m$   such that  $\beta_m\in2^{-m}[2^{m}]$ and $\beta_m\to \beta$, 
 where $\alpha$ is the mean of $f$. 
Combining this  inequality  with Proposition \ref{prop:majorization-1-1} implies the following  relation between noise stability and  majorization for the noise operator.

\begin{proposition}[Relation Between Majorization and Noise Stability]
\label{prop:-Let-} Let $\Theta$ be an upper bound on $\mathbf{S}_{n},\forall n\ge 1$, which satisfies that given $\alpha$, $\Theta(\alpha,\beta)$ is absolutely continuous
and nondecreasing in $\beta$, and $\Theta(\alpha,0)=0$. Then, for
any $\rho,\alpha\in[0,1]$, it holds that $\mathcal{F}_{\alpha}\prec\theta_{\alpha}$,
where for $\beta\in[0,1]$, 
\begin{equation}
\theta_{\alpha}(\beta):=\frac{\partial\Theta(\alpha,\beta)}{\partial\beta}.\label{eq:theta}
\end{equation}
\end{proposition}
 
\begin{remark}
Although this proposition is stated for the Bonami-Beckner operator
on the hypercube, it in fact can be straightforwardly generalized
to an arbitrary semigroup operator on an arbitrary space. 
\end{remark}
\begin{proof}
Since  $\beta \in [0,1] \mapsto \Theta(\alpha,\beta)$ is nondecreasing, it is differentiable almost everywhere. By the fundamental theorem of Lebesgue integral calculus,    $\Theta(\alpha,\beta)=\int_{0}^{\beta}\frac{\partial\Theta(\alpha,t)}{\partial t}\d t$.

Let $m\ge n$. A Boolean function $f$ on the (discrete) $n$-cube can be viewed as a Boolean function 
on the $m$-cube  which only depends on the first $n$ coordinates. 
This yields that      
$T_{\rho} f$ (with $T_{\rho}=T_{\rho}^{(n)}$), if viewed as the function $T_{\rho}^{(m)} f$ on the $m$-cube, also only depends on the first $n$ coordinates. From this view, we can rewrite 
\begin{align}
\mathfrak{C}_{T_{\rho}f}(\beta) & =\max_{\phi:\mu^{(m)}(\phi)\le \beta}\int\phi T_{\rho}f\d\mu^{(m)}, \label{eq:concentration-3}
\end{align}
where the maximization in \eqref{eq:concentration-3} is over all  
$\phi:\{0,1\}^{m}\to[0,1]$ such that $\mu^{(m)}(\phi) \le \beta$. 
This formula follows since, on one hand, if we   view a feasible  solution  $\phi$ to the maximization in \eqref{eq:concentration-1}   as the corresponding function on $\{0,1\}^{m}$ that only depends  on the first $n$ coordinates, then  $\phi$ is  also feasible to the maximization in \eqref{eq:concentration-3}, which implies that  the LHS in \eqref{eq:concentration-3} is smaller than or equal to the RHS; on the other hand, for a feasible solution $\phi$ to 
the maximization in \eqref{eq:concentration-3}, if we convert it to the function $(x_1,...,x_n)\mapsto \mathbb{E}[\phi (x_1,...,x_n,X_{n+1},...,X_m)]$ on the $n$-cube, then the resulting function is feasible to  the maximization in \eqref{eq:concentration-1} and also admits the same objective value as the one in \eqref{eq:concentration-3}, which implies that  the LHS in \eqref{eq:concentration-3} is larger than or equal to the RHS.

Obviously, the optimization in  \eqref{eq:concentration-3} is a linear program (viewing $\phi$ as a vector $(\phi(\mathbf{x}):\mathbf{x}\in\{0,1\}^{m})$), and by the fundamental theorem of linear programming, the optimal value is attained at some vertex (or extreme point) of the feasible region. For $\beta\in2^{-m}[2^{m}]$, a vertex (i.e., an extreme point) of the feasible region must satisfy that $\phi(\mathbf{x})=0$ or $1$ for all $\mathbf{x}\in\{0,1\}^{m}$, which corresponds to a Boolean function on $\{0,1\}^{m}$. 
That is, for such $\beta$, $\mathfrak{C}_{T_{\rho}f}(\beta)$  is attained by a Boolean function $\phi$ with mean $\beta$, and hence, 
$\mathfrak{C}_{T_{\rho}f}(\beta)\le \mathbf{S}_{m}(\alpha,\beta)$, where $\alpha$ is the mean of $f$.  By  assumption, 
for such $\beta$, 
$\mathfrak{C}_{T_{\rho}f}(\beta)\le \Theta(\alpha,\beta)$. 
Since $\bigcup_{m \ge n}2^{-m}[2^{m}]$ is dense in $[0,1]$ and   
both these two functions are continuous in $\beta$,   this inequality can be extended to all 
 $\beta \in [0,1]$. 
Furthermore,  by the definition of smooth concentration spectrum, $\mathfrak{C}_{\theta_{\alpha}}(\beta)\ge\int_{0}^{\beta}\frac{\partial\Theta(\alpha,t)}{\partial t}\d t=\Theta(\alpha,\beta).$
Combining the two inequalities above yields  $\mathfrak{C}_{T_{\rho}f}(\beta)\le\mathfrak{C}_{\theta_{\alpha}}(\beta)$
for all   $\beta \in [0,1]$. Combined with Proposition \ref{prop:majorization-1-1},
these imply $\mathcal{F}_{\alpha}\prec\theta_{\alpha}.$

\end{proof}

A well-known upper bound on the maximal noise stability is given by
the small-set expansion theorem \cite[p.280]{O'Donnell14analysisof},
which states that 
\begin{align}
\mathbf{S}_{n}(\alpha,\beta) & \le \Theta_{0}(\alpha,\beta):= \begin{cases}
\exp\left(-\frac{s^{2}+t^{2}-2\rho st}{2(1-\rho^{2})}\right), &  \rho s\le t\le s/\rho, \\
\exp(-\frac{s^{2}}{2}), & 0\le t\le\rho s,\\
\exp(-\frac{t^{2}}{2}), & t\ge s/\rho,
\end{cases}\label{eq:Sn}
\end{align}
where $\alpha=\exp(-\frac{s^{2}}{2})$ and $\beta=\exp(-\frac{t^{2}}{2})$
with $s,t\ge0$.

\begin{figure}
\centering \includegraphics[scale=0.5]{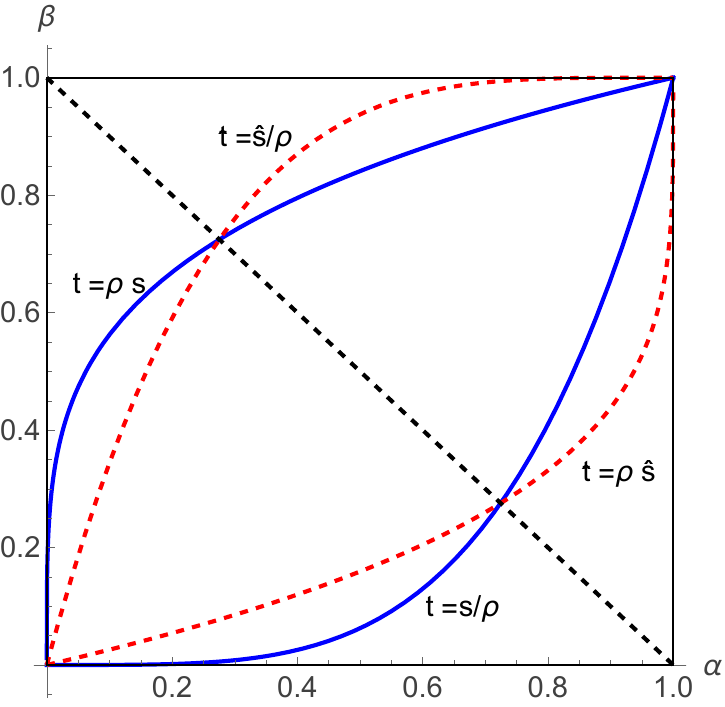}
 \includegraphics[height=0.4 \textwidth , width=0.5 \textwidth]{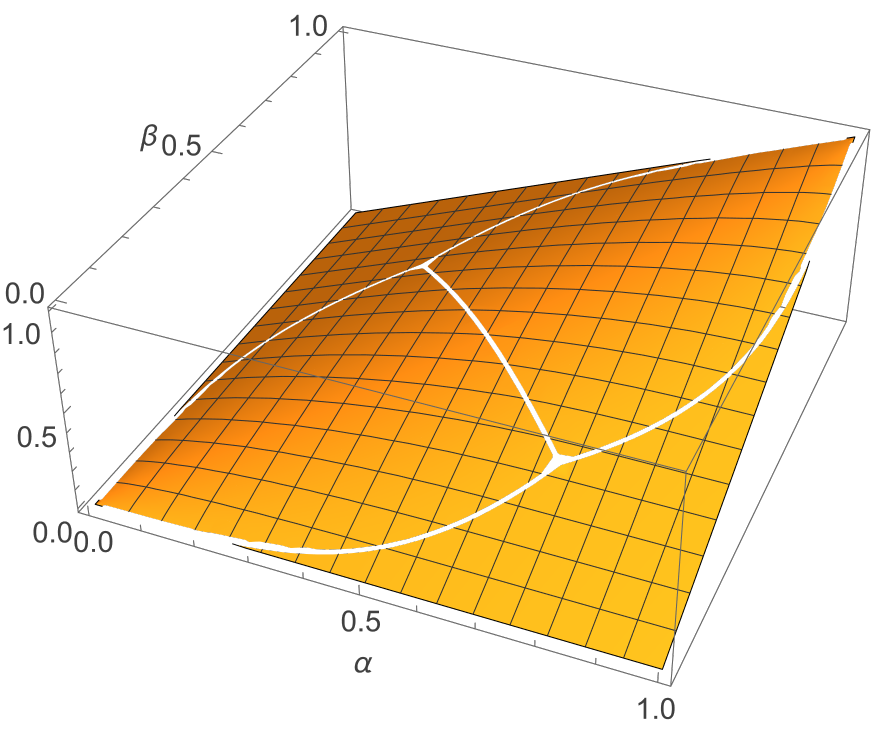}
\caption{\label{fig:alphabeta}Illustration of the  regions  $A_i,B_i, i=0,1,2$ in the left subfigure and the function $\Theta$ given in 
\eqref{eq:Sn2} in the right subfigure.}
\end{figure}

Combining  \eqref{eq:Sn} with the fact that 
\[
\int\phi T_{\rho}f\d\mu=\mu(\phi)+\mu(f)-1+\int(1-\phi)T_{\rho}(1-f)\d\mu,
\]
we obtain  that
\begin{align*}
\mathbf{S}_{n}(\alpha,\beta) \le\Theta(\alpha,\beta)  
& := \min \{ \Theta_{0}(\alpha,\beta), \alpha +\beta -1 + \Theta_{0}(1-\alpha,1-\beta)\}.
\end{align*}
The upper bound $\Theta$  is simplified in the following proposition.
Denote $s=\sqrt{-2\ln\alpha}$, $t=\sqrt{-2\ln\beta}$,
$\hat{s}=\sqrt{-2\ln(1-\alpha)}$, and $\hat{t}=\sqrt{-2\ln(1-\beta)}$, i.e.,  $\alpha=\exp(-\frac{s^{2}}{2})=1-\exp(-\frac{\hat{s}^{2}}{2})$
and $\beta=\exp(-\frac{t^{2}}{2})=1-\exp(-\frac{\hat{t}^{2}}{2})$.
Denote
\begin{align*}
A_0 &= \{ 0\le t\le\rho s\}, A_1=\{\rho s\le t\le s/\rho\}, A_2=\{t\ge s/\rho\}, \\
B_0 &= \{ 0\le \hat{t}\le\rho \hat{s}\}, B_1=\{\rho \hat{s}\le \hat{t}\le \hat{s}/\rho\}, B_2=\{\hat{t}\ge \hat{s}/\rho\}.
\end{align*}
\begin{proposition}\label{prop:Theta}
For $\alpha,\beta\in[0,1]$, it holds that 
\[
\Theta_{0}(\alpha,\beta)\le\alpha+\beta-1+\Theta_{0}(1-\alpha,1-\beta)
\]
if $\alpha+\beta\le1$; the inequality is reversed if $\alpha+\beta>1$. As a consequence, 
\begin{align}
\Theta(\alpha,\beta)  & = 
\begin{cases}
\Theta_{0}(\alpha,\beta), &  \alpha+\beta\le 1 \\  
\alpha +\beta -1 + \Theta_{0}(1-\alpha,1-\beta),   &  \alpha+\beta > 1
\end{cases} \\
& =  \begin{cases}
\exp(-\frac{t^{2}}{2}), &   A_2 \cap   B_0,\\
\exp\left(-\frac{s^{2}+t^{2}-2\rho st}{2(1-\rho^{2})}\right), &  A_1  \cap \{\alpha+\beta\le 1\},  \\ 
\exp(-\frac{s^{2}}{2})+\exp(-\frac{t^{2}}{2})-1  +\exp\left(-\frac{\hat{s}^{2}+\hat{t}^{2}-2\rho\hat{s}\hat{t}}{2(1-\rho^{2})}\right), &  B_1 \cap \{\alpha+\beta> 1\}, \\  
\exp(-\frac{s^{2}}{2}), &  A_0 \cap   B_2.
\end{cases}  \label{eq:Sn2}
\end{align}
\end{proposition}
The proof of this proposition is given in Appendix \ref{sec:Proof-of-Proposition-Theta}. 
These regions  $A_i,B_i, i=0,1,2$   and the function $\Theta$ are  plotted in Fig. \ref{fig:alphabeta}.  
Here $\alpha +\beta \le 1$ is equivalent to $t \ge  \hat{s}$ (or $s \ge  \hat{t}$). Moreover, $\Theta$ is continuous and piecewise continuously differentiable, and hence it is absolutely continuous. 

Combining  Proposition \ref{prop:-Let-} and the explicit bound on noise stability given  in \eqref{eq:Sn2} yields an explicit function that majorizes the noise operator. 

\begin{proposition}[Majorization for Noise Operator]\label{prop:majorizationNO}
For any $\rho,\alpha\in[0,1]$, it holds that $\mathcal{F}_{\alpha}\prec\theta_{\alpha}$,
where for $\beta\in[0,1]$, 

\begin{align}
\theta_{\alpha}(\beta) & :=\frac{\partial\Theta(\alpha,\beta)}{\partial\beta} =\begin{cases}
1, & A_2 \cap  B_0,\\
\displaystyle \frac{t-\rho s}{(1-\rho^{2})t}\exp\left(-\frac{(s-\rho t)^{2}}{2(1-\rho^{2})}\right), &  A_1  \cap \{\alpha+\beta\le 1\},  \\ 
\displaystyle 1-\frac{\hat{t}-\rho \hat{s}}{(1-\rho^{2})\hat{t}}\exp\left(-\frac{(\hat{s}-\rho\hat{t})^{2}}{2(1-\rho^{2})}\right), & B_1 \cap \{\alpha+\beta> 1\}, \\
0, &  A_0 \cap B_2.
\end{cases}\label{eq:theta-1}
\end{align}
with $\alpha=\exp(-\frac{s^{2}}{2})=1-\exp(-\frac{\hat{s}^{2}}{2})$
and $\beta=\exp(-\frac{t^{2}}{2})=1-\exp(-\frac{\hat{t}^{2}}{2})$.
\end{proposition}

We have investigated the majorization of  $T_\rho f$ for Boolean functions $f$ in the above. We next generalize it to the case with $f$ replaced by a mixture of Boolean functions. 
For Boolean functions $f_{i},i\in[k]$ and a vector $\vec{\lambda}:=(\lambda_{1},\lambda_{2},...,\lambda_{k})$
with nonnegative real numbers $\lambda_{i}$, define $f_{\vec{\lambda}}:=\sum_{i=1}^{k}\lambda_{i}f_{i}$.
We now investigate the majorization for $T_{\rho}f_{\vec{\lambda}}$.

\begin{definition}
Given $\rho,\vec{\lambda}:=(\lambda_{1},\lambda_{2},...,\lambda_{k}),\vec{\alpha}:=(\alpha_{1},\alpha_{2},...,\alpha_{k})$,
we denote $\mathcal{F}_{\vec{\lambda},\vec{\alpha}}$ as the set of
$T_{\rho}f_{\vec{\lambda}}$ for Boolean functions $f_{i},i\in[k]$
with mean $\mu(f_{i})=\alpha_{i},\forall i\in[k]$. 
\end{definition}
Following the proofs  of Propositions  \ref{prop:-Let-} and \ref{prop:majorizationNO} (specifically, showing that  $\mathfrak{C}_{T_{\rho}f_{\vec{\lambda}}}(\beta)\le \sum_{i=1}^{k}\lambda_{i}\mathbf{S}_{n}(\alpha_i,\beta)\le \sum_{i=1}^{k}\lambda_{i}\Theta(\alpha_i,\beta)$), we obtain the following proposition. 
\begin{proposition}
Given any $\rho\in[0,1],\vec{\lambda}\in[0,1]^{k},\vec{\alpha}\in[0,1]^{k}$,
it holds that $\mathcal{F}_{\vec{\lambda},\vec{\alpha}}\prec\theta_{\vec{\lambda},\vec{\alpha}}$,
where
\[
\theta_{\vec{\lambda},\vec{\alpha}}(\beta)=\sum_{i=1}^{k}\lambda_{i}\theta_{\alpha_{i}}(\beta)
\]
with $\theta_{\alpha}$ given in \eqref{eq:theta-1}.
\end{proposition}
As a consequence of this proposition, a bound on the $\Phi$-stability for the mixture function $f_{\vec{\lambda}}$ can be obtained. 
\begin{corollary}[Bound on $\Phi$-Stability]
\label{cor:generalbound-1} Let $\Phi:[0,1]\to\mathbb{R}$ be a continuous
convex function. For any Boolean functions $f_{i},i\in[k]$ with
mean $\mu(f_{i})=\alpha_{i}$, it holds that 
\[
\mathbb{E}_{\mu}[\Phi(T_{\rho}f_{\vec{\lambda}})]\le\int_{0}^{1}\Phi(\theta_{\vec{\lambda},\vec{\alpha}}(\beta))\d\beta.
\]
\end{corollary}

\section{Local Optimality of Dictatorships}

\subsection{$\Phi$-Stability}

We next apply the majorization results derived above to investigate the local optimality
of dictator functions in the noise stability problem. 
For ease of presentation, we rewrite the hypercube as $\{\pm1\}^{n}$, where $0,1$ in the original hypercube $\{0,1\}^{n}$ are respectively  mapped to $1,-1$. 
We use $\mathbf{x}=x_{[n]}=(x_{1},...,x_{n})$
to denote a length-$n$ vector. We denote the Bonami-Beckner operator
by $T_{\rho}=T_{\rho}^{(n)}$, and the uniform measure on $\{\pm1\}^{n}$
by $\mu=\mu^{(n)}$. For a Boolean function $f$, denote $A=f^{-1}(1)$,
i.e., the support of $f$. Denote the subcubes $C_{i}:=\{\mathbf{x}:x_{i}=1\},i\in[n]$.
The distance between $A$ and $C_{i}$  is 
\[
d_{i}(f):=\mu(A\Delta C_{i}),
\]
where $\Delta$ is the symmetric difference. It is easy to verify
that 
\[
d_{i}(f)=\|f-1_{C_{i}}\|_{1},
\]
where $\|\phi\|_{1}:=\mathbb{E}_{\mu}[|\phi|]$ is the $L_1$-norm. Hence, $d_{i}(f)$
represents the distance between $f$ and the dictator function $1_{C_{i}}$.

Denote the  Fourier coefficients of $f$  as 
$$\hat{f}_{S}:=\mathbb{E}[f(\mathbf{X})\chi_{S}(\mathbf{X})],S\subseteq[n]$$
where  $\chi_{S}(\mathbf{X})=\Pi_{i\in S}X_{i},S\subseteq[n]$ form an orthogonal Fourier basis (on $\{-1,1\}^n$).
Then, in terms of Fourier analysis, $d_{i}(f)$ can be also expressed as
\[
d_{i}(f)=\frac{1}{2}-\hat{f}_{\{i\}}.
\]
Denote 
\[
\tilde{d}_{i}(f):=\min\{d_{i}(f),1-d_{i}(f)\},
\]
which represents the distance in the $L_1$ sense   between $f$ and the set of dictator functions  $\{1_{C_{i}}, 1_{C_{i}^{c}}\}$. Similarly, the quantity
\begin{equation}  
\tilde{d}(f):=\min_{i\in[n]}\tilde{d}_{i}(f) \label{eq:df}
\end{equation}
is the (minimum) distance between the balanced Boolean function $f$ and the set of dictator functions.
We are ready to state our first main result, which provides a bound on the $\Phi$-stability in terms of the distance $\tilde{d}_{i}(f)$ or $\tilde{d}(f)$.

\begin{theorem}[Bound on $\Phi$-Stability]
\label{thm:generalbound-2} Let $\Phi:[0,1]\to\mathbb{R}$ be a continuous
convex function. For any balanced Boolean function $f$, it holds
that 
\begin{equation}
\mathbb{E}_{\mu}[\Phi(T_{\rho}f)]\le\min_{i\in[n]}\Gamma(\tilde{d}_{i}(f))\le \Gamma(\tilde{d}(f)),\label{eq:-2}
\end{equation}
where for $\epsilon\in[0,1]$, 
\[
\Gamma(\epsilon):=\frac{1}{2}\int_{0}^{1}\left(\Phi(\frac{1+\rho}{2}\theta_{1-\epsilon}(\beta)+\frac{1-\rho}{2}\theta_{\epsilon}(\beta))+\Phi(\frac{1-\rho}{2}\theta_{1-\epsilon}(\beta)+\frac{1+\rho}{2}\theta_{\epsilon}(\beta))\right)\d\beta,
\]
with $\theta_{\epsilon}$ and $\theta_{1-\epsilon}$ given in \eqref{eq:theta}. 
\end{theorem}
\begin{proof}
Firstly, note that for any $i\in[n]$, any balanced Boolean function
$f$ can be written as 
\begin{equation}
f(\mathbf{x})=1\{x_{i}=1\}f_{+}(x_{\backslash i})+1\{x_{i}=-1\}f_{-}(x_{\backslash i}),\label{eq:f}
\end{equation}
where $x_{\backslash i}:=(x_{1},...,x_{i-1},x_{i+1},...,x_{n})\in \{\pm 1\}^{n-1}$,
and $f_{+}(x_{\backslash i}):=f(x_{i\mapsto1})$, $f_{-}(x_{\backslash i}):=f(x_{i\mapsto-1})$  with $x_{i\mapsto a}=(x_{1},...,x_{i-1},a,x_{i+1},...,x_{n}),a\in \{\pm1\}$ are Boolean functions defined on $\{\pm 1\}^{n-1}$.
So, $\mu^{(n-1)}(f_{+})=1-d_{i}(f)$, $\mu^{(n-1)}(f_{-})=d_{i}(f)$, and we can write
\begin{equation}
T_{\rho}f(\mathbf{x})=\begin{cases}
T_{\rho}^{(n-1)}g_{+}(x_{\backslash i}), & x_{i}=1\\
T_{\rho}^{(n-1)}g_{-}(x_{\backslash i}), & x_{i}=-1
\end{cases},\label{eq:Tf}
\end{equation}
where 
\[
g_{+}:=\frac{1+\rho}{2}f_{+}+\frac{1-\rho}{2}f_{-},\quad g_{-}:=\frac{1-\rho}{2}f_{+}+\frac{1+\rho}{2}f_{-}.
\]
The $\Phi$-stability of $f$ is then 
\begin{align}
\mathbb{E}_{\mu}[\Phi(T_{\rho}f)] & =\frac{1}{2}\mathbb{E}\Bigl[\Phi(T_{\rho}^{(n-1)}g_{+}(X_{\backslash i}))+\Phi(T_{\rho}^{(n-1)}g_{-}(X_{\backslash i}))\Bigr]. \label{eq:EPhi}
\end{align}
By Corollary \ref{cor:generalbound-1}, we obtain 
\begin{equation}
\mathbb{E}_{\mu}[\Phi(T_{\rho}f)]\le\min_{i\in[n]}\Gamma(d_{i}(f)).\label{eq:-2-4}
\end{equation}
 Since $\Gamma(\epsilon)=\Gamma(1-\epsilon)$, \eqref{eq:-2-4}
implies \eqref{eq:-2}. 
\end{proof}

When we take $\epsilon\downarrow0$, Theorem \ref{thm:generalbound-2}
implies the following asymptotic bound. The proof is provided in Appendix
\ref{sec:Proof-of-Theorem}.

\begin{theorem}[Asymptotic Bound]
\label{thm:asymptotics} Let $\Phi:[0,1]\to\mathbb{R}$ be continuous
convex on $[0,1]$ and twice continuously differentiable on $(0,1)$.
Let $\rho\in(0,1)$. Then, for any $c\in(\rho,1/\rho)$ it holds
that as $\epsilon\downarrow0$,
\[
\mathbf{Stab}_{\Phi}[f]\le A_{0}-A_{1}\epsilon+A_{2}(1+o(1))\epsilon^{\frac{2}{1+\rho^{2}}}\sqrt{2\ln\frac{1}{\epsilon}}+O(\epsilon^{c^{2}}),
\]
for any balanced Boolean function $f$ such that $\tilde{d}_{i}(f)=\epsilon$
for some $i\in[n]$, where 
\begin{align*}
A_{0} & =\frac{1}{2}\left(\Phi(\frac{1+\rho}{2})+\Phi(\frac{1-\rho}{2})\right),\\
A_{1} & =\frac{\rho}{2}\left(\Phi'(\frac{1+\rho}{2})-\Phi'(\frac{1-\rho}{2})\right),\\
A_{2} & =\sqrt{2\pi}\frac{\rho}{4}\sqrt{\frac{1-\rho^{2}}{1+\rho^{2}}}\left(\Phi''(\frac{1+\rho}{2})+\Phi''(\frac{1-\rho}{2})\right).
\end{align*}
  
\end{theorem}
\begin{remark}
If $\Phi''$ is continuous  and hence  bounded on $[0,1]$, the bound above can be sharpened into that
\begin{align*}
\mathbf{Stab}_{\Phi}[f] & =A_{0}-A_{1}\epsilon+O(\epsilon^{\frac{2}{1+\rho^{2}}})
\end{align*}
for any balanced Boolean function $f$. This can be seen as follows.
Any balanced Boolean function $f$ can be expressed as $f(\mathbf{x})=1\{x_{i}=1\}+g(\mathbf{x})$,
where $g(\mathbf{x})=1\{x_{i}=-1\}1_{B_{2}}(x_{\backslash i})-1\{x_{i}=1\}1_{B_{1}}(x_{\backslash i})$
for some sets $B_{1},B_{2}\subseteq\{\pm1\}^{n-1}$ with same probability
$\epsilon$. Given $f$, by Taylor's theorem,
\begin{align*}
\Phi(T_{\rho}f(\mathbf{x})) & =\Phi(\frac{1+\rho x_{i}}{2})+\Phi'(\frac{1+\rho x_{i}}{2})T_{\rho}g(\mathbf{x})+\frac{1}{2}\Phi''(\xi)(T_{\rho}g(\mathbf{x}))^{2},
\end{align*}
where $\xi$ depends on $\mathbf{x}$. If $\Phi''$ is bounded, then
taking expectation for two sides of the equation above yields that
\begin{align*}
\mathbf{Stab}_{\Phi}[f] & =A_{0}-A_{1}\epsilon+O(\mathbb{E}[(T_{\rho}g)^{2}]),
\end{align*}
where by the hypercontractivity inequality (stated in the following lemma), $\mathbb{E}[(T_{\rho}g)^{2}]\le \|g\|_{1+\rho^2}^{2} = (2\epsilon)^{\frac{2}{1+\rho^{2}}}.$
So, for this case, 
\begin{align*}
\mathbf{Stab}_{\Phi}[f] & =A_{0}-A_{1}\epsilon+O(\epsilon^{\frac{2}{1+\rho^{2}}}).
\end{align*}
 Although this asymptotic equality is sharper, it requires the boundness
of $\Phi''$. Unfortunately, the most interesting cases, $\Phi_{q}$
and $\Phi_{q}^{\mathrm{sym}}$ with $q<2$, do not satisfy this condition. 
\end{remark}
The hypercontractivity inequality   mentioned above is given below.
\begin{lemma}[Hypercontractivity Inequality; see e.g., \cite{O'Donnell14analysisof}] \label{lem:HC} Let $T_\rho$ be the noise operator defined on $\{\pm 1\}^n$. Then, 
for $q>1$, it holds that 
 $\|T_{\rho}g\|_{q}\le\|g\|_{p}$ for all function $g:\{\pm 1\}^n \to \mathbb{R}$
where $p=1+(q-1)\rho^{2}$. 
\end{lemma}

Note that $A_{0}$ in Theorem \ref{thm:asymptotics} is exactly the $\Phi$-stability of dictator functions.
Moreover, if $\Phi$ is strictly convex, then $A_{1}>0$. So,  Theorem
\ref{thm:asymptotics} implies the local optimality of dictator functions. 
\begin{theorem}[Local Optimality of Dictator Functions]
\label{thm:asymptotics-1} Let $\Phi:[0,1]\to\mathbb{R}$ be a twice continuously
differentiable and strictly convex function. Let $\rho\in(0,1)$.
Then, there exists $\epsilon>0$ such that for any balanced Boolean
function $f$ satisfying $0<\min_{i\in[n]}\tilde{d}_{i}(f)\le\epsilon$,
it holds that 
\[
\mathbf{Stab}_{\Phi}[f]<\mathbf{Stab}_{\Phi}[f_{\mathrm{d}}].
\]
In other words, dictator functions are locally optimal in maximizing
the $\Phi$-stability. 
\end{theorem}

\subsection{$q$-Stability}

Theorem \ref{thm:asymptotics-1} shows that for sufficiently small
$\epsilon>0$, dictator functions have larger $\Phi$-stability than
any balanced Boolean function $f$ satisfying $0<\min_{i\in[n]}\tilde{d}_{i}(f)\le\epsilon$.
The value of $\epsilon$ here  can be determined by checking our
proof, but the final expression seems very complicated.  However,
if we focus on the $q$-stability, then we can obtain an explicit
 expression for such an $\epsilon$ by using hypercontractivity
inequalities.  Before stating this result, we need introduce some notations.

For $\epsilon\in[0,1]$, define 
\[
\Gamma_{q}(\epsilon):=\frac{\gamma_{q}(\epsilon)-\frac{1}{2}}{q-1},
\]
for $q>0$ but $q\neq1$, where
\[
\gamma_{q}(\epsilon):=\frac{1}{2}\left(\epsilon+\left(\frac{1+\rho}{2}\right)^{p}(1-2\epsilon)\right)^{q/p}+\frac{1}{2}\left(\epsilon+\left(\frac{1-\rho}{2}\right)^{p}(1-2\epsilon)\right)^{q/p}
\]
with $p=1+(q-1)\rho^{2}$. For $q=1$, define 
\begin{align*}
\Gamma_{1}(\epsilon) & :=\lim_{q\to1}\Gamma_{q}(\epsilon) =\frac{1}{2}(1-\rho^{2})h\left(\frac{1-\rho}{2}+\rho\epsilon\right)+(\frac{1}{2}-\epsilon)\rho^{2}h\left(\frac{1-\rho}{2}\right),
\end{align*}
where recall that $h(t)=t\ln t+(1-t)\ln(1-t).$ 

\begin{theorem}[Bound on  $q$-Stability]
\label{thm:qstability} Let $f$ be a balanced Boolean function.
Then, for any $q>0$ and $i\in[n]$,  it holds that 
\begin{equation}
\mathbf{Stab}_{q}[f]\le\Gamma_{q}(\tilde{d}_{i}(f)).\label{eq:-2-2}
\end{equation}
\end{theorem}

\begin{remark}
This theorem immediately implies a similar result for the symmetric
$q$-stability, i.e., 
\begin{align*}
\mathbf{Stab}_{q}^{\mathrm{sym}}[f] & \le2\Gamma_{q}(\tilde{d}_{i}(f)),\;q>0.
\end{align*}
\end{remark}

\begin{proof}
We first consider the case $q>1$. Using \eqref{eq:EPhi} with $\Phi$ set to $t\mapsto t^q$, we can write
\begin{align}
\mathbb{E}_{\mu}[(T_{\rho}f)^{q}]=\frac{1}{2}\left(\mathbb{E}_{\mu^{(n-1)}}[(T_{\rho}^{(n-1)}g_{+})^{q}]+\mathbb{E}_{\mu^{(n-1)}}[(T_{\rho}^{(n-1)}g_{-})^{q}]\right). \label{eq:ETf}
\end{align}
To upper bound the RHS above, we need the hypercontractivity inequality, restated in Lemma \ref{lem:HC}. That is,
\[
\mathbb{E}_{\mu}[(T_{\rho}g)^{q}]\le\mathbb{E}_{\mu}[g^{p}]^{q/p},\;\forall g,
\]
where $p=1+(q-1)\rho^{2}$.
By this inequality, 
\begin{align}
\mathbb{E}_{\mu^{(n-1)}}[(T_{\rho}^{(n-1)}g_{+})^{q}] & \le\mathbb{E}_{\mu^{(n-1)}}[g_{+}^{p}]^{q/p},\\
\mathbb{E}_{\mu^{(n-1)}}[(T_{\rho}^{(n-1)}g_{-})^{q}] & \le\mathbb{E}_{\mu^{(n-1)}}[g_{-}^{p}]^{q/p}. \label{eq:ETf2}
\end{align}

We now upper bound $\mathbb{E}_{\mu^{(n-1)}}[g_{+}^{p}]$ and $\mathbb{E}_{\mu^{(n-1)}}[g_{-}^{p}]$.
Let $A_{+}=\supp(f_{+})$ and $A_{-}=\supp(f_{-})$.  So, $\mu^{(n-1)}(A_{+})=1-d_{i}(f)$,
$\mu^{(n-1)}(A_{-})=d_{i}(f)$. Then, 
\begin{align*}
g_{+} & =1_{A_{+}\cap A_{-}}+\frac{1+\rho}{2}1_{A_{+}\backslash A_{-}}+\frac{1-\rho}{2}1_{A_{-}\backslash A_{+}},\\
g_{-} & =1_{A_{+}\cap A_{-}}+\frac{1-\rho}{2}1_{A_{+}\backslash A_{-}}+\frac{1+\rho}{2}1_{A_{-}\backslash A_{+}}.
\end{align*}
So,
\[
\mathbb{E}_{\mu^{(n-1)}}[g_{+}^{p}]=\mu^{(n-1)}(A_{+}\cap A_{-})+\left(\frac{1+\rho}{2}\right)^{p}\mu^{(n-1)}(A_{+}\backslash A_{-})+\left(\frac{1-\rho}{2}\right)^{p}\mu^{(n-1)}(A_{-}\backslash A_{+}).
\]
Denote $\kappa=\mu^{(n-1)}(A_{+}\cap A_{-})\in[0,\tilde{d}_{i}(f)]$.
Then, 
\[
\mathbb{E}_{\mu^{(n-1)}}[g_{+}^{p}]=\kappa+\left(\frac{1+\rho}{2}\right)^{p}(1-d_{i}(f)-\kappa)+\left(\frac{1-\rho}{2}\right)^{p}(d_{i}(f)-\kappa).
\]
Note that $\left(\frac{1+\rho}{2}\right)^{p}+\left(\frac{1-\rho}{2}\right)^{p}<1$
for $p>1$. Hence, the RHS is maximized at $\kappa=\tilde{d}_{i}(f)$,
i.e., 
\begin{align*}
\mathbb{E}_{\mu^{(n-1)}}[g_{+}^{p}] & \le\tilde{d}_{i}(f)+\left(\frac{1+\rho}{2}\right)^{p}(1-d_{i}(f)-\tilde{d}_{i}(f))+\left(\frac{1-\rho}{2}\right)^{p}(d_{i}(f)-\tilde{d}_{i}(f))\\
 & =\begin{cases}
d_{i}(f)+\left(\frac{1+\rho}{2}\right)^{p}(1-2d_{i}(f)), & d_{i}(f)\le1/2\\
1-d_{i}(f)+\left(\frac{1-\rho}{2}\right)^{p}(2d_{i}(f)-1), & d_{i}(f)>1/2
\end{cases}.
\end{align*}
Similarly,
\begin{align*}
\mathbb{E}_{\mu^{(n-1)}}[g_{-}^{p}] & \le\begin{cases}
d_{i}(f)+\left(\frac{1-\rho}{2}\right)^{p}(1-2d_{i}(f)), & d_{i}(f)\le1/2\\
1-d_{i}(f)+\left(\frac{1+\rho}{2}\right)^{p}(2d_{i}(f)-1), & d_{i}(f)>1/2
\end{cases}.
\end{align*}

Combining these estimates with \eqref{eq:ETf} and \eqref{eq:ETf2} yields that
\begin{align*}
\mathbb{E}_{\mu}[(T_{\rho}f)^{q}] & \le\frac{1}{2}\left(\mathbb{E}_{\mu^{(n-1)}}[g_{+}^{p}]^{q/p}+\mathbb{E}_{\mu^{(n-1)}}[g_{-}^{p}]^{q/p}\right)\\
 & \le\frac{1}{2}\left(\tilde{d}_{i}(f)+\left(\frac{1+\rho}{2}\right)^{p}(1-2\tilde{d}_{i}(f))\right)^{q/p} \\
 & \qquad +\frac{1}{2}\left(\tilde{d}_{i}(f)+\left(\frac{1-\rho}{2}\right)^{p}(1-2\tilde{d}_{i}(f))\right)^{q/p}.
\end{align*}

The inequality in \eqref{eq:-2-2} for $0<q<1$  follows similarly, but  
the reverse hypercontractivity inequality will be applied. Taking the
limit $q\downarrow1$ in \eqref{eq:-2-2}, we obtain the bound on
 $1$-stability.
\end{proof}

Based on the bound given in Theorem \ref{thm:qstability}, we next provide an  explicit range of $\tilde{d}(f)$ for which  dictator functions are optimal in maximizing $q$-stability. Here, recall that $\tilde{d}(f)$ is the minimum distance in the $L_1$ sense  between the balanced Boolean function $f$ and the set of dictator functions.  
That is,  given $q > 0$, our interest lies in determining the values of $\epsilon$ for which the   inequality
\begin{align}
\Gamma_{q}(\epsilon) & \le\mathbf{Stab}_{q}[f_{\mathrm{d}}]\label{eq:-6-1}
\end{align}
holds. 
Theorem \ref{thm:qstability} implies the optimality of dictator functions over all balanced Boolean functions with $\epsilon=\tilde{d}(f)$ satisfying \eqref{eq:-6-1}. This recovers and strengthens Theorem \ref{thm:asymptotics-1}.

To this end, we introduce the following notations. For $q \in (0,\infty)  \backslash \{1\}$
and $\rho\in(0,1)$, let $\epsilon_{q}^{*}(\rho)$ be the unique value
$\epsilon$ in $(0,1/2]$ attaining the equality in \eqref{eq:-6-1}
 if such $\epsilon$ exists, i.e., the unique solution
in $(0,1/2]$ (if the solution exists) to the equation 
\begin{equation}
\left(\epsilon+\left(\frac{1+\rho}{2}\right)^{p}(1-2\epsilon)\right)^{q/p}+\left(\epsilon+\left(\frac{1-\rho}{2}\right)^{p}(1-2\epsilon)\right)^{q/p}=\left(\frac{1+\rho}{2}\right)^{q}+\left(\frac{1-\rho}{2}\right)^{q}\label{eq:-34}
\end{equation}
with $p=1+(q-1)\rho^{2}$. By continuous extension, let $\epsilon_{q}^{*}(0):=\lim_{\rho\downarrow0}\epsilon_{q}^{*}(\rho)=1-\ln2$,
and let $\epsilon_{q}^{*}(1):=\lim_{\rho\uparrow1}\epsilon_{q}^{*}(\rho)$
which is the unique solution in $(0,1/2]$ to the equation 
\[
2(q-1)h(\epsilon)=\epsilon q^{2}
\]
(if the solution exists). Note that the base of the logarithm in the definition
of $h$ is $e$. 

By continuous extension, for $q=1$, let $\epsilon^{*}(\rho):=\epsilon_{1}^{*}(\rho)$
be the unique solution in $(0,1/2]$ to the equation 
\begin{equation}
h\left(\frac{1-\rho}{2}+\rho\epsilon\right)=\left(1+\frac{2\rho^{2}\epsilon}{1-\rho^{2}}\right)h\left(\frac{1-\rho}{2}\right),\label{eq:-33}
\end{equation}
where recall that $-h$ is the binary entropy function.

\begin{figure}
\centering \includegraphics[scale=0.5]{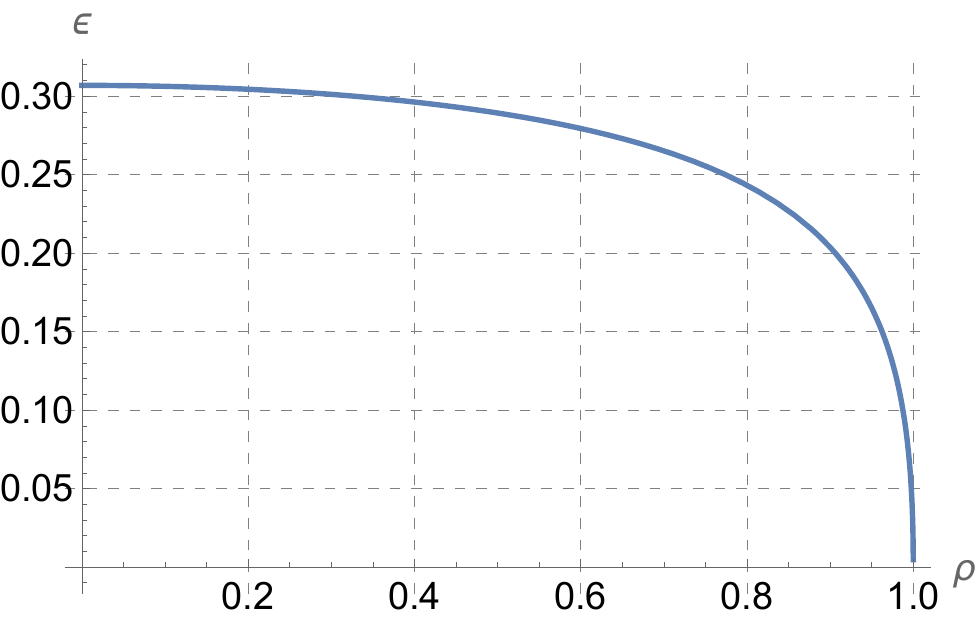}

\caption{\label{fig:The-region-of}The plot of the  function   $ \epsilon_{1}^{*}$.}
\end{figure}

Since the LHS in \eqref{eq:-34} and the LHS in \eqref{eq:-33} are
strictly convex in $\epsilon$ for $q\ge1$  and strictly
concave in $\epsilon$ for $0<q<1$, there are at most two solutions
to the equation in \eqref{eq:-34} or the equation in \eqref{eq:-33}.
Obviously, $\epsilon=0$ is one  of them. So, there is at most one
solution other than $0$. Moreover, this solution must be $\le1/2$
if it exists. This is because, the LHS in \eqref{eq:-34} with $\epsilon=1/2$
is larger than the RHS for $q>1$, and is smaller than the
RHS for $0<q<1$, which follows by substituting dictator functions
into the forward and reverse hypercontractivity inequalities. In fact,
numerical results show that for $q\ge6.5$ or $q<1$, there always
exists $\rho\in(0,1)$ for which no solution in $(0,1/2]$ to the
equation in \eqref{eq:-34} exists.  However, for $1\le q\le6.4$,
it seems that there always exists a solution in $(0,1/2]$ to the
equation in \eqref{eq:-34} or \eqref{eq:-33} (and hence, it is unique).
For $q\in[1,2]$, we provide a rigorous proof for this point in the
following lemma.  The function $\epsilon_{1}^{*}(\rho)$ is plotted
in Fig. \ref{fig:The-region-of}. 

\begin{lemma}
\label{lem:unique}The following hold. 
\begin{enumerate}
\item For $q=1$ and $\rho\in(0,1)$, there is a unique solution in $(0,1/2]$
to the equation in \eqref{eq:-33}.
\item For $q\in(1,2]$ and $\rho\in(0,1)$, there is a unique solution in
$(0,1/2]$ to the equation in \eqref{eq:-34}. 
\end{enumerate}
\end{lemma}

Moreover, for $q=1$, $\epsilon_{1}^{*}(\rho)$ satisfies monotonicity,
as observed from Fig. \ref{fig:The-region-of}. 
\begin{lemma}
\label{lem:monotonicity} For $q=1$, $\epsilon_{1}^{*}(\rho)$ is
decreasing in $\rho\in(0,1)$.
\end{lemma}
The proofs of Lemmas \ref{lem:unique} and \ref{lem:monotonicity}
are respectively provided in Appendices \ref{sec:Proof-of-Lemma}
and \ref{sec:Proof-of-Lemma-1}.

As mentioned above, Theorem \ref{thm:qstability}   implies the local optimality of
dictator functions, as shown in the following theorem.   Different from Theorem \ref{thm:asymptotics-1}, the following theorem
provides an explicit expression for a threshold $\epsilon$. 
\begin{theorem}[Local Optimality of Dictator Functions]
\label{thm:optimality-2} Let $q >0$ and
$\rho\in(0,1)$ such that $\epsilon_{q}^{*}(\rho)$ exists and $\epsilon_{q}^{*}(\rho)>0$.
Let $f$ be a balanced Boolean function $f$ satisfying $0\le \tilde{d}(f)\le\epsilon_{q}^{*}(\rho)$, where recall that $\tilde{d}(f):=\min_{i\in[n]}\tilde{d}_{i}(f)$.
Then, it holds that 
\begin{align}
\mathbf{Stab}_{q}[f] & \le\mathbf{Stab}_{q}[f_{\mathrm{d}}].
\end{align}
In other words, dictator functions maximize the $q$-stability   over all balanced
Boolean functions $f$ such that $0\le \tilde{d}(f)\le\epsilon_{q}^{*}(\rho)$. 
\end{theorem}

This theorem  still holds when the asymmetric $q$-stability is replaced
by the symmetric $q$-stability. 

\section{Global Optimality of Dictatorships}

We next investigate the global optimality of dictator functions in
maximizing the $q$-stability. Specifically, we prove improved bounds
for the Courtade--Kumar conjecture and the Li--M\'edard conjecture, which
respectively correspond to maximizing the $q$-stability with $q=1$
and $1<q<2$. 

\subsection{$q$-Stability with $q=1$ (Courtade--Kumar Conjecture)}

In this subsection, we focus on the case $q=1$ and will prove   the Courtade--Kumar Conjecture   for $\rho \in [0,0.914]$.
Our proof idea is as follows. By Theorem \ref{thm:optimality-2}, dictator functions are locally optimal in maximizing the $q$-stability. Hence, we only need to consider the Boolean functions far from dictator functions. 
Fortunately, due to the spread of the spectrum of such Boolean functions $f$, the energy of $T_{\rho}f$ 
is significantly smaller than the energy of $T_{\rho}f_{\mathrm{d}}$ (especially for small $\rho$) where $f_{\mathrm{d}}$ is a dictator function. Utilizing this property, the present author established a bound on $\Phi$-stability in \cite{yu2023phi}, which proves effective for such Boolean functions and hence will serve as a key tool in proving our new bounds.

\begin{figure}
\centering \includegraphics[scale=0.5]{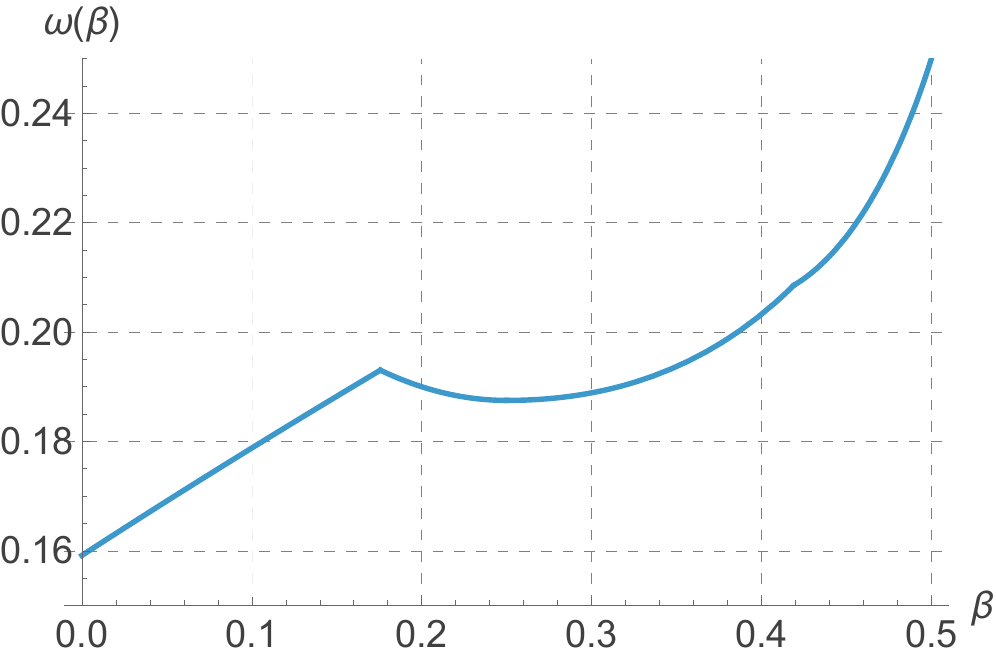}
\caption{\label{fig:omega}Illustration of the function   $\omega$.}
\end{figure}

We first introduce  the previous bound derived by the present author in \cite{yu2023phi}.
Define for $0\leq\beta\leq\frac{1}{2}$,
\begin{equation}
\omega(\beta):=\min\Big\{\beta^{2}+\varphi(\frac{1}{2}-\beta),\,\frac{\left(1+\sqrt{1+4(\pi-\sqrt{2\pi})\beta}\right)^{2}}{8\pi}\Big\},\label{eq:-21-1}
\end{equation}
where  for $0\leq t\leq1$, 
\[
\varphi(t):=\min\{\varphi_{\mathrm{C}}(\tilde{t}),\varphi_{\mathrm{LP}}(\tilde{t})\},
\]
with  $\tilde{t}:=\min\{t,1-t\}$ and\footnote{Here, the subscript ``$\mathrm{C}$'' refers to ``Chang'', since
$\varphi_{\mathrm{C}}$ was used to bound the first-level Fourier
weight of a Boolean function in Chang's lemma \cite{chang2002polynomial,O'Donnell14analysisof}.
The subscript ``$\mathrm{LP}$'' refers to ``Linear Programming'',
since $\varphi_{\mathrm{LP}}$ was used to bound the first-level Fourier
weight by linear programming methods \cite{fu2001minimum,yu2019improved}.} 
\begin{align*}
\varphi_{\mathrm{C}}(t) & :=2t^{2}\ln\frac{1}{t},\;0\leq t\leq\frac{1}{2},\\
\varphi_{\mathrm{LP}}(t) & :=\begin{cases}
2t^{2}(\frac{1}{\sqrt{t}}-1), & 0<t\leq\frac{1}{4}\\
\frac{t}{2}, & \frac{1}{4}\leq t\leq\frac{1}{2}
\end{cases}. 
\end{align*}
The quantity $\omega(\beta)$ in fact forms an upper bound on the   first-level Fourier weight of a balanced Boolean function for which a particular first-level Fourier coefficient has absolute value  $\beta$ \cite{yu2023phi}. That is,
\begin{equation}
\mathbf{W}_{1}:= \sum_{i=1}^n \hat{f}_{\{i\}}^2\le\omega(\beta)\leq\frac{1}{4} \label{eq:FKN}
\end{equation}
holds for all balanced Boolean function $f$ with $|\hat{f}_{\{i\}}|=\beta $ for some $i\in [n]$.  This is a quantitative version of Friedgut-Kalai-Naor (FKN) theorem.
The function   $\omega$ is plotted in Fig. \ref{fig:omega}. 

Define for $0\leq\beta\leq\frac{1}{2}$, 
\begin{equation}
\Upsilon_{\rho}(\beta):=\sup_{(z_{1},z_{2})\in\mathcal{Z}}\gamma(z_{1},z_{2},\beta),\label{eq:-49}
\end{equation}
where
\begin{equation}
\mathcal{Z}:=\{(z_{1},z_{2}):-\frac{1}{2\rho}< z_{1}\leq z_{2} < \frac{1}{2\rho},0\leq p_{1}\leq\frac{1}{4}+\frac{\beta}{2},0\leq p_{2}\leq\frac{1}{4}-\frac{\beta}{2}\},\label{eq:-19}
\end{equation}
and 
\[
\gamma(z_{1},z_{2},\beta):=(1-2p_{1}-2p_{2})\Phi(0)+2p_{1}\Phi(\frac{1}{2}+\rho z_{1})+2p_{2}\Phi(\frac{1}{2}+\rho z_{2}),
\]
with 
\begin{align}
 & p_{1}=\frac{(1-\rho)(1+\rho-4\rho^{2}\omega(\beta)+2\beta(1+2\rho z_{2}-\rho^{2}))}{4(1+2\rho z_{1})(1+\rho z_{2}-\rho z_{1}-\rho^{2})},\label{eq:p-2}\\
 & p_{2}=\frac{(1-\rho)(1+\rho-4\rho^{2}\omega(\beta)-2\beta(1-2\rho z_{1}-\rho^{2}))}{4(1-2\rho z_{2})(1+\rho z_{2}-\rho z_{1}-\rho^{2})}.\label{eq:q-2}
\end{align}

\begin{theorem}[\cite{yu2023phi}]
\label{thm:optimality-1} If $\Phi$ is symmetric
and continuous on $[0,1]$ and differentiable on $(0,1)$ whose derivative
$\Phi'$ is increasing and continuous on $(0,1)$, and strictly concave
on $(0,\frac{1}{2}]$, then for any $\rho\in(0,1)$ and any balanced
Boolean function $f$, it holds that $\mathbf{Stab}_{\Phi}[f]\leq\Upsilon_{\rho}(\beta),$
where $\beta=\max_{i\in[n]}|\hat{f}_{\{i\}}|$. 
\end{theorem}
This theorem was  proven by Boolean Fourier analysis,  particularly
by using an energy control technique and a quantitative version of
the Friedgut--Kalai--Naor (FKN) theorem.  To be self-contained, here we provide a proof sketch of this theorem. The complete 
proof can be found in \cite{yu2023phi}.
\begin{proof}[Proof sketch of Theorem \ref{thm:optimality-1}]
Let $f$ be a balanced Boolean function.  Recall that $\{\hat{f}_{T}:T\subseteq[n]\}$ are the Fourier coefficients of $f$. 
Denote the \emph{degree-$k$ part} of $f$ as 
\[
f_{k}(\mathbf{x}):=\sum_{T\subseteq[n]:|T|=k}\hat{f}_{T}\chi_{T}(\mathbf{x}),
\]
and the \emph{degree-$k$ Fourier weight }of $f$ as 
\begin{align}
\mathbf{W}_{k}:=\mathbf{W}_{k}[f] & := \mathbb{E}[ f_{k}^2(\mathbf{X})] =\sum_{T:|T|=k}\hat{f}_{T}^{2},\quad k\in[0:n].\label{eq:FourierWeight}
\end{align}
It is evident that the functions $f_k$ for $k \in [0:n]$ are mutually orthogonal, a direct consequence of the orthogonality inherent to the Fourier basis.

Without loss of generality, we assume $|\hat{f}_{\{1\}}|\geq|\hat{f}_{\{i\}}|$
for all $2\le i\le n$. If $\hat{f}_{\{1\}}\ge0$, denote $\beta=\hat{f}_{\{1\}}$
and $X:=X_{1}$; otherwise, $\beta:=-\hat{f}_{\{1\}}$ and $X:=-X_{1}$.
Hence $\beta\ge0$. Define several other random variables: 
\begin{align*}
S & :=f(\mathbf{X})-\frac{1}{2}=\sum_{k=1}^{n}f_{k}(\mathbf{X}), \quad Z  :=\sum_{k=1}^{n}\rho^{k-1}f_{k}(\mathbf{X}).
\end{align*}
Since $\mathbb{E}[SX]=\beta$, $(S,X)$ satisfies the following joint distribution 
\begin{equation*}
P_{SX}(s,x)=\begin{cases}
\frac{1+2\beta}{4}, & (s,x)= \pm (1/2,1)\\
\frac{1-2\beta}{4}, & (s,x)= \pm (1/2,-1)
\end{cases}. 
\end{equation*} 
The energy of $Z$ can be controlled as follows:
\begin{align*}
\mathbb{E}[Z^{2}] =\mathbf{W}_{1}+\rho^{2}\sum_{k=2}^{n}\rho^{2(k-2)}\mathbf{W}_{k}   & \leq(1-\rho)\mathbf{W}_{1}+\rho\left(\mathbf{W}_{1}+\rho \sum_{k=2}^{n}\rho^{k-2}\mathbf{W}_{k}\right)\\
 & \leq(1-\rho)\omega(\beta)+\rho\mathbb{E}[SZ],  
\end{align*}
where in the last line, the quantitative version of FKN theorem in \eqref{eq:FKN} is applied. Moreover, by definition, 
\begin{align*}
\mathbb{E}Z & =0, \quad \mathbb{E}[XZ]=\beta, \quad  T_{\rho}f(\mathbf{X})=\sum_{k=0}^{n}\rho^{k}f_{k}(\mathbf{X})=\frac{1}{2}+\rho Z, 
\end{align*}

Summarizing all the properties listed above on $Z$, we obtain the following bound on $\Phi$-stability:  $\mathbf{Stab}_{\Phi}[f] \leq \tilde{\Upsilon}_{\rho}(\beta)$,
where $\tilde{\Upsilon}_{\rho}(\beta)$ is defined as the supremum of $\mathbb{E}[\Phi(\frac{1}{2}+\rho Z)]$ over all random variables $Z$ (equivalently, over all conditional distributions $P_{Z|SX}$) defined on any finite set such that 
\begin{align}
 & 0\leq \frac{1}{2} +\rho Z\leq1\textrm{ a.s.}, \quad  \mathbb{E}Z=0, \quad  \mathbb{E}[XZ]=\beta, \quad  \mathbb{E}[Z^{2}]\leq(1-\rho)\omega(\beta)+\rho\mathbb{E}[SZ]. \label{eq:constraints}
\end{align}

We next simplify the bound $\tilde{\Upsilon}_{\rho}(\beta)$ and prove that $\tilde{\Upsilon}_{\rho}(\beta)= \Upsilon_{\rho}(\beta)$. By a proper version of Carathéodory's theorem,  the optimization in the definition of $\tilde{\Upsilon}_{\rho}(\beta)$ can be converted into   a finite-dimensional program and the supremum therein is  a maximum, since the feasible region is compact and the objective
and constraint functions are continuous. Moreover, the maximum is attained only 
when the last constraint in \eqref{eq:constraints} is an equality. 
We can assume that  the linear independence constraint qualification (LICQ), a sufficient condition for the Karush--Kuhn--Tucker
(KKT) conditions to be valid, is satisfied by all optimal solutions to this finite program, since  all feasible solutions violating the  LICQ are  suboptimal (more precisely, they cannot be better than the best one satisfying the LICQ).  Applying KKT conditions and excluding all suboptimal solutions satisfying these conditions,  one can find that the optimal solutions 
have the following form:
\begin{equation*}
P_{SXZ}(s,x,z)=\begin{cases}
\frac{1+2\beta}{4}-p_1, & (s,x,z)=\pm(-\frac{1}{2},-1,-\frac{1}{2\rho})\\
p_1, & (s,x,z)=\pm(-\frac{1}{2},-1,z_{1})\\
\frac{1-2\beta}{4}-p_2, & (s,x,z)=\pm(\frac{1}{2},-1,\frac{1}{2\rho})\\
p_2, & (s,x,z)=\pm(\frac{1}{2},-1,z_{2})
\end{cases}\label{eq:PSXZ2}
\end{equation*}
for some $z_{1},z_{2}$ such that $-\frac{1}{2\rho}\leq z_{1}\leq z_{2}\leq\frac{1}{2\rho}$. 
For the case of $-\frac{1}{2\rho}< z_{1}\leq z_{2}<\frac{1}{2\rho}$, solving the last two constraints in \eqref{eq:constraints} (with the last one set to an equality) gives the solution $p_1,p_2$
in \eqref{eq:p-2} and \eqref{eq:q-2}, which yields $\Upsilon_{\rho}$. The case 
  $ z_{1} = -\frac{1}{2\rho}$ or the case $z_{2}=\frac{1}{2\rho}$ is  covered via taking limits. 
\end{proof}

We in fact use this theorem
to prove that the Courtade--Kumar conjecture (as well as the symmetrized
Li--M\'edard conjecture) is true for all $\rho\in[0,0.46]$ \cite{yu2023phi}.
The intuition behind the proof given above  is that if $f$ is far from dictator
functions, then $\beta$ would be small, due to the following identity $|\hat{f}_{\{i\}}|=\frac{1}{2}-\tilde{d}_{i}(f)$, i.e., 
\[
\beta=\frac{1}{2}-\tilde{d}(f).
\]
By the quantitative version of FKN theorem in \eqref{eq:FKN} and  
the energy control  (i.e., the last constraint in \eqref{eq:constraints}), 
in this case, the energy of $Z$ (defined in the proof sketch above), or equivalently, the energy of $T_{\rho}f$ (i.e., the $2$-stability  $\mathbb{E}[(T_{\rho}f)^2]$) 
would be small, especially for small $\rho$. This in
turn decreases the $\Phi$-stability (due to the facts that $\mathbf{Stab}_{\Phi}[f]=\mathbb{E}[\Phi(\frac{1}{2}+\rho Z)]$ and $\Phi$ is convex). Hence, to maximize  the $\Phi$-stability, 
$f$ must be close to or exactly equal to a dictator function.

Bounding the energy of $T_{\rho}f$ is easy for small $\rho$, but
not easy for large $\rho$. This is because, when  $\rho$ is close to
$1$, the energy of $T_{\rho}f$ is close to the energy of $f$ (i.e.,
$1/2$) which is always large. This is the reason why the bound in Theorem \ref{thm:optimality-1} does not work effectively for small $\rho$, in the aspect of excluding the Boolean functions close to dictator
functions from being optimal (i.e., showing the local 
optimality of dictator functions).
Fortunately, Theorem  \ref{thm:optimality-2}  can also serve the same purpose.  Specifically, by Theorem  \ref{thm:optimality-2}, if
$0\le \tilde{d}(f)\le\epsilon^{*}(\rho)$, then
$\mathbf{Stab}_{1}^{\mathrm{sym}}[f]\le\mathbf{Stab}_{1}^{\mathrm{sym}}[f_{\mathrm{d}}]$.
On the other hand, if $\tilde{d}(f)>\epsilon^{*}(\rho)$,
i.e., 
$
\beta=\frac{1}{2}-\tilde{d}(f)<\frac{1}{2}-\epsilon^{*}(\rho),
$
then by Theorem \ref{thm:optimality-1} with $\Phi=\Phi_{1}^{\mathrm{sym}}$,
\begin{equation}
\mathbf{Stab}_{1}^{\mathrm{sym}}[f]\leq\max_{\beta\in[0,\frac{1}{2}-\epsilon^{*}(\rho)]}\Upsilon_{\rho}(\beta).
\end{equation}
So, combining these two cases, we get a bound for the $1$-stability:
\begin{equation}
\mathbf{Stab}_{1}^{\mathrm{sym}}[f]\leq \max \left\{ \mathbf{Stab}_{1}^{\mathrm{sym}}[f_{\mathrm{d}}], \max_{\beta\in[0,\frac{1}{2}-\epsilon^{*}(\rho)]}\Upsilon_{\rho}(\beta)\right\},\label{eq:-3}
\end{equation}
which can be relaxed to the following simple one. 
\begin{theorem}
\label{thm:bound} For
any $\rho\in(0,1)$ and any balanced Boolean function $f$, it holds
that 
\begin{equation}
\mathbf{Stab}_{1}^{\mathrm{sym}}[f]\leq \max \left\{ \mathbf{Stab}_{1}^{\mathrm{sym}}[f_{\mathrm{d}}], \max_{t\in[0,1)}  \eta_1(t), \max_{t\in[0,1)} \eta_2(t)\right\},\label{eq:-3-1}
\end{equation}
where 
\begin{align*}
\eta_1(t) & := \frac{\left(1-\rho\right)\left(1+\rho-4\rho^{2}\omega^{*}(\rho)\right)}{2\left(1+t-\rho^{2}\right)}\phi(\frac{1-t}{2}), \\
\eta_2(t) & :=\frac{1+\rho-4\rho^{2}\omega^{*}(\rho)- (1- 2 \epsilon^{*}(\rho))\rho^{2} t}{2(1+\rho)(1+t)}  \phi(\frac{1-t}{2}),\\
\omega^{*}(\rho) & :=\max_{\beta\in[0,\frac{1}{2}-\epsilon^{*}(\rho)]}\omega(\beta),\\
\phi(s) & :=\frac{\Phi_{1}^{\mathrm{sym}}(s)}{s}=\frac{s\ln s+(1-s)\ln(1-s)}{s},\,s\in(0,1].
\end{align*}
\end{theorem}
The proof of this theorem  is
provided in Appendix \ref{sec:Proof-of-Theorem-1}. The last two constraints in \eqref{eq:-19} have been discarded in  proving 
the bound in \eqref{eq:-3-1}, which might make it weaker than the one in \eqref{eq:-3}. However, numerical results indicate that as $\rho$ increases beyond a critical threshold ($\rho = 0.914...$), the value $\max\left\{\max_{t\in[0,1)} \eta_1(t), \max_{t\in[0,1)} \eta_2(t)\right\}$ surpasses $\mathbf{Stab}_{1}^{\mathrm{sym}}[f_{\mathrm{d}}]$. Notably, the discarded constraints remain inactive near this critical point. In other words, this simplified bound retains full optimality for the specific purposes of this paper.
Numerical results 
also indicate that  the bound in \eqref{eq:-3-1} is better than  the bound in Theorem \ref{thm:optimality-1} for  large $\rho$, which is intuitively obvious due to  the incorporation of  Theorem \ref{thm:optimality-2}.

By this theorem, if for a given $\rho$, the maximum at the RHS of
\eqref{eq:-3-1} is no larger than $\mathbf{Stab}_{1}^{\mathrm{sym}}[f_{\mathrm{d}}]=\Phi_{1}^{\mathrm{sym}}(\frac{1+\rho}{2})$,
then dictator functions are globally optimal in maximizing the $1$-stability
over all balanced Boolean functions.  In \cite{yu2023phi}, we have
already proven that the Courtade--Kumar conjecture is true for $\rho\in[0,0.46]$.
So, here we only focus on the case $\rho\ge0.46$. 
\begin{corollary}
\label{cor:If-for-a} If for a given $\rho$, 
\begin{align}
1+\rho-4\rho^{2}\omega^{*}(\rho)- (1- 2 \epsilon^{*}(\rho)) (2-\rho^2) & \ge 0,\label{eq:condition}  \\ 
\max_{\substack{t\in[0,1)}
}\left(1+\rho-4\rho^{2}\omega^{*}(\rho)\right)\phi(\frac{1-t}{2})-\left(1+t-\rho^{2}\right)\phi(\frac{1-\rho}{2})\, & \le0,\label{eq:-26}
\end{align}
then $\mathbf{Stab}_{1}^{\mathrm{sym}}[f]\le\mathbf{Stab}_{1}^{\mathrm{sym}}[f_{\mathrm{d}}]$
for all balanced Boolean functions $f$. In particular, if for some
$\rho^{*}\in[0.46,1]$, 
\begin{align}
\min_{\rho\in[0.46,\rho^{*}]} 1+\rho-4\rho^{2}\omega_{0}- (1- 2 \epsilon_{0}) (2-\rho^2) &\ge 0, \label{eq:condition2} \\
\max_{\rho\in[0.46,\rho^{*}]}\max_{\substack{t\in[0,1)}
}\chi(\rho,t)\,& \le0,\label{eq:-29}
\end{align}
where the function $\chi$ is given by 
\[
\chi(\rho,t):=\left(1+\rho-4\rho^{2}\omega_{0}\right)\phi(\frac{1-t}{2})-\left(1+t-\rho^{2}\right)\phi(\frac{1-\rho}{2}),
\]
$\epsilon_{0}$ is a number such that $\epsilon^{*}(\rho)\ge \epsilon_{0},\forall \rho\in[0.46,\rho^{*}]$, and  $\omega_{0}$ is a number such that $\omega^{*}(\rho)\le\omega_{0},\forall \rho\in[0.46,\rho^{*}]$, then for all $\rho\in[0.46,\rho^{*}]$,
$\mathbf{Stab}_{1}^{\mathrm{sym}}[f]\le\mathbf{Stab}_{1}^{\mathrm{sym}}[f_{\mathrm{d}}]$
for all balanced Boolean functions $f$. 
\end{corollary}
\begin{remark}
\label{rem:Note-that-since}Note that   $t \in [0,1) \mapsto\chi(\rho,t)$
is concave (see Lemma \ref{lem:phi}), and hence the optimal $t$ attaining the maximum in \eqref{eq:-29}
is the unique solution $t_{\rho}\in(0,1)$ to the equation
\begin{equation}
-\frac{1}{2}\left(1+\rho-4\rho^{2}\omega_{0}\right)\phi'(\frac{1-t}{2})-\phi(\frac{1-\rho}{2})=0.\label{eq:-30}
\end{equation}
\end{remark}
\begin{proof}
It is straightforward to see that $\eta_2(t) \le \eta_1(t),\forall t\in [0,1)$ (defined in Theorem \ref{thm:bound}) if \eqref{eq:condition} holds. 
If we can further show   $\eta_1(t)\le\Phi_{1}^{\mathrm{sym}}(\frac{1+\rho}{2}),\forall t\in [0,1)$,
or equivalently, the condition in \eqref{eq:-26},
then by Theorem \ref{thm:bound}, $\mathbf{Stab}_{1}^{\mathrm{sym}}[f]\le\mathbf{Stab}_{1}^{\mathrm{sym}}[f_{\mathrm{d}}]$
for all balanced Boolean functions $f$.  Since $\phi$ is negative,
the second statement follows obviously. 
\end{proof}
In the following, we analytically verify \eqref{eq:condition2} and numerically verify \eqref{eq:-29} for proper
choices of $\rho^{*}$ and $\omega_{0}$. Specifically, we choose
\begin{equation}
\rho^{*}=0.914,\label{eq:-21}
\end{equation}
for which we observe that $\epsilon^{*}(\rho )\ge \epsilon_0:=\epsilon^{*}(0.914)=0.195055...$ and
$\omega^{*}(\rho)\le\omega_{0}:=\omega^{*}(0.914)=\omega(\beta_{0})=0.193026...$
for all $\rho\in[0.46,0.914]$, with $\beta_{0}=0.175661...$ 
Note that the objective function in \eqref{eq:condition2} is a  quadratic polynomial. We hence easily prove that the inequality in  \eqref{eq:condition2} holds in this case. On the other hand, using
Matlab, we numerically verify that 
\begin{equation}
\max_{\rho\in[0.46,0.914]}\max_{\substack{t\in[0,1)}
}\chi(\rho,t)=-0.00169063...<0,\label{eq:max}
\end{equation}
with the maximum attained at $\rho=0.914$ and $t=0.663100...$
  Hence, for all $0.46\le\rho\le0.914$, dictator functions maximize
the $1$-stability over all balanced Boolean functions.  Combining
it with our previous result in \cite{yu2023phi}, we conclude that
for all $0\le\rho\le0.914$, dictator functions maximize the $1$-stability
over all balanced Boolean functions.  

We now use a computer-assisted method to prove the numerical observation
above. 
\begin{theorem}
\label{thm:CK} The optimal value of \eqref{eq:max} with $\Phi=\Phi_{1}^{\mathrm{sym}}$
and $\rho^{*}=0.914$ is negative.  In other words, for any $\rho\in[0.46,0.914]$,
it holds that $\mathbf{Stab}_{1}^{\mathrm{sym}}[f]\le\mathbf{Stab}_{1}^{\mathrm{sym}}[f_{\mathrm{d}}]$
for all balanced Boolean functions $f$. 
\end{theorem}
Combining this theorem with the result in \cite{yu2023phi} indicates
that dictator functions maximize the $1$-stability over all balanced
Boolean functions for all $\rho\in[0,0.914]$. 
\begin{theorem}[Restatement of Theorem \ref{thm:CK2-2}]
\label{thm:CK2}The Courtade--Kumar
conjecture is true for all $\rho\in[0,0.914]$.
\end{theorem}
The Courtade--Kumar conjecture for the case of $\rho\in(0.914,1)$
is still open.

Our proof idea for Theorem \ref{thm:CK} is as follows. We first show
that $\partial_{\rho}\chi$ in \eqref{eq:max} is bounded on the feasible
region (say, belonging to $[-M,M]$). Then, for any $\rho$, $\hat{\rho}$,
and $t$, it holds that 
\begin{equation}
\chi(\hat{\rho},t)\le\chi(\rho,t)+M|\hat{\rho}-\rho|.\label{eq:-24}
\end{equation}
That is, if for any $\delta>0$, we verify $\chi(\rho,t)<-\delta$
for all $\rho\in[0.46,0.914]\cap(\frac{\delta}{M}\mathbb{Z}$), then
we can conclude that $\chi(\hat{\rho},t)<0$ for all $\hat{\rho}\in[0.46,0.914]$.
The detailed proof is provided in Appendix \ref{sec:Proof-of-Theorem-1-1}. 

\subsection{$q$-Stability with $1<q<2$ (Li--M\'edard Conjecture)}

We next focus on $q$-stability for $1 < q < 2$ and prove the Li--M\'edard Conjecture for $q \in [1.36, 2)$. We build on the approach used for the case $q = 1$, extending it to the present setting. Specifically, we combine Theorems \ref{thm:optimality-2} and \ref{thm:optimality-1} to establish the global optimality of dictator functions. Theorem \ref{thm:optimality-2} ensures that any extremizer maximizing $q$-stability among Boolean functions close to some dictator function must itself be a dictator function. 
Meanwhile, Theorem \ref{thm:optimality-1} ensures that any Boolean function sufficiently distant from all dictator functions cannot exhibit a higher $q$-stability than a dictator function.

We first establish an analogue of Theorem \ref{thm:bound}, i.e., a simplified version of Theorem \ref{thm:optimality-1}, for $q$-stability. By Theorem 
\ref{thm:optimality-2}, if $0\le\tilde{d}(f)\le\epsilon_{q}^{*}(\rho)$,
then $\mathbf{Stab}_{q}^{\mathrm{sym}}[f]\le\mathbf{Stab}_{q}^{\mathrm{sym}}[f_{\mathrm{d}}]$.
On the other hand, if $\tilde{d}(f)>\epsilon_{q}^{*}(\rho)$,
i.e., 
$
\beta=\max_{i\in[n]}|\hat{f}_{\{i\}}|=\frac{1}{2}-\tilde{d}(f)<\frac{1}{2}-\epsilon_{q}^{*}(\rho),
$
then by Theorem \ref{thm:optimality-1} with $\Phi=\Phi_{q}^{\mathrm{sym}}$,
\begin{equation}
\mathbf{Stab}_{q}^{\mathrm{sym}}[f]\leq\max_{\beta\in[0,\frac{1}{2}-\epsilon_{q}^{*}(\rho)]}\Upsilon_{\rho}(\beta).\label{eq:-3-2}
\end{equation}
Similarly to the case $q=1$, we  can simplify this bound as follows.

\begin{theorem}
\label{thm:bound-1} Let $\Phi=\Phi_{q}^{\mathrm{sym}}$. Then, for
any $q\in(1,2)$, $\rho\in(0,1)$, and any balanced Boolean function
$f$, it holds that 
\begin{equation}
\mathbf{Stab}_{q}^{\mathrm{sym}}[f]\leq \max \left\{ \mathbf{Stab}_{q}^{\mathrm{sym}}[f_{\mathrm{d}}], \max_{t\in[0,1]}  \eta_{q,1}(t), \max_{t\in[0,1]} \eta_{q,2}(t)\right\},\label{eq:-3-1-1}
\end{equation}
where 
\begin{align*}
\eta_{q,1}(t) & := \frac{\left(1-\rho\right)\left(1+\rho-4\rho^{2}\omega_{q}^{*}(\rho)\right)}{2\left(1+t-\rho^{2}\right)}\phi_{q}(\frac{1-t}{2}), \\
\eta_{q,2}(t) & :=\frac{1+\rho-4\rho^{2}\omega_{q}^{*}(\rho)- (1- 2 \epsilon_{q}^{*}(\rho))\rho^{2} t}{2(1+\rho)(1+t)}  \phi_{q}(\frac{1-t}{2}),\\
\omega_{q}^{*}(\rho) & :=\max_{\beta\in[0,\frac{1}{2}-\epsilon_{q}^{*}(\rho)]}\omega(\beta),\\
\phi_{q}(s) & :=\frac{\Phi_{q}^{\mathrm{sym}}(s)}{s}=\begin{cases}
\frac{s\ln_{q}s+(1-s)\ln_{q}(1-s)}{s}, & s\in(0,1]\\
\frac{-q}{q-1}, & s=0
\end{cases}.
\end{align*}
\end{theorem}

The proof of this theorem  is
similar to that of Theorem \ref{thm:bound}, and is provided in Appendix
\ref{sec:Proof-of-Theorem-1-2}.  Specifically,  Lemma \ref{lem:phi} and the proof steps of Theorem \ref{thm:bound} given in  Section \ref{sec:Proof-of-Theorem-1}  hold verbatim with the substitutions $\phi \leftarrow \phi_q, \ln \leftarrow \ln_q,  \omega^{*}\leftarrow \omega_{q}^{*}, \epsilon^{*} \leftarrow \epsilon_{q}^{*}$. Similarly to the case $q=1$, in our proof, the last two constraints in \eqref{eq:-19} have been discarded, but numerical results indicate that this simplified bound retains full optimality for the specific purposes of this paper.

By this theorem, if for a given $\rho$, the maximum at the RHS of
\eqref{eq:-3-1-1} is no larger than $\mathbf{Stab}_{q}^{\mathrm{sym}}[f_{\mathrm{d}}]=\Phi_{q}^{\mathrm{sym}}(\frac{1+\rho}{2})$,
then dictator functions are globally optimal in maximizing the $q$-stability
over all balanced Boolean functions. Since given $\rho$, if dictator
functions maximize the $q$-stability for some $q\in[1,2)$, then dictator
functions maximize the $\hat{q}$-stability for any $\hat{q}\in[q,2)$;
see \cite{barnes2020courtade,yu2023phi}. So, Theorem \ref{thm:CK}
in fact implies that for any $\rho\in[0,0.914]$, dictator functions
maximize the $q$-stability for any $q\in(1,2)$. Here, we only need focus
on the case of $\rho\in[0.914,1]$. 
\begin{corollary}
\label{cor:If-for-a-1} Given $q\in(1,2)$ and $\rho\in(0,1)$, if
\begin{align}
1+\rho-4\rho^{2}\omega_{q}^{*}(\rho)- (1- 2 \epsilon_{q}^{*}(\rho)) (2-\rho^2) & \ge 0, \label{eq:condition3} \\
\max_{\substack{t\in[0,1]}
}\left(1+\rho-4\rho^{2}\omega_{q}^{*}(\rho)\right)\phi_{q}(\frac{1-t}{2})-\left(1+t-\rho^{2}\right)\phi_{q}(\frac{1-\rho}{2})\,&\le0,\label{eq:-27}
\end{align}
then 
\begin{equation}
\mathbf{Stab}_{q}^{\mathrm{sym}}[f]\le\mathbf{Stab}_{q}^{\mathrm{sym}}[f_{\mathrm{d}}]\label{eq:-28}
\end{equation}
for all balanced Boolean functions $f$. In particular, given $q\in(1,2)$,
if 
\begin{align}
\min_{\rho\in[0.914,1]} 1+\rho-4\rho^{2}\omega_{0}- (1- 2 \epsilon_{0}) (2-\rho^2) &\ge 0, \label{eq:condition4} \\
\max_{\rho\in[0.914,1]}\max_{\substack{t\in[0,1]}
}\chi_{q}(\rho,t)\,&\le0,\label{eq:-29-1}
\end{align}
where  the function $\chi_{q}$ is given by 
\[
\chi_{q}(\rho,t):=\left(1+\rho-4\rho^{2}\omega_{0}\right)\phi_{q}(\frac{1-t}{2})-\left(1+t-\rho^{2}\right)\phi_{q}(\frac{1-\rho}{2}),
\]
$\epsilon_{0}$ is a number such that $\epsilon_{q}^{*}(\rho)\ge \epsilon_{0},\forall \rho\in[0.914,1]$, and $\omega_{0}$ is  a number such that $\omega_{q}^{*}(\rho)\le\omega_{0},\forall \rho\in[0.914,1]$, then for all $\rho\in[0.914,1]$, $\mathbf{Stab}_{q}^{\mathrm{sym}}[f]\le\mathbf{Stab}_{q}^{\mathrm{sym}}[f_{\mathrm{d}}]$
for all balanced Boolean functions $f$. 
\end{corollary}

\begin{remark}
\label{rem:Note-that-since-1}Note that   $t \in [0,1] \mapsto\chi_{q}(\rho,t)$
is concave, and  hence the optimal $t$ attaining the maximum in \eqref{eq:-29-1}
is the unique solution $t_{q,\rho}\in(0,1)$ to the equation
\begin{equation}
-\frac{1}{2}\left(1+\rho-4\rho^{2}\omega_{0}\right)\chi_{q}'(\frac{1-t}{2})-\chi_{q}(\frac{1-\rho}{2})=0.\label{eq:-30-1}
\end{equation}
\end{remark}

\begin{figure}
\centering \includegraphics[scale=0.5]{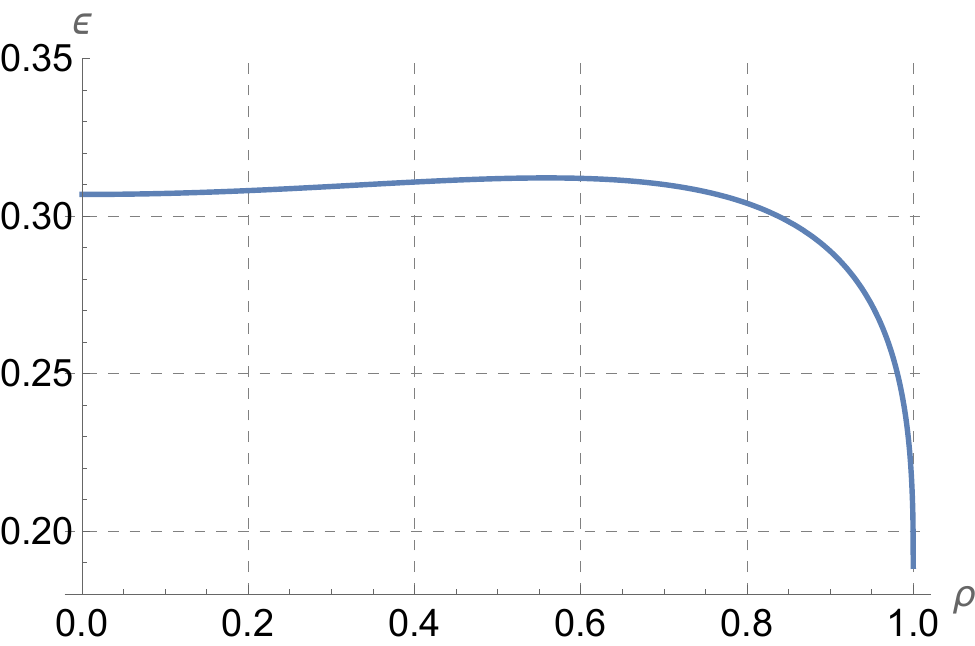}

\caption{\label{fig:The-region-of-1}The plot of the  function   $ \epsilon_{q}^{*} $
with $q=1.36$.}
\end{figure}

\begin{proof}
The proof is similar  to that of Corollary \ref{cor:If-for-a}. It is straightforward to see that $\eta_{q,2}(t) \le \eta_{q,1}(t),\forall t\in [0,1]$ (defined in Theorem \ref{thm:bound-1}) if \eqref{eq:condition3} holds. 
If we can further show   $\eta_{q,1}(t)\le\Phi_{q}^{\mathrm{sym}}(\frac{1+\rho}{2}),\forall t\in [0,1]$, 
or equivalently, the condition in \eqref{eq:-27},
then by Theorem \ref{thm:bound-1}, \eqref{eq:-28} holds for all
balanced Boolean functions $f$.  Since $\phi_{q}$ is negative,
the second statement follows obviously. 
\end{proof}
In the following,  we analytically verify \eqref{eq:condition4} and numerically verify \eqref{eq:-29-1} for $q=1.36$.
For $q=1.36$, we observe that for all $\rho\in[0,1]$,
$
\epsilon_{q}^{*}(\rho)\ge \epsilon_{0}:=\epsilon_{q}^{*}(1)=0.188324...
$
(see Fig. \ref{fig:The-region-of-1}), and for all $\rho\in[0,1]$,
$
\omega_{q}^{*}(\rho)\le\omega_{0}:=\max_{\beta\in[0,\frac{1}{2}-\epsilon_{q}^{*}(1)]}\omega(\beta)=\omega(\beta_{0})=0.193026...
$
where $\beta_{0}=0.175661...$ Note that the objective function in \eqref{eq:condition4} is a  quadratic polynomial. We hence easily prove that the inequality in  \eqref{eq:condition4} holds in this case. On the other hand, using Matlab, we numerically verify
that  for $q=1.36$, 
\begin{equation}
\max_{\rho\in[0.914,1]}\max_{\substack{t\in[0,1]}
}\chi_{q}(\rho,t)=-0.00169063...<0,\label{eq:max-1}
\end{equation}
with the maximum attained at $\rho=1$ and $t=0.858686...$  
Hence, for all $0.914\le\rho\le1$, dictator functions maximize the
$q$-stability with $q=1.36$  over all balanced Boolean functions.
 Combining it with Theorem  \ref{thm:CK2},
we conclude that for all $\rho\in[0,1]$, dictator functions maximize
the $q$-stability with $q=1.36$ (and hence for all $q\in[1.36,2)$)
over all balanced Boolean functions.  

We employ a computer-assisted approach to rigorously verify the numerical observation stated above. The proof of the following theorem follows an argument nearly identical to that of Theorem \ref{thm:CK}; the complete derivation  is provided in Appendix \ref{sec:Proof-of-Theorem-1-1-1}.

\begin{theorem}
\label{thm:CK-1} The optimal value of \eqref{eq:max-1} for $\Phi=\Phi_{q}^{\mathrm{sym}}$
with $q=1.36$ is negative.  In other words, for any $\rho\in[0,1]$
and $q\in[1.36,2)$, it holds that $\mathbf{Stab}_{q}^{\mathrm{sym}}[f]\le\mathbf{Stab}_{q}^{\mathrm{sym}}[f_{\mathrm{d}}]$
for all balanced Boolean functions $f$. 
\end{theorem}
\begin{remark}
It is already known \cite{yu2023phi,mossel2005coin} that for $q\in[2,5]$,
dictator functions maximize the symmetric $q$-stability over all
balanced Boolean functions. So, combining it with this theorem, we
know that for $q\in[1.36,5]$, dictator functions maximize the symmetric
$q$-stability.
\end{remark}
Li and M\'edard conjectured that dictator functions maximize the symmetric
$q$-stability for all $q\in(1,2)$. The theorem above implies that
the Li--M\'edard conjecture is true for all $q\in[1.36,2)$. Furthermore, it was observed that dictator functions maximizing the $\hat{q}$-stability
for some $\hat{q}\in(0,2)$ implies dictator functions maximizing the $q$-stability
for all $q\in[\hat{q},2)$ \cite{barnes2020courtade}. So, Theorem  \ref{thm:CK2}
implies that the Li--M\'edard conjecture is also true for $q\in(1,1.36)$
and $\rho\in[0,0.914]$. 
\begin{theorem}[Restatement of Theorem \ref{thm:Li-Medard-2}]
\label{thm:Li-Medard}The (symmetrized) Li--M\'edard conjecture is
true for all $q\in[1.36,2)$ and $\rho\in[0,1]$, and also true for all $q\in(1,1.36)$ and $\rho\in[0,0.914]$.
\end{theorem}
The Li--M\'edard conjecture is still open for the case of $q\in(1,1.36)$
and $\rho\in(0.914,1)$. The bound here for $q\in[1.36,2)$ is obtained
from Corollary \ref{cor:If-for-a-1}. In fact, for $q\in(1,1.36)$,
 we numerically evaluate our bound in Corollary \ref{cor:If-for-a-1}
and find a region of $(q,\rho)$ for which dictator functions maximize
the $q$-stability (i.e., for which the Li--M\'edard conjecture is
true). We plot this region in Fig. \ref{fig:rho_q}. 

\subsection{$q$-Stability with $0<q<1$}

It is natural to further investigate whether dictator functions maximize
$q$-stability for  $0<q<1$. In this subsection, we connect   $q$-stability of a 
Boolean function with the edge-boundary of its support. Via this connection, we  observe that 
dictator functions are not extremizers anymore for  $0<q<1/2$, but we conjecture that 
they are indeed extremizers for $1/2\le q<1$.

Let $\mathbf{X}$ be a random vector
uniformly distributed on the discrete cube $\{\pm1\}^{n}$. Let $\mathbf{Y}\in\{\pm1\}^{n}$
be the random vector obtained by independently changing the sign of
each component of $\mathbf{X}$ with the same probability $\frac{1-\rho}{2}$.
Let $A$ be the support of $f$. Then, 
\begin{align*}
T_{\rho}f(\mathbf{x}) & =\mathbb{P}(\mathbf{Y}\in A|\mathbf{X}=\mathbf{x})  =\sum_{\mathbf{y}\in A}\left(\frac{1+\rho}{2}\right)^{n-d(\mathbf{y},\mathbf{x})}\left(\frac{1-\rho}{2}\right)^{d(\mathbf{y},\mathbf{x})},
\end{align*}
where $d(\mathbf{y},\mathbf{x})=|\{i\in[n]:y_{i}\neq x_{i}\}|$ is
the Hamming distance between $\mathbf{x}$ and $\mathbf{y}$. 
\begin{lemma}
\label{lem:Denote} Let $f:\{\pm1\}^{n}\to\{0,1\}$ be any Boolean
function with support $A$. Then, the following hold. 
\begin{enumerate}
\item Denote $t=(\frac{1-\rho}{2})^{q}$, i.e., $\rho=1-2t^{1/q}$. It holds
that for $0<q<1$, 
\begin{align*}
\partial_{t}\mathbb{E}[(T_{\rho}f)^{q}]\big|_{t=0} & =\mathbb{E}[h_{A^{c}}^{q}],
\end{align*}
where for $B\subset\{\pm1\}^{n}$, $h_{B}:\{\pm1\}^{n}\to\mathbb{R}$
is defined by $h_{B}(\mathbf{x})$ equaling the number of edges joining
$\mathbf{x}$ with a vertex in $B^{c}$ if $\mathbf{x}\in B$, and
$h_{B}(\mathbf{x})=0$ if $\mathbf{x}\in B^{c}$. 
\item Denote $t=\rho^{2}$, i.e., $\rho=\sqrt{t}$. It holds that for $0<q<1$,
\begin{align}
\partial_{t}\mathbb{E}[(T_{\rho}f)^{q}]\big|_{t=0} & =\frac{q(q-1)}{2}\alpha^{q-2}\mathbf{W}_{1}[f],\label{eq:-31}
\end{align}
where $\alpha=\mathbb{E}[f]$ is the mean of $f$ and $\mathbf{W}_{1}[f]=\sum_{|S|=1}\hat{f}_{S}^{2}$
is the first-level Fourier weight of $f$ with $\hat{f}_{S},S\subseteq[n]$
denoting the Fourier coefficients of $f$. 
\end{enumerate}
\end{lemma}
\begin{remark}
In contrast, for $q>1$, it is known \cite{li2020boolean} that 
\begin{align*}
\partial_{\rho}\mathbb{E}[(T_{\rho}f)^{q}]\big|_{\rho=1} & =q\mathbb{E}[h_{A}],
\end{align*}
where $\mathbb{E}[h_{A}]=\mathbb{E}[h_{A^{c}}]$ is the classic edge-boundary
of $A$ (or $A^{c}$) normalized by $2^{n}$. Note that $\mathbb{E}[h_{A^{c}}^{q}]$
is a variant of edge-boundary of $A^{c}$ which reduces to the edge-boundary
for $q=1$ and reduces to the vertex-boundary for $q=0$. This quantity
has attracted considerable attention from the probability community
and discrete analysis community \cite{talagrand1993isoperimetry,bobkov1999discrete,kahn2020isoperimetric,beltran2023sharp,durcik2024sharp}.
\end{remark}
\begin{proof}
Statement 1: Observe that 
\begin{align*}
\partial_{t}\mathbb{E}[(T_{\rho}f)^{q}] & =\frac{q}{2^{n}}\sum_{\mathbf{x}\in\{\pm1\}^{n}}\left(T_{\rho}f(\mathbf{x})\right)^{q-1}\partial_{t}T_{\rho}f(\mathbf{x})\\
 & =\frac{q}{2^{n}}\sum_{\mathbf{x}\in\{\pm1\}^{n}}\left(T_{\rho}f(\mathbf{x})\right)^{q-1}\left(\partial_{\rho}T_{\rho}f(\mathbf{x})\big|_{\rho=1-2t^{1/q}}\right)\left(\frac{-2}{q}t^{1/q-1}\right)\\
 & =\frac{-1}{2^{n-1}}\sum_{\mathbf{x}\in\{\pm1\}^{n}}\left(T_{\rho}f(\mathbf{x})\right)^{q-1}\left(\partial_{\rho}T_{\rho}f(\mathbf{x})\big|_{\rho=1-2t^{1/q}}\right)t^{1/q-1}\\
 & =\frac{-1}{2^{n-1}}\sum_{\mathbf{x}\in\{\pm1\}^{n}}\left(\frac{T_{\rho}f(\mathbf{x})}{t^{1/q}}\right)^{q-1}\left(\partial_{\rho}T_{\rho}f(\mathbf{x})\big|_{\rho=1-2t^{1/q}}\right)\\
 & \to\frac{-1}{2^{n-1}}\sum_{\mathbf{x}\in\{\pm1\}^{n}}\lim_{\rho\to1}\left(\frac{2T_{\rho}f(\mathbf{x})}{1-\rho}\right)^{q-1}\left(\partial_{\rho}T_{\rho}f(\mathbf{x})\big|_{\rho=1}\right)\\
 & =\frac{-1}{2^{n-1}}\sum_{\mathbf{x}\in\{\pm1\}^{n}}\left(\partial_{\rho}T_{\rho}f(\mathbf{x})\big|_{\rho=1}\right)\times\begin{cases}
0, & \mathbf{x}\in A\\
\left(-2\partial_{\rho}T_{\rho}f(\mathbf{x})\big|_{\rho=1}\right)^{q-1}, & \mathbf{x}\in A^{c}
\end{cases}\\
 & =\frac{2^{q-1}}{2^{n-1}}\sum_{\mathbf{x}\in A^{c}}\left(-\partial_{\rho}T_{\rho}f(\mathbf{x})\big|_{\rho=1}\right)^{q}.
\end{align*}
For $\mathbf{x}\in A^{c}$, 
\begin{align*}
\partial_{\rho}T_{\rho}f(\mathbf{x}) & =\frac{1}{2}\sum_{\mathbf{y}\in A}(n-d(\mathbf{y},\mathbf{x}))\left(\frac{1+\rho}{2}\right)^{n-d(\mathbf{y},\mathbf{x})-1}\left(\frac{1-\rho}{2}\right)^{d(\mathbf{y},\mathbf{x})}\\
 & \qquad-d(\mathbf{y},\mathbf{x})\left(\frac{1+\rho}{2}\right)^{n-d(\mathbf{y},\mathbf{x})}\left(\frac{1-\rho}{2}\right)^{d(\mathbf{y},\mathbf{x})-1}\\
 & \to-\frac{1}{2}|\{\mathbf{y}\in A:d(\mathbf{y},\mathbf{x})=1\}|  =-\frac{1}{2}h_{A^{c}}(\mathbf{x}),
\end{align*}
where note that $d(\mathbf{y},\mathbf{x})\ge1$ for $\mathbf{y}\in A$
and $\mathbf{x}\in A^{c}$.  Therefore, 
\begin{align*}
\partial_{t}\mathbb{E}[(T_{\rho}f)^{q}] & \to\frac{2^{q-1}}{2^{n-1}}\sum_{\mathbf{x}\in A^{c}}\left(\frac{1}{2}h_{A^{c}}(\mathbf{x})\right)^{q}  =\frac{1}{2^{n}}\sum_{\mathbf{x}\in A^{c}}h_{A^{c}}^{q}(\mathbf{x}) =\mathbb{E}[h_{A^{c}}^{q}].
\end{align*}

Statement 2: Using Fourier expansion $T_{\rho}f=\sum_{S\subset[n]}\rho^{|S|}\hat{f}_{S}\chi_{S}$
where $\chi_{S}$'s are Fourier basis and $\hat{f}_{S}$'s are Fourier
coefficients of $f$, one can easily obtain \eqref{eq:-31}.

\end{proof}
By comparing the $q$-stability for dictator functions and majority
functions, it has been found that for $0<q<1/2$, 
majority functions (i.e., 
$f_{\mathrm{maj}}(\mathbf{x}) = 1\{ \sum_{i=1}^n x_i \ge 0\}$ for odd $n$)
have larger $q$-stability than dictator functions for $\rho$ sufficiently
close to $1$ \cite{beltran2023sharp}. But for $1/2\le q<1$, dictator
functions seem to have the largest $q$-stability for all $0<\rho<1$.
So, we make the following conjectures. 
\begin{conjecture}
\label{conj:Dictator-functions-maximize}Dictator functions maximize
the asymmetric $\frac{1}{2}$-stability over all balanced Boolean
functions. That is, dictator functions minimize $\mathbb{E}[\sqrt{T_{\rho}f}]$
over all balanced Boolean functions $f$. 
\end{conjecture}
A weaker version of this conjecture is the following one with the
asymmetric $\frac{1}{2}$-stability replaced by the symmetric one. 
\begin{conjecture}
\label{conj:Dictator-functions-maximize-1}Dictator functions maximize
the symmetric $\frac{1}{2}$-stability over all balanced Boolean functions.
That is, dictator functions minimize $\mathbb{E}[\sqrt{T_{\rho}f}+\sqrt{1-T_{\rho}f}]$
over all balanced Boolean functions $f$. 
\end{conjecture}
These two conjectures  are  equivalent to that dictator functions maximize
the (asymmetric or symmetric) $q$-stability for all $q\in[1/2,2)$, since the optimality of dictator functions
in maximizing the $\hat{q}$-stability for some $\hat{q}\in(0,2)$ implies
the optimality of dictator functions  for all $q\in[\hat{q},2)$
\cite{barnes2020courtade}. Hence, the conjectures above are stronger than 
the Courtade--Kumar and Li--M\'edard conjectures. 

If we let $\rho=1-2t^{1/q}$ and let $t\downarrow0$, then by Lemma
\ref{lem:Denote}, Conjecture \ref{conj:Dictator-functions-maximize}
implies that for any set $A$ of size $2^{n-1}$, 
\begin{equation}
\mathbb{E}[\sqrt{h_{A}}]\ge\mathbb{E}[\sqrt{h_{C}}]=1/2,\label{eq:-32}
\end{equation}
where $C$ is a subcube of the same size. This is an important conjecture
in discrete analysis \cite[Conjecture 1.3]{durcik2024sharp}. The
quantity $\mathbb{E}[\sqrt{h_{A}}]$ was first considered by Talagrand
in \cite{talagrand1993isoperimetry}, and he provided an isoperimetric
inequality weaker than the conjectured one in \eqref{eq:-32}. The conjectured inequality in \eqref{eq:-32} can be seen as a  sharp
version of Talagrand's isoperimetric inequality. Currently, it is known in \cite{durcik2024sharp} that 
$
\mathbb{E}[\sqrt{h_{A}}]\ge0.4985, 
$
and for $q\ge0.50057$, 
$\mathbb{E}[h_{A}^{q}]\ge\mathbb{E}[h_{C}^{q}]=1/2.$

A similar  $\Phi$-stability, $\mathbb{E}[\sqrt{1-(T_{\rho}g)^{2}}]$,
was considered by Anantharam et al. \cite{anantharam2017conjecture},
where $g=2f-1$ is a Boolean function taking values in $\{\pm1\}$.
This $\Phi$-stability is essentially equivalent to $\mathbb{E}[\sqrt{(T_{\rho}f)(1-T_{\rho}f)}]$
for a Boolean function $f$ taking values in $\{0,1\}$. They conjectured
that this kind of $\Phi$-stability is also maximized by dictator
functions over all balanced Boolean functions, and found that their
conjecture can recover the inequality in \eqref{eq:-32} as well.
 In fact, it is easy to see that Conjecture \ref{conj:Dictator-functions-maximize-1}
is stronger than Anantharam et al.'s conjecture, since 
\begin{align*}
\mathbb{E}[\sqrt{T_{\rho}f}+\sqrt{1-T_{\rho}f}] & =\mathbb{E}[\sqrt{1+2\sqrt{(T_{\rho}f)(1-T_{\rho}f)}}] \le\sqrt{1+2\mathbb{E}[\sqrt{(T_{\rho}f)(1-T_{\rho}f)}]}.
\end{align*}

Although our proof strategy used for $q\ge1$ can be applied to  the
case $1/2\le q<1$, we numerically compute the bound induced by our
strategy and found that  the obtained bound is not very nice. So,
we do not present it here. 

The symmetric $q$-stability problem can be also seen as a problem
of locating roots of $q\mapsto N_{q}(f):=\mathbf{Stab}_{q}^{\mathrm{sym}}[f]-\mathbf{Stab}_{q}^{\mathrm{sym}}[f_{\mathrm{d}}]$,
where $f$ is a balanced Boolean function but not a dictator, and
$f_{\mathrm{d}}$ is a dictator. Barnes and \"Ozg\"ur \cite{barnes2020courtade}
observed that given $\rho\in(0,1)$, $N_{q}(f)$ has at most three
roots: $0$, $q_{1}(f)$, and $q_{2}(f)$, where $q_{1}(f)<2<q_{2}(f)$.
They also observed that the symmetric $q$-stability is maximized
by dictator functions if and only if $q_{1}(f)\le q\le q_{2}(f)$
for all balanced Boolean functions $f$. It is interesting to find
a sharp upper bound on $q_{1}(f)$ and a sharp lower bound on $q_{2}(f)$.
Theorem \ref{thm:CK} implies $q_{1}(f)\le1$  for any $\rho\in(0,0.914]$,
and Theorem \ref{thm:CK-1} implies $q_{1}(f)\le1.36$  for any $\rho\in(0.914,1)$.
Conjecture \ref{conj:Dictator-functions-maximize-1} implies $q_{1}(f)\le1/2$
 for any $\rho\in(0,1)$. If Conjecture \ref{conj:Dictator-functions-maximize-1} is true, then the bound $1/2$
would be sharp uniformly for all $\rho\in(0,1)$. 

\section{Extensions}

\subsection{Extensions to Arbitrary Boolean Functions }

We next extend our bounds on the $\Phi$-stability to arbitrary  (not necessarily balanced) Boolean
functions.  For a Boolean function $f$, denote $A=f^{-1}(1)$, i.e.,
the support of $f$. Here $\mu (A)$ is not necessarily $1/2$.  
Denote the subcubes $C_{S\mapsto a}:=\{\mathbf{x}:x_{S}=a\},S\subseteq[n],a\in\{\pm1\}^{|S|}$.
We aim to investigate the following question: given $1\le k\le n$ and the distances $|A \Delta C_{S \mapsto a}|$ between $A$ and each $k$-codimensional subcube $C_{S \mapsto a}$ (where $|S| = k$ and $a \in \{\pm1\}^k$), what can be said about the $\Phi$-stability of $f$?

For a Boolean function $f$, define 
\begin{align*}
p_{S\mapsto a}(f) & :=\mu(A\cap C_{S\mapsto a})  =\mu(f1_{C_{S\mapsto a}}) =\frac{1}{2}\left(\alpha+2^{-|S|}-\|f-1_{C_{S\mapsto a}}\|_{1}\right),
\end{align*}
where $\alpha=\mathbb{E}f$. Then, the   distance  between $A$ and  $C_{S \mapsto a}$ (or the distance between $f$ and  $1_{C_{S\mapsto a}}$) is determined by $p_{S\mapsto a}(f)$.
In terms of Fourier analysis, $p_{S\mapsto a}(f)$
can be also expressed as 
\[
p_{S\mapsto a}(f)=2^{-|S|}\sum_{T\subseteq S}a_{T}\hat{f}_{T},
\]
where $a_{T}:=\Pi_{i\in T}a_{i}$. Denote 
\[
\vec{p}_{S}(f)=(p_{S\mapsto a}(f))_{a\in\{\pm1\}^{|S|}}.
\]
We now extend Theorem \ref{thm:generalbound-2} to this general setting. 

\begin{theorem}[Bound on $\Phi$-Stability]
\label{thm:generalbound-2-1} Let $\Phi:[0,1]\to\mathbb{R}$ be a
continuous convex function and let $1\le k\le n$. For any Boolean function $f$, it holds
that 
\begin{equation}
\mathbb{E}_{\mu}[\Phi(T_{\rho}f)]\le\min_{S\subseteq[n]:|S|=k}\Gamma(\vec{p}_{S}(f)),\label{eq:-2-3}
\end{equation}
where for $\vec{\epsilon}=(\epsilon_{a})_{a\in\{\pm1\}^{k}}\in[0,1]^{2^{k}}$,
\[
\Gamma(\vec{\epsilon}):=2^{-k}\sum_{a\in\{\pm1\}^{k}}\int_{0}^{1}\Phi\left(\sum_{b\in\{\pm1\}^{k}}\left(\frac{1+\rho}{2}\right)^{n-d(a,b)}\left(\frac{1-\rho}{2}\right)^{d(a,b)}\theta_{\epsilon_{b}}(\beta)\right)\d\beta,
\]
with  $d(a,b):=|\{i:a_i \neq b_i\}|$ denoting the Hamming distance between $a$ and $b$. 
\end{theorem}
\begin{proof}
We first consider the case that $S=[n-k+1:n]$. Note that any Boolean
function $f$ can be written as 
\begin{equation}
f(\mathbf{x})=\sum_{a\in\{\pm1\}^{k}}1\{x_{[n-k+1:n]}=a\}f_{a}(x_{[n-k]}),\label{eq:f-1}
\end{equation}
where $f_{a}(x_{[n-k]}):=f(x_{[n-k]},a)$ is a restriction of $f$, and obviously, $\mu(f_{a})=p_{S\mapsto a}(f)$. We   can write 
\begin{equation}
T_{\rho}f(x_{[n-k]},a)=T_{\rho}^{(n-k)}g_{a}(x_{[n-k]}),\label{eq:Tf-1}
\end{equation}
where 
\[
g_{a}(x_{[n-k]}):=\sum_{b\in\{\pm1\}^{k}}\left(\frac{1+\rho}{2}\right)^{n-d(a,b)}\left(\frac{1-\rho}{2}\right)^{d(a,b)}f_{b}(x_{[n-k]}).
\]
The $\Phi$-stability of $f$ is then 
\begin{align*}
\mathbb{E}_{\mu}[\Phi(T_{\rho}f)] & =2^{-k}\sum_{a\in\{\pm1\}^{k}}\mathbb{E}\Bigl[T_{\rho}^{(n-1)}g_{a}(X_{[n-k]})\Bigr].
\end{align*}

By symmetry, one can obtain a similar formula for other $S$. By Corollary
\ref{cor:generalbound-1}, we obtain Theorem \ref{thm:generalbound-2-1}.
\end{proof}

Similarly to the balanced case, we can also apply hypercontractivity inequalities to 
bound the $q$-stability for non-balanced Boolean functions.
A subset $A\subseteq\{\pm1\}^{n}$ is called a lexicographic set if
it is the initial segment of $\{\pm1\}^{n}$ labeled in the lexicographic
ordering. For example, 
$$A=\{(-1,-1,-1),(-1,-1,1),(-1,1,-1),(-1,1,1)\}$$
is the lexicographic set of $\{\pm1\}^{3}$ with $4$ elements. We
call a Boolean function lexicographic if its support is a lexicographic
set. For a function $f$ and $S\subseteq[n]$, denote $f_{S\mapsto a}(\mathbf{y}):=f(\mathbf{x}),\,\mathbf{y}\in\{\pm1\}^{|S^{c}|}$
as  a restriction of $f$ on $x_{S}=a$, where $\mathbf{x}$ is a vector
given by $x_{S}=a,x_{S^{c}}=\mathbf{y}$. For a Boolean function
$f$ and $S\subseteq[n]$, the $T$-rearrangement $f_{S}^{*}$ of
$f$ is a Boolean function $g$ such that $\mathbb{E}g_{S\mapsto a}=\mathbb{E}f_{S\mapsto a}$
and $g_{S\mapsto a}$ is lexicographic for any $a\in\{\pm1\}^{|S|}$.
We obtain a generalization of Theorem \ref{thm:qstability} for arbitrary
Boolean functions. 
\begin{theorem}[Bound on  $q$-Stability]
\label{thm:generalbound-2-1-2} For  any Boolean function $f$,
it holds that for any $q>0$ and $S\subseteq[n]$, 
\begin{align}
\mathbf{Stab}_{q}[f] & \le\frac{\mathbb{E}\left[\mathbb{E}\left[(T_{\rho}^{S}f_{S}^{*}(\mathbf{X}))^{p}|X_{S}\right]^{q/p}\right]-\mathbb{E}[f]}{q-1},\label{eq:-2-2-4}
\end{align}
where $p=1+(q-1)\rho^{2}$, and  $T_{\rho}^{S}$ is the noise
operator acting on coordinates indexed by  $S$, i.e., $T_{\rho}^{S}g(\mathbf{X})=\mathbb{E}[g(\mathbf{Z})|X_{S}]$
with $\mathbf{Z}$  denoting a random vector such that $Z_{S^{c}}=X_{S^{c}}$
and $Z_{S}\in\{\pm1\}^{|S|}$  being $\rho$-correlated with $X_{S}$.
Here, the inequality in \eqref{eq:-2-2-4} for $q=1$ is understood
as the limit of this inequality as $q\to1$. 
\end{theorem}
 
\begin{proof}
Denote $\mu_{S}$ and $\mu_{S^{c}}$ as the marginal distributions
of $\mu$ respectively on the components indexed by  $S$ and $S^{c}$. By
comparing the concentration spectra of $(T_{\rho}^{S}f)_{S\mapsto a}$
and $(T_{\rho}^{S}f_{S}^{*})_{S\mapsto a}$ for $a\in\{\pm1\}^{|S|}$,
it is easy to see that $(T_{\rho}^{S}f_{S}^{*})_{S\mapsto a}$ has
a larger concentration spectrum (and also a larger smooth version)
than $(T_{\rho}^{S}f)_{S\mapsto a}$. So, $(T_{\rho}^{S}f)_{S\mapsto a}\prec(T_{\rho}^{S}f_{S}^{*})_{S\mapsto a}$
for any $a\in\{\pm1\}^{|S|}$. We then have that for $q>1$, 
\begin{align}
\mathbb{E}[(T_{\rho}f)^{q}] & =\mathbb{E}_{\mu_{S}}[\mathbb{E}_{\mu_{S^{c}}}[(T_{\rho}^{S^{c}}(T_{\rho}^{S}f))^{q}]]\\
 & \le\mathbb{E}_{\mu_{S}}\left[\mathbb{E}_{\mu_{S^{c}}}\left[(T_{\rho}^{S}f)^{p}\right]^{q/p}\right] \le\mathbb{E}_{\mu_{S}}\left[\mathbb{E}_{\mu_{S^{c}}}\left[(T_{\rho}^{S}f_{S}^{*})^{p}\right]^{q/p}\right],
\end{align}
where the first inequality follows by the hypercontractivity inequality
acting on $(T_{\rho}^{S^{c}},\mu_{S^{c}})$, and the last inequality
follows by the majorization inequality in Proposition \ref{prop:majorization-1}.
The case $0<q<1$ follows similarly. 
\end{proof}

\subsection{Extensions to Ornstein--Uhlenbeck Operators }

 We now consider the Ornstein--Uhlenbeck operator  given by 
\[
U_{\rho}f(\mathbf{x})=\mathbb{E}[f(\rho\mathbf{x}+\sqrt{1-\rho^{2}}\mathbf{Z})],\,\forall\mathbf{x}\in\mathbb{R}^{n},
\]
where $\mathbf{Z}$ is a standard $n$-dimensional Gaussian random
vector. The invariant distribution for this semigroup is the standard
Gaussian measure $\gamma$ on $\mathbb{R}^{n}$. In this subsection,
given $\rho,\alpha\in[0,1]$, we use $\mathcal{F}_{\alpha}=\mathcal{F}_{\alpha,\rho}$
to denote  the set of $U_{\rho}f$ for Boolean functions $f$ with
mean $\gamma(f)=\alpha$. 

 For $(\alpha,\beta)\in[0,1]^{2}$, the maximal noise stability at
$(\alpha,\beta)$ for the Ornstein--Uhlenbeck operator is 
\[
\mathbf{S}_{n}(\alpha,\beta)=\max_{\textrm{Bool }\phi,f:\,\gamma(f)=\alpha,\gamma(\phi)=\beta}\int\phi U_{\rho}f\d\gamma.
\]
It is known that the maximization is attained by half spaces, which yields that 
\[
\mathbf{S}_{n}(\alpha,\beta)=\Theta(\alpha,\beta):=\Psi_{\rho}(\Psi^{-1}(\alpha),\Psi^{-1}(\beta)),
\]
where $\Psi$ is the CDF of the standard univariate Gaussian, and
$\Psi_{\rho}$ is the CDF of the bivariate Gaussians with correlation
coefficient $\rho$. 
\begin{proposition}[Relation Between Majorization and Noise Stability for OU Operators]
\label{prop:-Let--1} Then, for any $\rho,\alpha\in[0,1]$, it holds
that $\mathcal{F}_{\alpha}\prec\theta_{\alpha}$, where for $\beta\in[0,1]$,
\begin{align}
\theta_{\alpha}(\beta) & :=\frac{\partial_{\beta}\Theta(\alpha,\beta)}{\partial\beta} =\Psi\left(\frac{\Psi^{-1}(\alpha)-\rho\Psi^{-1}(\beta)}{\sqrt{1-\rho^{2}}}\right).\label{eq:-7}
\end{align}
\end{proposition}
Note that by choosing $f(x)=1\{x\le\Psi^{-1}(\alpha)\}$ as the indicator
of a half-space of measure $\alpha$, one can find that for any $t\in[0,1]$,
$
U_{\rho}f(\Psi^{-1}(t))=\theta_{\alpha}(t),
$
which implies 
\begin{align}
\mathfrak{C}_{U_{\rho}f}(\beta)=\mathfrak{D}_{U_{\rho}f}(\beta) & =\sup_{\textrm{Bool }\phi:\,\gamma(\phi)=\beta}\int\phi U_{\rho}f\d\gamma\label{eq:concentration-2-1-1-1}\\
 & =\sup_{\textrm{Bool }\phi:\,\int_{0}^{1}\phi(\Psi^{-1}(t))\d t=\beta}\int_{0}^{1}\phi(\Psi^{-1}(t))U_{\rho}f(\Psi^{-1}(t))\d t \label{eq:change}\\
 & =\sup_{\textrm{Bool }\kappa:\,\int_{0}^{1}\kappa(t)\d t=\beta}\int_{0}^{1}\kappa(t)\theta_{\alpha}(t)\d t\label{eq:change2}\\
 & =\mathfrak{D}_{\theta_{\alpha}}(\beta)=\mathfrak{C}_{\theta_{\alpha}}(\beta),
\end{align}
where \eqref{eq:change} follows by a change of variables and in \eqref{eq:change2},  $\kappa :=\phi\circ \Psi^{-1}$. That is, $U_{\rho}f$ defined on $\mathbb{R}$ with the Gaussian
measure and $\theta_{\alpha}$  defined on $[0,1]$ with the uniform
measure have the same concentration spectrum, i.e., $U_{\rho}f\prec\theta_{\alpha}$
and $U_{\rho}f\succ\theta_{\alpha}$. So, the majorization relation given in Proposition \ref{prop:-Let--1}
is sharp and cannot be improved any more. 

Using Proposition \ref{prop:-Let--1}, one can recover the classic
Borell's result on $\Phi$-stability. More precisely, the maximal
noise stability problem is equivalent to the maximal $\Phi$-stability
problem. 

\begin{theorem}[Borell's $\Phi$-Stability Theorem \cite{borell1985geometric}]
\label{thm:generalbound-3} Let $\Phi:[0,1]\to\mathbb{R}$ be a 
convex function. For any Boolean function $f$ with mean $\gamma(f)=\alpha$,
it holds that
\[
\mathbb{E}_{\gamma}[\Phi(U_{\rho}f)]\le\int_{0}^{1}\Phi(\theta_{\alpha}(\beta))\d\beta=\mathbb{E}_{Z\sim\gamma}\Bigg[\Phi\bigg(\Psi\left(\frac{\Psi^{-1}(\alpha)-\rho Z}{\sqrt{1-\rho^{2}}}\right)\bigg)\Bigg],
\]
where $\theta_{\alpha}$ is given in \eqref{eq:-7}. Moreover, this
upper bound is attained by $f(x)=1\{x\le\Psi^{-1}(\alpha)\}$. 
\end{theorem}

\appendix

\section{\label{sec:Proof-of-Proposition-Theta}Proof of Proposition \ref{prop:Theta}}

In fact, the small set expansion theorem is obtained by substituting Boolean
functions of  into the hypercontractivity inequality (stated in Lemma \ref{lem:HC}). That is, it
was shown on \cite[p.280]{O'Donnell14analysisof} that 
\begin{equation}
\mathbf{S}_{n}(\alpha,\beta)\le\min_{p,q\ge1:(p-1)(q-1)=\rho^{2}}\alpha^{1/p}\beta^{1/q}=\Theta_{0}(\alpha,\beta).\label{eq:-333}
\end{equation}

Due to the equality in \eqref{eq:-333}, to prove Proposition \ref{prop:Theta},
it suffices to show that for all $\alpha,\beta\in(0,1)$ satisfying
$\alpha+\beta\le1$ and all $p,q\ge1$ satisfying $\frac{1}{p}+\frac{1}{q}\ge1$,
\[
\alpha^{1/q}\beta^{1/p}\le\alpha+\beta-1+(1-\alpha)^{1/p}(1-\beta)^{1/q}.
\]
Changing variables $s=\frac{1}{p},t=\frac{1}{q}$, $u=1-\alpha,v=\beta$,
it is equivalent to that for all $u,v\in(0,1)$ satisfying $u\ge v$
and all $s,t\in[0,1]$ satisfying $s+t\ge1$, 
\[
f(s):=u^{s}(1-v)^{t}-v^{s}(1-u)^{t}+v-u\ge0.
\]
We now prove this inequality. 

Observe that 
\[
f'(s)=u^{s}(1-v)^{t}\ln u-v^{s}(1-u)^{t}\ln v,
\]
which has the same sign as 
\[
g(s):=\left(\frac{u}{v}\right)^{s}(1-v)^{t}\ln u-(1-u)^{t}\ln v.
\]
Obviously, $g$ is nonincreasing. Thus, $g$  (and also $f'$) is
either always positive, always negative, or positive first and then
negative as $s\in[1-t,1]$ increases, which means that $f$ is either
always monotone or increasing first and then decreasing as $s\in[1-t,1]$
increases. As a consequence, to prove $f(s)\ge0$ for all $s\in[1-t,1]$,
it suffices to verify  that $f(1-t)\ge0$ and $f(1)\ge0$. That is, for
all $t\in[0,1]$,
\begin{align*}
\varphi(t) & :=f(1)=u(1-v)^{t}-v(1-u)^{t}+v-u\ge0,\\
\psi(t) & :=f(1-t)=u^{1-t}(1-v)^{t}-v^{1-t}(1-u)^{t}+v-u\ge0.
\end{align*}
By the same idea used above, one can find that both $\varphi$ and
$\psi$ are either always monotone or increasing first and then decreasing
as $t\in[0,1]$ increases. Thus, we only need verify that $\varphi(0)\ge0,\varphi(1)\ge0,\psi(0)\ge0,\psi(1)\ge0$.
By definition, these are obviously true.

\section{\label{sec:Proof-of-Theorem}Proof of Theorem \ref{thm:asymptotics}}

From \eqref{eq:theta-1} (more clearly, from Fig. \ref{fig:alphabeta}), one can observe that given $\rho\in (0,1)$, for sufficiently
large $s$ (i.e., for sufficiently small $\epsilon>0$),  
\begin{align}
(\theta_{\epsilon}(\beta),\theta_{1-\epsilon}(\beta)) & =\begin{cases}
(0,0), & 0\le t\le\rho s, \; \hat{t}\ge s/\rho\\
(0,1-g(s,\hat{t})), & 0\le t\le\rho s, \; \rho s\le\hat{t}\le s/\rho\\
(0,1), & 0\le t\le\rho s,\; 0\le\hat{t}\le\rho s\\
(g(s,t),1), & \rho s\le t\le s/\rho, \; 0\le\hat{t}\le\rho s\\
(1,1), & t\ge s/\rho, \; 0\le\hat{t}\le\rho s
\end{cases},\label{eq:-8}
\end{align}
where  $\epsilon=\exp(-\frac{s^{2}}{2})=1-\exp(-\frac{\hat{s}^{2}}{2})$,
$\beta=\exp(-\frac{t^{2}}{2})=1-\exp(-\frac{\hat{t}^{2}}{2})$, and 
\begin{align}
g(s,t):=\frac{t-\rho s}{(1-\rho^{2})t}\exp\left(-\frac{(s-\rho t)^{2}}{2(1-\rho^{2})}\right). \label{eq:g}
\end{align}
Note that here $s,\hat{s}$ are defined in the   way given in  Proposition  \ref{prop:majorizationNO} for $\alpha=\epsilon$, and hence  $s$ plays the role of $\hat{s}$ for $\alpha=1-\epsilon$.  
On the other hand, 
\begin{align*}
0\le\hat{t}\le\rho s & \Longleftrightarrow t\ge t_{1}(s),\\
\hat{t}\ge s/\rho & \Longleftrightarrow t\le t_{0}(s),
\end{align*}
where for sufficiently large $s$, 
\begin{align*}
t_{1}(s) & :=\sqrt{-2\ln(1-e^{-\frac{s^{2}}{2}\rho^{2}})}\sim\sqrt{2}e^{-\frac{s^{2}}{4}\rho^{2}},\\
t_{0}(s) & :=\sqrt{-2\ln(1-e^{-\frac{s^{2}}{2}\rho^{-2}})}\sim\sqrt{2}e^{-\frac{s^{2}}{4}\rho^{-2}}.
\end{align*}
Hence, \eqref{eq:-8} can be rewritten as 
\begin{align*}
(\theta_{\epsilon}(\beta),\theta_{1-\epsilon}(\beta)) & =\begin{cases}
(0,0), & t\le t_{0}(s)\\
(0,1-g(s,\hat{t})), & t_{0}(s)\le t\le t_{1}(s)\\
(0,1), & t_{1}(s)\le t\le\rho s\\
(g(s,t),1), & \rho s\le t\le s/\rho\\
(1,1), & t\ge s/\rho
\end{cases}.
\end{align*}
Substituting this into Theorem \ref{thm:generalbound-2} yields
\begin{align}
\Gamma(\epsilon) & =\left(\Phi(0)+\Phi(1)\right)e^{-\frac{s^{2}}{2\rho^{2}}}+\frac{1}{2}\left(\Phi(\frac{1+\rho}{2})+\Phi(\frac{1-\rho}{2})\right)\left(1-2e^{-\frac{s^{2}\rho^{2}}{2}}\right)+A(s)+B(s),\label{eq:-46}
\end{align}
where $\epsilon=e^{-\frac{s^{2}}{2}}$, 
\begin{align}
A(s) & :=\frac{1}{2}\int_{t_{0}(s)}^{t_{1}(s)}\Big\{\Phi\left(\frac{1+\rho}{2}(1-g(s,\hat{t}))\right)+\Phi\left(\frac{1-\rho}{2}(1-g(s,\hat{t}))\right)\Big\} te^{-\frac{t^{2}}{2}}\d t\nonumber \\
 & =\frac{1}{2}\int_{\rho s}^{s/\rho}\Big\{\Phi\left(\frac{1+\rho}{2}(1-g(s,\hat{t}))\right)+\Phi\left(\frac{1-\rho}{2}(1-g(s,\hat{t}))\right)\Big\}\hat{t}e^{-\frac{\hat{t}^{2}}{2}}\d\hat{t},\label{eq:-44}
\end{align}
and 
\begin{align}
B(s) & :=\frac{1}{2}\int_{\rho s}^{s/\rho}\Big\{\Phi\left(\frac{1+\rho}{2}+\frac{1-\rho}{2}g(s,t)\right)+\Phi\left(\frac{1-\rho}{2}+\frac{1+\rho}{2}g(s,t)\right)\Big\} te^{-\frac{t^{2}}{2}}\d t.\label{eq:-45}
\end{align}

We divide the integral in \eqref{eq:-44}  into two parts $\int_{\rho s}^{cs}+\int_{cs}^{s/\rho}$
where $c\in(\rho,1/\rho)$ is a fixed number chosen  close to $1/\rho$.
Applying Taylor's theorem, the integral $\int_{\rho s}^{cs}$  is
equal to 
\begin{align*}
 & \frac{1}{2}\int_{\rho s}^{cs}\Big\{\Phi\left(\frac{1+\rho}{2}\right)-\Phi'\left(\frac{1+\rho}{2}\right)\frac{1+\rho}{2}g(s,\hat{t})+\frac{1}{2}\Phi''\left(\xi_{1,\hat{t}}\right)\left(\frac{1+\rho}{2}g(s,\hat{t})\right)^{2}\\
 & \qquad+\Phi\left(\frac{1-\rho}{2}\right)-\Phi'\left(\frac{1-\rho}{2}\right)\frac{1-\rho}{2}g(s,\hat{t})+\frac{1}{2}\Phi''\left(\xi_{2,\hat{t}}\right)\left(\frac{1-\rho}{2}g(s,\hat{t})\right)^{2}\Big\}\hat{t}e^{-\frac{\hat{t}^{2}}{2}}\d\hat{t}
\end{align*}
where $\xi_{1,\hat{t}}\in[\frac{1+\rho}{2}(1-g(s,\hat{t})),\frac{1+\rho}{2}]$
and $\xi_{1,\hat{t}}\in[\frac{1-\rho}{2}(1-g(s,\hat{t})),\frac{1-\rho}{2}]$.
Since by the definition of $g$   in \eqref{eq:g}, $g(s,\hat{t})\to0$ as $s\to\infty$ uniformly for all $\hat{t}\in[\rho s,cs]$,
the integral above is then equal to 
\begin{align*}
 & \frac{1}{2}\left(\Phi(\frac{1+\rho}{2})+\Phi(\frac{1-\rho}{2})\right)\left(e^{-\frac{\rho^{2}s^{2}}{2}}-e^{-\frac{c^{2}s^{2}}{2}}\right)\\
 & \quad-\frac{1}{2}\left(\frac{1+\rho}{2}\Phi'(\frac{1+\rho}{2})+\frac{1-\rho}{2}\Phi'(\frac{1-\rho}{2})\right)\left(e^{-\frac{s^{2}}{2}}-e^{-\left(1+\frac{(c-\rho)^{2}}{1-\rho^{2}}\right)\frac{s^{2}}{2}}\right)\\
 & \quad+\frac{1+o(1)}{2}\left(\left(\frac{1+\rho}{2}\right)^{2}\Phi''(\frac{1+\rho}{2})+\left(\frac{1-\rho}{2}\right)^{2}\Phi''(\frac{1-\rho}{2})\right)\int_{\rho s}^{cs}g(s,\hat{t})^{2}\hat{t}e^{-\frac{\hat{t}^{2}}{2}}\d\hat{t},
\end{align*}
where $o(1)$ is a term vanishing as $s\to\infty$. 

On the other hand, 
\begin{align*}
 & \left|\frac{1}{2}\int_{cs}^{s/\rho}\Big\{\Phi\left(\frac{1+\rho}{2}(1-g(s,\hat{t}))\right)+\Phi\left(\frac{1-\rho}{2}(1-g(s,\hat{t}))\right)\Big\}\hat{t}e^{-\frac{\hat{t}^{2}}{2}}\d\hat{t}\right|\\
 & \le\frac{1}{2}\max_{cs\le t\le s/\rho}\left|\Phi\left(\frac{1+\rho}{2}(1-g(s,t))\right)+\Phi\left(\frac{1-\rho}{2}(1-g(s,t))\right)\right|\int_{cs}^{s/\rho}\hat{t}e^{-\frac{\hat{t}^{2}}{2}}\d\hat{t}\\
 & =O\left(e^{-\frac{c^{2}s^{2}}{2}}-e^{-\frac{s^{2}}{2\rho^{2}}}\right).
\end{align*}
Here, note that $\Phi$ is bounded, since it is continuous on the closed interval $[0,1]$.

Similarly, the integral in \eqref{eq:-45} is divided into two parts
$\int_{\rho s}^{cs}+\int_{cs}^{s/\rho}$, the first part of which
is equal to 
\begin{align*}
 & \frac{1}{2}\left(\Phi(\frac{1+\rho}{2})+\Phi(\frac{1-\rho}{2})\right)\left(e^{-\frac{\rho^{2}s^{2}}{2}}-e^{-\frac{c^{2}s^{2}}{2}}\right)\\
 & \quad+\frac{1}{2}\left(\frac{1-\rho}{2}\Phi'(\frac{1+\rho}{2})+\frac{1+\rho}{2}\Phi'(\frac{1-\rho}{2})\right)\left(e^{-\frac{s^{2}}{2}}-e^{-\left(1+\frac{(c-\rho)^{2}}{1-\rho^{2}}\right)\frac{s^{2}}{2}}\right)\\
 & \quad+\frac{1+o(1)}{2}\left(\left(\frac{1-\rho}{2}\right)^{2}\Phi''(\frac{1+\rho}{2})+\left(\frac{1+\rho}{2}\right)^{2}\Phi''(\frac{1-\rho}{2})\right)\int_{\rho s}^{cs}g(s,t)^{2}te^{-\frac{t^{2}}{2}}\d t,
\end{align*}
and the second of which is bounded as follows:
\begin{align*}
 & \left|\frac{1}{2}\int_{cs}^{s/\rho}\Big\{\Phi\left(\frac{1+\rho}{2}+\frac{1-\rho}{2}g(s,t)\right)+\Phi\left(\frac{1-\rho}{2}+\frac{1+\rho}{2}g(s,t)\right)\Big\} te^{-\frac{t^{2}}{2}}\d t\right|\\
 & \le\frac{1}{2}\max_{cs\le t\le s/\rho}\left|\Phi\left(\frac{1+\rho}{2}+\frac{1-\rho}{2}g(s,t)\right)+\Phi\left(\frac{1-\rho}{2}+\frac{1+\rho}{2}g(s,t)\right)\right|\int_{cs}^{s/\rho}te^{-\frac{t^{2}}{2}}\d t\\
 & =O\left(e^{-\frac{c^{2}s^{2}}{2}}-e^{-\frac{s^{2}}{2\rho^{2}}}\right).
\end{align*}

Substituting all of these into \eqref{eq:-46} and noticing that $c^{2}<1+\frac{(c-\rho)^{2}}{1-\rho^{2}}<1/\rho^{2}$,
we obtain that 
\begin{align*}
\Gamma(\epsilon) & =\frac{1}{2}\left(\Phi(\frac{1+\rho}{2})+\Phi(\frac{1-\rho}{2})\right)-\frac{\rho}{2}\left(\Phi'(\frac{1+\rho}{2})-\Phi'(\frac{1-\rho}{2})\right)e^{-\frac{s^{2}}{2}}\\
 & \qquad+(1+o(1))\frac{1+\rho^{2}}{4}\left(\Phi''(\frac{1+\rho}{2})+\Phi''(\frac{1-\rho}{2})\right)\int_{\rho s}^{cs}g(s,t)^{2}te^{-\frac{t^{2}}{2}}\d t+O(e^{-\frac{c^{2}s^{2}}{2}}).
\end{align*}
It seems not easy to estimate the exact asymptotics for $\int_{\rho s}^{cs}g(s,t)^{2}te^{-\frac{t^{2}}{2}}\d t$.
But we can estimate it by using 
\[
\int_{\rho s}^{cs}g(s,t)^{2}t^{2}e^{-\frac{t^{2}}{2}}\d t\sim\sqrt{2\pi}\frac{\rho^{2}}{1+\rho^{2}}\sqrt{\frac{1-\rho^{2}}{1+\rho^{2}}}s^{2}e^{-\frac{s^{2}}{1+\rho^{2}}}  \; (\text{as } s\to \infty),
\]
which induces that 
\[
\int_{\rho s}^{cs}g(s,t)^{2}te^{-\frac{t^{2}}{2}}\d t\sim C\sqrt{2\pi}\frac{\rho^{2}}{1+\rho^{2}}\sqrt{\frac{1-\rho^{2}}{1+\rho^{2}}}se^{-\frac{s^{2}}{1+\rho^{2}}},
\]
where $C\in[\rho,1/\rho]$. Therefore, 
\begin{align*}
\Gamma(\epsilon) & =\frac{1}{2}\left(\Phi(\frac{1+\rho}{2})+\Phi(\frac{1-\rho}{2})\right)-\frac{\rho}{2}\left(\Phi'(\frac{1+\rho}{2})-\Phi'(\frac{1-\rho}{2})\right)e^{-\frac{s^{2}}{2}}\\
 & \qquad+(1+o(1))C\sqrt{2\pi}\frac{\rho^{2}}{4}\sqrt{\frac{1-\rho^{2}}{1+\rho^{2}}}\left(\Phi''(\frac{1+\rho}{2})+\Phi''(\frac{1-\rho}{2})\right)se^{-\frac{s^{2}}{1+\rho^{2}}}+O(e^{-\frac{c^{2}s^{2}}{2}}).
\end{align*}

\section{\label{sec:Proof-of-Lemma}Proof of Lemma \ref{lem:unique}}

Statement 1: Denote $f(\epsilon):=h\left(\frac{1-\rho}{2}+\rho\epsilon\right)-\left(1+\frac{2\rho^{2}\epsilon}{1-\rho^{2}}\right)h\left(\frac{1-\rho}{2}\right).$
Then, $f$ is convex given any $\rho\in(0,1)$. Observe that $f'(0)=\frac{-2\rho}{1-\rho^{2}}g(\rho)$,
where $g(\rho):=\frac{1+\rho}{2}\log\frac{1+\rho}{2}-\frac{1-\rho}{2}\log\frac{1-\rho}{2}$.
Moreover, $g''(\rho)=\frac{-2\rho}{1-\rho^{2}}<0$. So, $g$ is strictly
concave. Observe that $g(0)=g(1)=0$. So, $g>0$ on $(0,1)$, i.e.,
$f'(0)<0$ for any $\rho\in(0,1)$. On the other hand, by hypercontractivity
inequalities, $f(\frac{1}{2})\ge0$. Combining all these points yields
that there is a unique solution in $(0,1/2]$ to the equation in \eqref{eq:-33}. 

Statement 2: Denote $f(\epsilon):=\left(\epsilon+\left(\frac{1+\rho}{2}\right)^{p}(1-2\epsilon)\right)^{q/p}+\left(\epsilon+\left(\frac{1-\rho}{2}\right)^{p}(1-2\epsilon)\right)^{q/p}-\left(\frac{1+\rho}{2}\right)^{q}-\left(\frac{1-\rho}{2}\right)^{q}.$
Since for $q>1$, $f$ is convex and $f(\frac{1}{2})\ge0$ (by hypercontractivity
inequalities), $f(\epsilon)=0$ has a unique solution in $(0,1/2]$
if and only if $f'(0)<0$. 

Observe that $f'(0)=\frac{q}{p}g(q)$, where 
\begin{align}
	g(q) & :=\left(\frac{1+\rho}{2}\right)^{(q-1)(1-\rho^{2})}+\left(\frac{1-\rho}{2}\right)^{(q-1)(1-\rho^{2})}-2\left(\frac{1+\rho}{2}\right)^{q}-2\left(\frac{1-\rho}{2}\right)^{q}\nonumber \\
	& =\left(\frac{1+\rho}{2}\right)^{-(1-\rho^{2})}e^{q(1-\rho^{2})\ln\frac{1+\rho}{2}}+\left(\frac{1-\rho}{2}\right)^{-(1-\rho^{2})}e^{q(1-\rho^{2})\ln\frac{1-\rho}{2}} \nonumber \\ 
	& \qquad-2e^{q\ln\frac{1+\rho}{2}}-2e^{q\ln\frac{1-\rho}{2}}.\label{eq:-38}
\end{align}
It suffices to prove $g(q)<0$ for all $q\in(1,2]$. In the following,
we use Laguerre's lemma to prove this point. 

Observe that for $\rho\in(0,1)$, the exponents in \eqref{eq:-38}
satisfy that 
\begin{equation}
	\ln\frac{1-\rho}{2}<(1-\rho^{2})\ln\frac{1-\rho}{2}<\ln\frac{1+\rho}{2}<(1-\rho^{2})\ln\frac{1+\rho}{2},\label{eq:-37}
\end{equation}
which follows by the following arguments. The first and the last inequalities
are obvious. For the second inequality, we compute the derivatives
of $r(\rho):=(1-\rho^{2})\ln\frac{1-\rho}{2}-\ln\frac{1+\rho}{2}$
as follows: 
\begin{align*}
	r'''(\rho) & =\frac{2\left(1+7\rho+2\rho^{2}-\rho^{3}-\rho^{4}\right)}{(1-\rho)^{2}(1+\rho)^{3}}>0,\\
	r''(0) & =\ln4>0.
\end{align*}
So, $r''(\rho)>0$ for all $\rho\in(0,1)$, i.e., $r$ is strictly
convex. Verifying that $r(0)=r(1)=0$, we have that $r(\rho)<0$ for
all $\rho\in(0,1)$. That is, the second inequality in \eqref{eq:-37}
holds. 

The corresponding coefficients of the terms \eqref{eq:-38} in the
same order as \eqref{eq:-37} are $-2$, $\left(\frac{1-\rho}{2}\right)^{-(1-\rho^{2})}$,
$-2$, and $\left(\frac{1+\rho}{2}\right)^{-(1-\rho^{2})}$, whose
signs are $-,+,-,+$. By Laguerre's lemma \cite{laguerre1883sur}\cite[Lemma 1]{barnes2020courtade},
there are at most three zeros for $g(q)$.  Obviously, one of the
zeros is $1$. Moreover, for sufficiently large $q$, the term $\left(\frac{1+\rho}{2}\right)^{-(1-\rho^{2})}e^{q(1-\rho^{2})\ln\frac{1+\rho}{2}}$
dominates others, and hence, $g(q)>0$ for sufficiently large $q$.
Similarly, the term $-2e^{q\ln\frac{1-\rho}{2}}$ dominates others
as $q$ approaches $-\infty$, and hence, $g(q)<0$ as $q$ approaches
$-\infty$. Observe that 
\begin{align*}
	g'(1) & =(1-\rho^{2})\ln\frac{1+\rho}{2}+(1-\rho^{2})\ln\frac{1-\rho}{2}-2\left(\frac{1+\rho}{2}\right)\ln\frac{1+\rho}{2}-2\left(\frac{1-\rho}{2}\right)\ln\frac{1-\rho}{2}\\
	& =-2\rho\left(\frac{1+\rho}{2}\log\frac{1+\rho}{2}-\frac{1-\rho}{2}\log\frac{1-\rho}{2}\right)  <0,
\end{align*}
where the last line follows by the proof for the case $q=1$ given
above. Combining all these points, we know that $g$ has exactly three
roots $1,q_{1},q_{2}$ such that $q_{1}<1<q_{2}$. Moreover, $g(q)<0$
for all $q\in(1,q_{2})$ and $g(q)>0$ for all $q>q_{2}$. So, to
prove $g(q)<0$ for all $q\in(1,2]$, it suffices to prove that $g(2)<0$. 

Observe that 
\begin{align*}
	g(2) & =\left(\frac{1+\rho}{2}\right)^{1-\rho^{2}}+\left(\frac{1-\rho}{2}\right)^{1-\rho^{2}}-2\left(\frac{1+\rho}{2}\right)^{2}-2\left(\frac{1-\rho}{2}\right)^{2}\\
	& =\left(\frac{1+\rho}{2}\right)^{1-\rho^{2}}+\left(\frac{1-\rho}{2}\right)^{1-\rho^{2}}-1-\rho^{2}\\
	& \le2\left(\frac{1}{2}\right)^{1-\rho^{2}}-1-\rho^{2} =2^{\rho^{2}}-1-\rho^{2}  <0,
\end{align*}
where the first inequality follows by Jensen's inequality (since $t\mapsto t^{1-\rho^{2}}$
is concave), and the last inequality follows by the convexity of the
exponential function. Therefore, we have $g(q)<0$ for all $q\in(1,2]$,
i.e., $f'(0)<0$ for all $q\in(1,2]$ and $\rho\in(0,1)$. 

\section{\label{sec:Proof-of-Lemma-1}Proof of Lemma \ref{lem:monotonicity}}

Consider the function 
\begin{equation}
	f(\rho,t):=(1-\rho^{2})h\left(\frac{1-t\rho}{2}\right)-(1-t\rho^{2})h\left(\frac{1-\rho}{2}\right).\label{eq:-33-1}
\end{equation}
So, by definition of $\epsilon_{1,\rho}^{*}$, it holds that $f(\rho,1-2\epsilon_{1,\rho}^{*})=0$.
We then compute 
\[
\partial_{\rho}^{4}f(\rho,t)=\frac{2(1-t)a(\rho,t)}{(1-\rho)^{3}(1+\rho)^{3}(1-\rho t)^{3}(1+\rho t)^{3}}
\]
where 
\begin{align*}
	a(\rho,t) & :=\left(\rho^{4}-3\rho^{2}+6\right)\rho^{6}t^{6}+\left(-7\rho^{6}+10\rho^{4}-3\rho^{2}\right)t^{5}-\left(3\rho^{6}-2\rho^{4}+8\rho^{2}+3\right)\rho^{2}t^{4}\\
	& \qquad+\left(3\rho^{6}+\rho^{4}-3\rho^{2}-1\right)t^{3}+\left(6\rho^{6}-8\rho^{4}+15\rho^{2}-1\right)t^{2} \\
	&\qquad+\left(\rho^{4}-6\rho^{2}+5\right)t-3\rho^{2}-1.
\end{align*}
We then compute the partial derivative of $a$:
\[
\partial_{\rho}a(\rho,t)=2\rho b(\rho,t)
\]
where 
\begin{align*}
	b(\rho,t) & :=\rho^{4}\left(5\rho^{4}-12\rho^{2}+18\right)t^{6}+\left(-21\rho^{4}+20\rho^{2}-3\right)t^{5}-\left(12\rho^{6}-6\rho^{4}+16\rho^{2}+3\right)t^{4}\\
	& \qquad+\left(9\rho^{4}+2\rho^{2}-3\right)t^{3}+\left(18\rho^{4}-16\rho^{2}+15\right)t^{2}+2\left(\rho^{2}-3\right)t-3.
\end{align*}
Compute the partial derivative of $b$:
\[
\partial_{\rho}b(\rho,t)=4\rho tc(\rho,t)
\]
where
\begin{align*}
	c(\rho,t) & :=2\rho^{2}\left(5\rho^{4}-9\rho^{2}+9\right)t^{5}+\left(10-21\rho^{2}\right)t^{4}-2\left(9\rho^{4}-3\rho^{2}+4\right)t^{3}\\
	& \qquad+\left(9\rho^{2}+1\right)t^{2}+2\left(9\rho^{2}-4\right)t+1.
\end{align*}
Compute the partial derivative of $c$:
\[
\partial_{\rho}c(\rho,t)=6\rho td(\rho,t)
\]
where
\[
d(\rho,t)=2\left(5\rho^{4}-6\rho^{2}+3\right)t^{4}-7t^{3}+\left(2-12\rho^{2}\right)t^{2}+3t+6.
\]
Compute the partial derivative of $d$:
\begin{equation}
	\partial_{\rho}d(\rho,t)=8\rho t^{2}\left(\left(5\rho^{2}-3\right)t^{2}-3\right)\le8\rho t^{2}\left(2t^{2}-3\right)\le0.\label{eq:-43}
\end{equation}

Equation \eqref{eq:-43} implies that $d$ is decreasing in $\rho$.
Combining this with 
\[
d(0,t)=6t^{4}-7t^{3}+2t^{2}+3t+6\ge0
\]
yields that $d$ is either positive or first-positive-then-negative.
So, $c$ is either increasing or first-increasing-then-decreasing.
Moreover, we observe that 
\[
c(1,t)=10t^{5}-11t^{4}-20t^{3}+10t^{2}+10t+1\ge0.
\]
So, $c$ is either positive or first-negative-then-positive, which
implies that $b$ is either increasing or first-decreasing-then-increasing.
We also observe that 
\[
b(0,t)=-3\left(t^{5}+t^{4}+t^{3}-5t^{2}+2t+1\right)\le0.
\]
So, $b$ is either negative or first-negative-then-positive, which
implies that $a$ is either decreasing or first-decreasing-then-increasing.
We also observe that $a(1,t)=-4(1-t^{2})^{3}\le0.$ So, $a$ is either
negative or first-positive-then-negative, and hence, so is $\partial_{\rho}^{4}f(\rho,t)$.
This implies that $\partial_{\rho}^{3}f(\rho,t)$ is either decreasing
or first-increasing-then-decreasing. Observe that $\partial_{\rho}^{3}f(0,t)=0$
and $\partial_{\rho}^{3}f(1,t)=-\infty.$ So, $\partial_{\rho}^{3}f(\rho,t)$
is either negative or first-positive-then-negative, which implies
that $\partial_{\rho}^{2}f(\rho,t)$ is either decreasing or first-increasing-then-decreasing.
Observe that $\partial_{\rho}^{2}f(1,t)=-\infty.$ So, $\partial_{\rho}^{2}f(\rho,t)$
is negative, or first-positive-then-negative, or first-negative-next-positive-then-negative,
which implies that $\partial_{\rho}f(\rho,t)$ is decreasing, or first-increasing-then-decreasing,
or first-decreasing-next-increasing-then-decreasing.  Observe that
$\partial_{\rho}f(0,t)=0$ and $\partial_{\rho}f(1,t)=-\infty.$ So,
$\partial_{\rho}f(\rho,t)$ is negative, or first-positive-then-negative,
or first-negative-next-positive-then-negative, which implies that
$f(\rho,t)$ is decreasing, or first-increasing-then-decreasing, or
first-decreasing-next-increasing-then-decreasing. 

We consider $t=t_{1}$ such that $f(\rho_{1},t_{1})=0$ for some $\rho_{1}\in(0,1)$,
which implies that there is at least one $\rho$ such that $f(\rho,t_{1})=0$.
Moreover, observe that $f(0,t)=f(1,t)=0.$ So, $f(\rho,t)$ is first-decreasing-next-increasing-then-decreasing,
which implies that $\rho=\rho_{1}$ is the unique solution on $(0,1)$
to $f(\rho,t_{1})=0$. Moreover,  for any $\rho_{2}>\rho_{1}$, it
holds that $f(\rho_{2},t_{1})>0$. Let $t_{2}$ be the solution to
$f(\rho_{2},t)=0$. 

Since $f(\rho_{2},t)$ is convex in $t$, and moreover, $f(\rho_{2},1)=f(\rho_{2},t_{2})=0$,
we have that $f(\rho_{2},t)>0$ for $t<t_{2}$, and $f(\rho_{2},t)<0$
for $t_{2}<t<1$. So, $f(\rho_{2},t_{1})>0$ implies $t_{1}<t_{2}$,
i.e., $\epsilon_{1,\rho_{1}}^{*}>\epsilon_{1,\rho_{2}}^{*}$ for any
$\rho_{1}<\rho_{2}$. That is, $\epsilon_{1,\rho}^{*}$ is decreasing
in $\rho\in(0,1)$. 

\section{\label{sec:Proof-of-Theorem-1}Proof of Theorem \ref{thm:bound}}

Recall that 
\[
\phi(s) =\frac{s\ln s+(1-s)\ln(1-s)}{s},\,s\in(0,1].
\]
We need the following basic properties of the function $\phi$. 
\begin{lemma}
\label{lem:phi}The following hold. 
\begin{enumerate}
\item $\phi''$ is increasing on $(0,1)$ and negative on $(0,1/2]$. In
particular, $\phi'$ is decreasing on $(0,1/2]$ and positive on $(0,1)$. 
\item $\phi$ is increasing on $(0,1]$, and 
\begin{equation}
\conc\,\phi(t)=\begin{cases}
\phi(t), & t\in(0,1/2]\\
4(t-1)\ln2, & t\in[1/2,1]
\end{cases}.\label{eq:-23}
\end{equation}
\item If $\phi'(t_{1})=\phi'(t_{2})$ for distinct $t_{1},t_{2}\in(0,1)$,
then $t_{1}+t_{2}>1$. 
\end{enumerate}
\end{lemma}
\begin{proof}
Statement 1: Compute derivatives of $\phi$: 
\begin{align*}
\phi'(t) & =-\frac{\ln(1-t)}{t^{2}}, \quad 
\phi''(t)  =\frac{\frac{t}{1-t}+2\ln(1-t)}{t^{3}},\quad 
\phi'''(t)  =\frac{\frac{t(5t-4)}{(1-t)^{2}}-6\ln(1-t)}{t^{4}}.
\end{align*}
It is easy to observe that $\phi'(t)>0$ and $\phi'''(t)>0$ for all
$t\in(0,1)$, which implies that $\phi''$ is increasing. Observe
that $\phi''(1/2)<0$, which implies that $\phi''$ is negative on
$(0,1/2]$. In other words, $\phi'$ is decreasing for $(0,1/2]$.

Statement 2: Since $\phi'(t)>0$ for all $t\in(0,1)$, we have that
$\phi$ is increasing. Furthermore, since $\phi''$ is increasing,
$\phi''$ is first-negative-then-positive. So, $\phi$ is first-concave-then-convex,
and moreover, $\phi$ is concave on $(0,1/2]$. Observe that $\phi'(1/2)=4\ln2=\frac{\phi(1)-\phi(1/2)}{1-1/2}$.
So, the graph of $\phi$ on $[1/2,1]$ is below the tangent line of
$\phi$ at $t=1/2$, which implies \eqref{eq:-23}. 

Statement 3: Denote $t_{1}=t-\delta,t_{2}=t+\delta$. Without loss
of generality, we assume $\delta>0$. Then, we need show that $\phi'(t-\delta)=\phi'(t+\delta)$
for some $\delta$ implies $t>1/2$. We now prove this by contradiction.
That is, we are aim at showing that $t\le1/2$ implies $\phi'(t-\delta)\neq\phi'(t+\delta)$
for all $\delta\in(0,1-t)$. We first assume $t=1/2$. For this case,
\begin{align*}
\phi'(\frac{1}{2}+\delta)-\phi'(\frac{1}{2}-\delta) & =-\frac{\ln(\frac{1}{2}-\delta)}{(\frac{1}{2}+\delta)^{2}}+\frac{\ln(\frac{1}{2}+\delta)}{(\frac{1}{2}-\delta)^{2}} =\frac{f(\delta)}{(\frac{1}{4}-\delta^{2})^{2}},
\end{align*}
with $f(\delta):=(\frac{1}{2}+\delta)^{2}\ln(\frac{1}{2}+\delta)-(\frac{1}{2}-\delta)^{2}\ln(\frac{1}{2}-\delta)$.
Observe that 
\begin{align*}
f'(\delta) & =2(\frac{1}{2}+\delta)\ln(\frac{1}{2}+\delta)+(\frac{1}{2}+\delta)+2(\frac{1}{2}-\delta)\ln(\frac{1}{2}-\delta)+(\frac{1}{2}-\delta) =1-2h(\frac{1}{2}-\delta),\\
f''(\delta) & =2h'(\frac{1}{2}-\delta)>0,\,\forall\delta\in(0,1/2),
\end{align*}
the latter of which implies that $f$ is strictly convex on $(0,1/2)$.
Combined with $\lim_{t\downarrow0}f(t)=\lim_{t\uparrow1/2}f(t)=0$,
this implies $f(\delta)<0$ for all $\delta\in(0,1/2)$. So, $\phi'(\frac{1}{2}-\delta)>\phi'(\frac{1}{2}+\delta)$
for all $\delta\in(0,1-t)$.

We next focus on the case $t<1/2$. Observe that $\phi'$ is decreasing
on $(0,t_{0}]$ and increasing on $[t_{0},1)$ for some $t_{0}\in(1/2,1)$.
So, $\phi'(t-\delta)>\phi'(t+\delta)$ for all $\delta\in(0,t_{0}-t]$,
and $\phi'(t-\delta)>\phi'(\frac{1}{2}-\delta)>\phi'(\frac{1}{2}+\delta)>\phi'(t+\delta)$
for all $\delta\in[t_{0}-t,1-t)$. That is, $\phi'(t-\delta)>\phi'(t+\delta)$
for all $\delta\in(0,1-t)$. 
\end{proof}
Denoting $t_{1}=-2\rho z_{1},t_{2}=2\rho z_{2}$, and discarding the
constraints $0\leq p_{1}\leq\frac{1}{4}+\frac{\beta}{2},\,0\leq p_{2}\leq\frac{1}{4}-\frac{\beta}{2}$
in \eqref{eq:-19}, we obtain an upper bound on the second term at RHS of 
\eqref{eq:-3} as follows: 
\begin{equation}
\max_{\beta\in[0,\frac{1}{2}-\epsilon^{*}(\rho)]}\Upsilon_{\rho}(\beta)\le \max_{\beta\in[0,\frac{1}{2}-\epsilon^{*}(\rho)]}\max_{t_{1},t_{2}\in(-1,1):t_{1}+t_{2}\ge0}\frac{1-\rho}{2}\eta(t_{1},t_{2},\beta),\label{eq:-4-1}
\end{equation}
where 
\[
\eta(t_{1},t_{2},\beta):=q_{1}\phi(\frac{1-t_{1}}{2})+q_{2}\phi(\frac{1-t_{2}}{2}),
\]
with 
\begin{align*}
 & q_{1}=\frac{1+\rho-4\rho^{2}\omega(\beta)+2\beta(1+t_{2}-\rho^{2})}{2+t_{1}+t_{2}-2\rho^{2}},\\
 & q_{2}=\frac{1+\rho-4\rho^{2}\omega(\beta)-2\beta(1+t_{1}-\rho^{2})}{2+t_{1}+t_{2}-2\rho^{2}}.
\end{align*}
Here, 
we observe that $\eta(t_{1},t_{2},\beta)\to-\infty$
as $t_{1}\uparrow1$ or $t_{2}\uparrow1$. So, the inner maximum in
\eqref{eq:-4-1} is attained at some $(t_{1},t_{2})$, and is indeed a maximum. 

We first determine
strictly interior stationary points for this maximization. Compute
the following derivatives: 
\begin{align*}
\partial_{t_{1}}\eta(t_{1},t_{2},\beta) & =-\frac{1}{2}q_{1}\phi'(\frac{1-t_{1}}{2})-q_{1}\frac{\phi(\frac{1-t_{1}}{2})+\phi(\frac{1-t_{2}}{2})}{2+t_{1}+t_{2}-2\rho^{2}},\\
\partial_{t_{2}}\eta(t_{1},t_{2},\beta) & =-\frac{1}{2}q_{2}\phi'(\frac{1-t_{2}}{2})-q_{2}\frac{\phi(\frac{1-t_{1}}{2})+\phi(\frac{1-t_{2}}{2})}{2+t_{1}+t_{2}-2\rho^{2}}.
\end{align*}
The strictly interior stationary points satisfy that  $\partial_{t_{1}}\eta(t_{1},t_{2},\beta)=\partial_{t_{2}}\eta(t_{1},t_{2},\beta)=0$,
which implies 
\begin{align*}
\phi'(\frac{1-t_{1}}{2}) & =\phi'(\frac{1-t_{2}}{2}),
\end{align*}
by noting that $q_{1},q_{2}>0$. Since $t_{1}+t_{2}\ge0$, by Lemma
\ref{lem:phi}, we have $t_{1}=t_{2}$. So, the inner maximum in \eqref{eq:-4-1}
reduces to
\begin{align}
 & \max_{t\in[0,1)}\frac{\left(1-\rho\right)\left(1+\rho-4\rho^{2}\omega(\beta)\right)}{2\left(1+t-\rho^{2}\right)}\phi(\frac{1-t}{2}).\label{eq:-22}
\end{align}

We next consider points on the boundary $t_{1}+t_{2}=0$ with
$t_{1}\in(-1,1)$.   For this case, the inner
maximum in \eqref{eq:-4-1} reduces to 
\begin{equation}
\max_{t\in(-1,1)}\frac{1-\rho}{2}\eta(-t,t,\beta),\label{eq:-4-1-3}
\end{equation}
where 
\begin{align*}
\eta(-t,t,\beta) & =q_{1}\phi(\frac{1+t}{2})+q_{2}\phi(\frac{1-t}{2})\\
 & = \frac{1}{2-2\rho^{2}}\Bigl[\left(A-B\right)\phi(\frac{1+t}{2}) +\left(A+B\right)\phi(\frac{1-t}{2})\Bigr],
\end{align*}
with  $A=1+\rho-4\rho^{2}\omega(\beta)$ and $B=2\beta\rho^{2}$. Here the identity $(1+t)\phi(\frac{1+t}{2})=(1-t)\phi(\frac{1-t}{2})$ is used. 


We first consider the case $t\ge 0$. By Lemma \ref{lem:phi}, 
\begin{align*}
\phi(\frac{1+t}{2})+\phi(\frac{1-t}{2}) & \le2\conc\,\phi(\frac{1}{2})=2\phi(\frac{1}{2}),\\
\phi(\frac{1+t}{2})-\phi(\frac{1-t}{2}) & \ge0.
\end{align*}
So,  the optimization in \eqref{eq:-4-1-3} is upper bounded by 
\[
\frac{1-\rho}{2}\frac{1+\rho-4\rho^{2}\omega(\beta)}{1-\rho^{2}}\phi(\frac{1}{2}),
\]
which corresponds to the value of the objective function in \eqref{eq:-22}
at $t=0$, and hence, is no larger than the optimal value of the optimization
in \eqref{eq:-22}. 

We next consider the case $t<0$. 
For this case,  we rewrite 
\begin{align*}
\eta(-t,t,\beta) & = \frac{A+Bt}{(1-\rho^{2})(1-t)} \phi(\frac{1+t}{2}) = \frac{A-Bs}{(1-\rho^{2})(1+s)} \phi(\frac{1-s}{2}),
\end{align*}
where $s=-t \in (0,1)$.
So,  the optimization in \eqref{eq:-4-1-3} is upper bounded by 
\[
\frac{1+\rho-4\rho^{2}\omega(\beta)-2\beta\rho^{2} s}{2(1+\rho)(1+s)}\phi(\frac{1-s}{2}).
\]

Summarizing all the above, the whole optimization in \eqref{eq:-4-1}
is equal to 
\begin{align*}
 & \max_{\beta\in[0,\frac{1}{2}-\epsilon^{*}(\rho)]}\max_{t\in[0,1)}\max \left\{\frac{\left(1-\rho\right)\left(1+\rho-4\rho^{2}\omega(\beta)\right)}{2\left(1+t-\rho^{2}\right)},\frac{1+\rho-4\rho^{2}\omega(\beta)-2\beta\rho^{2} t}{2(1+\rho)(1+t)}\right\} \phi(\frac{1-t}{2}) \nonumber \\
 & =\max_{t\in[0,1)}\max
 \left\{\eta_1(t),\eta_2(t)\right\}. 
\end{align*}

\section{\label{sec:Proof-of-Theorem-1-1}Proof of Theorem \ref{thm:CK}}

We need the following basic properties of $\omega$, which can be
verified easily. 
\begin{lemma}
\label{lem:omega} For $0\leq\beta\leq0.418$, it holds that $\varphi_{\mathrm{LP}}(\frac{1}{2}-\beta)\le\varphi_{\mathrm{C}}(\frac{1}{2}-\beta)$,
which implies 
\begin{align}
\omega(\beta) & =\min\Big\{\beta^{2}+\varphi_{\mathrm{LP}}(\frac{1}{2}-\beta),\,\frac{\left(1+\sqrt{1+4(\pi-\sqrt{2\pi})\beta}\right)^{2}}{8\pi}\Big\}\nonumber \\
 & =\begin{cases}
\frac{\left(1+\sqrt{1+4(\pi-\sqrt{2\pi})\beta}\right)^{2}}{8\pi}, & 0\leq\beta\leq\beta_{0}\\
\beta^{2}-\frac{\beta}{2}+\frac{1}{4}, & \beta_{0}<\beta\leq\frac{1}{4}\\
\beta^{2}+2(\frac{1}{2}-\beta)^{3/2}-2(\frac{1}{2}-\beta)^{2}, & \frac{1}{4}<\beta\leq0.418
\end{cases}\label{eq:-16}
\end{align}
with $\beta_{0}=0.175661...$ denoting the unique value $\beta\in(0,1/2)$
such that the expressions in the first and second clauses above are
equal. It further can be verified that for $0\leq\beta\leq0.345$,
$\omega(\beta)\le\omega(\beta_{0})=0.193026...$ 
\end{lemma}
We also need the following basic properties of $\epsilon^{*}(\rho)$,
$\omega^{*}(\rho)$, and $t_{\rho}$. 
\begin{lemma}
\label{lem:The-following-hold.}The following hold. 
\begin{enumerate}
\item For all $\rho\in[0.46,0.914]$, it holds that $\epsilon^{*}(\rho)\ge\epsilon^{*}(0.914)=0.195055...$
\item For all $\rho\in[0.46,0.914]$, it holds that $\omega^{*}(\rho)\le\omega^{*}(0.914) \le\omega(\beta_{0})=0.193026...$ 
\item For all $\rho\in[0.46,0.914]$, it holds that $t_{\rho}\in[0,0.75]$,
where $t_{\rho}\in(0,1)$ is the unique solution to the equation in
\eqref{eq:-30}. 
\end{enumerate}
\end{lemma}

\begin{proof}
Statement 1 follows by the fact that $\epsilon^{*}(\rho)$ is decreasing
in $\rho\in[0.46,0.914]$. The latter was shown in Lemma \ref{lem:monotonicity}.

Statement 2 can be verified by using Statement 1 and Lemma \ref{lem:omega}.
We now prove Statement 3. Observe that 
\begin{align*}
\phi'(\frac{1-t_{\rho}}{2}) & =\frac{-2\phi(\frac{1-\rho}{2})}{1+\rho-4\rho^{2}\omega^{*}(\rho)} \le-2\phi(\frac{1-\rho}{2})  \le-2\phi(\frac{1-0.914}{2})=-2\phi(0.043),
\end{align*}
which, combined with the monotonicity of $\phi'$ (given in Lemma
\ref{lem:phi}), implies $t_{\rho}\in[0,0.75]$ for all $\rho\in[0.46,0.914]$. 
\end{proof}

We next prove Theorem \ref{thm:CK} by using the proof strategy described
below Theorem  \ref{thm:CK2}. We
first swap two maximizations in \eqref{eq:max} and then apply Lemma
\ref{lem:The-following-hold.} to obtain that 
\begin{equation}
\max_{\rho\in[0.46,0.914]}\max_{t\in[0,1]}\chi(\rho,t)=\max_{t\in[0,0.75]}\max_{\rho\in[0.46,0.914]}\chi(\rho,t).\label{eq:-4-2-3}
\end{equation}
Observe that 
\begin{align*}
\partial_{\rho}\chi(\rho,t) & =\left(1-8\rho\omega_{0}\right)\phi(\frac{1-t}{2})+2\rho\phi(\frac{1-\rho}{2})+\frac{1}{2}\left(1+t-\rho^{2}\right)\phi'(\frac{1-\rho}{2}).
\end{align*}
We bound this derivative as follows: 
\begin{align*}
|\partial_{\rho}\chi(\rho,t)| & \le|1-8\rho\omega_{0}|\left|\phi(\frac{1-t}{2})\right|+\max\left\{ 2\rho\left|\phi(\frac{1-\rho}{2})\right|,\frac{1}{2}\left(1+t-\rho^{2}\right)\left|\phi'(\frac{1-\rho}{2})\right|\right\} \\
 & \le(8\times0.914\times\omega_{0}-1)\left|\phi(\frac{1-0.75}{2})\right|\\
 & \qquad+\max\left\{ 2\times0.914\left|\phi(\frac{1-0.914}{2})\right|,\frac{1}{2}\left(1+0.75-0.46^{2}\right)\left|\phi'(\frac{1-0.914}{2})\right|\right\} \\
 & \le20=:M.
\end{align*}
We choose $\delta=0.0016$. Then, $\Delta:=\frac{\delta}{M}=0.00008$.
We finally use Matlab to verify $\chi(\rho,t)<-\delta$ numerically
for all $t\in[0,0.75]$ and\footnote{Here, $[a:\Delta:b]$ denotes the set of numbers $a+k\Delta,k\in\mathbb{N}$
in $[a,b]$. } $\rho\in[0.46:\Delta:0.914]$, which, by \eqref{eq:-24}, yields
that $\chi(\hat{\rho},t)<0$ for all $\hat{\rho}\in[0.46,0.914]$.
Note that here, for each $\rho$, to verify $\chi(\rho,t)<-\delta$,
it suffices to verify $\chi(\rho,t_{\rho})<-\delta$ where $t_{\rho}\in(0,1)$
is the unique solution to the equation in \eqref{eq:-30}. The Matlab code for this verification is included within the TeX source files of this arXiv version.

\section{\label{sec:Proof-of-Theorem-1-2}Proof of Theorem \ref{thm:bound-1}}

We need the following basic properties of the function $\phi_{q}$. 
\begin{lemma}
	\label{lem:phi-1}For any $1<q<2$, the following hold. 
	\begin{enumerate}
		\item $\phi_{q}''$ is increasing on $(0,1)$ and negative on $(0,1/2]$.
		In particular, $\phi_{q}'$ is decreasing on $(0,1/2]$ and positive
		on $(0,1)$. 
		\item $\phi_{q}$ is increasing on $(0,1]$, and 
		\begin{equation}
			\conc\,\phi_{q}(t)=\begin{cases}
				\phi_{q}(t), & t\in(0,1/2]\\
				4(t-1)\ln_{q}2, & t\in[1/2,1]
			\end{cases}.\label{eq:-23-1}
		\end{equation}
		\item If $\phi_{q}'(t_{1})=\phi_{q}'(t_{2})$ for distinct $t_{1},t_{2}\in(0,1)$,
		then $t_{1}+t_{2}>1$. 
	\end{enumerate}
\end{lemma}
\begin{proof}
	Statement 1: Compute derivatives of $\phi_{q}$: 
	\begin{align*}
		\phi_{q}'(t) & =\frac{1+(q-1)t^{q}-\left(1+(q-1)t\right)(1-t)^{q-1}}{(q-1)t^{2}},\\
		\phi_{q}''(t) & =\frac{(q-2)(q-1)t^{q}+qt((q-3)t+2)(1-t)^{q-2}+2\left((1-t)^{q}-1\right)}{(q-1)t^{3}},\\
		\phi_{q}'''(t) & =\frac{1}{(q-1)t^{4}} \Big( (q-3)(q-2)(q-1)t^{q}\\
		&-\left((q-3)(q-2)(q-1)t^{3} +3(q-3)(q-2)t^{2}+6(q-3)t+6\right)(1-t)^{q-3}+6\Big).
	\end{align*}
	We observe that
	\begin{align*}
		(q-3)(q-2)(q-1)t^{q} & \ge0,\\
		-\left((q-3)(q-2)(q-1)t^{3}+3(q-3)(q-2)t^{2}+6(q-3)t+6\right)(1-t)^{q-3}+6 & \ge0,
	\end{align*}
	where the second inequality can be shown by analyzing derivatives
	of the expression at the LHS. So, we obtain that  $\phi_{q}'''(t)\ge0$
	for all $t\in(0,1)$, which implies that $\phi''$ is increasing.
	Observe that $\phi_{q}''(1/2)=\frac{2^{4-q}((q-1)q+2)-16}{q-1}\le0$,
	which implies that $\phi_{q}''$ is negative on $(0,1/2]$. In other
	words, $\phi_{q}'$ is decreasing for $(0,1/2]$.  Furthermore, observe
	that $\left(1+(q-1)t\right)(1-t)^{q-1}-1-(q-1)t^{q}$ is zero at $t=0$
	and is strictly increasing in $t\in(0,1)$, since its derivative is
	positive. So, we obtain that $\phi_{q}'(t)>0$.
	
	Statement 2: Since $\phi_{q}'(t)>0$ for all $t\in(0,1)$, we have
	that $\phi_{q}$ is increasing. Furthermore, since $\phi_{q}''$ is
	increasing, $\phi_{q}''$ is first-negative-then-positive. So, $\phi_{q}$
	is first-concave-then-convex, and moreover, $\phi_{q}$ is concave
	on $(0,1/2]$. Observe that $\phi_{q}'(1/2)=-\frac{4\left(2^{-(q-1)}-1\right)}{q-1}=\frac{\phi_{q}(1)-\phi_{q}(1/2)}{1-1/2}$.
	So, the graph of $\phi_{q}$ on $[1/2,1]$ is below the tangent line
	of $\phi_{q}$ at $t=1/2$, which implies \eqref{eq:-23-1}. 
	
	Statement 3: Denote $t_{1}=t-\delta,t_{2}=t+\delta$. Without loss
	of generality, we assume $\delta>0$. Then, we need show that $\phi_{q}'(t-\delta)=\phi_{q}'(t+\delta)$
	for some $\delta$ implies $t>1/2$. We now prove this by contradiction.
	That is, we are aim at showing that $t\le1/2$ implies $\phi_{q}'(t-\delta)\neq\phi_{q}'(t+\delta)$
	for all $\delta\in(0,1-t)$. We first assume $t=1/2$. For this case,
	\begin{align*}
		\phi_{q}'(\frac{1}{2}+\delta)-\phi_{q}'(\frac{1}{2}-\delta) & =\frac{8f(\delta)}{\left(1-4\delta^{2}\right)^{2}(q-1)},
	\end{align*}
	with 
	\[
	f(\delta):=(q-2\delta(q-2))\left(\frac{1}{2}+\delta\right)^{q}-(q+2\delta(q-2))\left(\frac{1}{2}-\delta\right)^{q}-4\delta.
	\]
	Observe that 
	\begin{align*}
		f'(\delta) & =q(q+2\delta(q-2))\left(\frac{1}{2}-\delta\right)^{q-1}+q(q-2\delta(q-2))\left(\frac{1}{2}+\delta\right)^{q-1}\\
		& \qquad-2(q-2)\left(\frac{1}{2}-\delta\right)^{q}-2(q-2)\left(\frac{1}{2}+\delta\right)^{q}-4,\\
		f''(\delta) & =q\left(q^{2}-3q+4-2\delta(q-2)(q+1)\right)\left(\frac{1}{2}+\delta\right)^{q-2}\\
		& \qquad-q\left(q^{2}-3q+4+2\delta(q-2)(q+1)\right)\left(\frac{1}{2}-\delta\right)^{q-2},\\
		f'''(\delta) & =(2-q)(q-1)q\left(3-q+2\delta(q+1)\right)\left(\frac{1}{2}-\delta\right)^{q-3}g(\delta),
	\end{align*}
	where 
	\[
	g(\delta):=\left(\frac{\frac{1}{2}-\delta}{\frac{1}{2}+\delta}\right)^{3-q}-\frac{2\delta(q+1)+q-3}{3-q+2\delta(q+1)}.
	\]
	Moreover, 
	\[
	g'(\delta)=-4(3-q)\left((1+2\delta)^{q-4}(1-2\delta)^{2-q}+\frac{q+1}{(3-q+2\delta(q+1))^{2}}\right)<0,
	\]
	which implies that $g$ is decreasing. Combining this with the facts
	that $g(0)=2>0$ and $g(1/2)=\frac{1-q}{2}<0$ yields that $g$ is
	first-positive-then-negative on $(0,1/2)$. Hence, $f'''$ is first-positive-then-negative,
	and hence, $f''$ is first-increasing-then-decreasing. Combining this
	with the facts that $f''(0)=0$ and $f''(1/2)=-\infty$ yields that
	$f''$ is first-positive-then-negative on $(0,1/2)$. Hence, $f'$
	is first-increasing-then-decreasing. Combining this with the facts
	that $f'(0)=2^{2-q}((q-1)q+2)-4\le0$ and $f'(1/2)=0$ yields that
	$f'$ is first-negative-then-positive on $(0,1/2)$. Hence, $f$ is
	first-decreasing-then-increasing. Combining this with the fact that
	$f(0)=f(1/2)=0$ yields that $f$ is negative on $(0,1/2)$. So,
	$\phi_{q}'(\frac{1}{2}-\delta)>\phi_{q}'(\frac{1}{2}+\delta)$ for
	all $\delta\in(0,\frac{1}{2})$.
	
	We next focus on the case $t<1/2$. Observe that $\phi_{q}'$ is decreasing
	on $(0,t_{0}]$ and increasing on $[t_{0},1)$ for some $t_{0}\in(1/2,1)$.
	So, $\phi_{q}'(t-\delta)>\phi_{q}'(t+\delta)$ for all $\delta\in(0,t_{0}-t]$,
	and $\phi_{q}'(t-\delta)>\phi_{q}'(\frac{1}{2}-\delta)>\phi_{q}'(\frac{1}{2}+\delta)>\phi_{q}'(t+\delta)$
	for all $\delta\in[t_{0}-t,1-t)$. That is, $\phi_{q}'(t-\delta)>\phi_{q}'(t+\delta)$
	for all $\delta\in(0,1-t)$. 
\end{proof}
Denoting $t_{1}=-2\rho z_{1},t_{2}=2\rho z_{2}$, and discarding the
constraints $0\leq p_{1}\leq\frac{1}{4}+\frac{\beta}{2},\,0\leq p_{2}\leq\frac{1}{4}-\frac{\beta}{2}$
in \eqref{eq:-19}, we obtain an upper bound on the term at RHS of 
\eqref{eq:-3-2} as follows:  
\begin{equation}
	\max_{\beta\in[0,\frac{1}{2}-\epsilon_{q}^{*}(\rho)]}\Upsilon_{\rho}(\beta) \le \max_{\beta\in[0,\frac{1}{2}-\epsilon_{q}^{*}(\rho)]}\max_{t_{1},t_{2}\in[-1,1]:t_{1}+t_{2}\ge0}\frac{1-\rho}{2}\eta_{q}(t_{1},t_{2},\beta),\label{eq:-4-1-2}
\end{equation}
where 
\[
\eta_{q}(t_{1},t_{2},\beta):=q_{1}\phi_{q}(\frac{1-t_{1}}{2})+q_{2}\phi_{q}(\frac{1-t_{2}}{2}),
\]
with 
\begin{align*}
	& q_{1}=\frac{1+\rho-4\rho^{2}\omega(\beta)+2\beta(1+t_{2}-\rho^{2})}{2+t_{1}+t_{2}-2\rho^{2}},\\
	& q_{2}=\frac{1+\rho-4\rho^{2}\omega(\beta)-2\beta(1+t_{1}-\rho^{2})}{2+t_{1}+t_{2}-2\rho^{2}}.
\end{align*}

Observe that the inner maximum in \eqref{eq:-4-1-2} is attained at
some $(t_{1},t_{2})$. We first determine strictly interior stationary
points for this maximization. Compute the following derivatives: 
\begin{align*}
	\partial_{t_{1}}\eta_{q}(t_{1},t_{2},\beta) & =-\frac{1}{2}q_{1}\phi_{q}'(\frac{1-t_{1}}{2})-q_{1}\frac{\phi_{q}(\frac{1-t_{1}}{2})+\phi_{q}(\frac{1-t_{2}}{2})}{2+t_{1}+t_{2}-2\rho^{2}},\\
	\partial_{t_{2}}\eta_{q}(t_{1},t_{2},\beta) & =-\frac{1}{2}q_{2}\phi_{q}'(\frac{1-t_{2}}{2})-q_{2}\frac{\phi_{q}(\frac{1-t_{1}}{2})+\phi_{q}(\frac{1-t_{2}}{2})}{2+t_{1}+t_{2}-2\rho^{2}}.
\end{align*}
The strictly interior stationary points satisfy that  $\partial_{t_{1}}\eta_{q}(t_{1},t_{2},\beta)=\partial_{t_{2}}\eta_{q}(t_{1},t_{2},\beta)=0$,
which implies 
\begin{align*}
	\phi_{q}'(\frac{1-t_{1}}{2}) & =\phi_{q}'(\frac{1-t_{2}}{2}),
\end{align*}
by noting that $q_{1},q_{2}>0$. Since $t_{1}+t_{2}\ge0$, by Lemma
\ref{lem:phi-1}, we have $t_{1}=t_{2}$. So, the inner maximum in
\eqref{eq:-4-1-2} reduces to
\begin{align}
	& \max_{t\in[0,1]}\frac{\left(1-\rho\right)\left(1+\rho-4\rho^{2}\omega(\beta)\right)}{2\left(1+t-\rho^{2}\right)}\phi_{q}(\frac{1-t}{2}).\label{eq:-22-2}
\end{align}

Observe that $\phi_{q}'(0)=\infty$, which implies 
\[
\partial_{t_{1}}\eta_{q}(t_{1},t_{2},\beta)\Big|_{t_{1}=1}=-q_{1}\left(\frac{1}{2}\phi_{q}'(0)+\frac{\frac{-q}{q-1}+\phi_{q}(\frac{1-t_{2}}{2})}{3+t_{2}-2\rho^{2}}\right)=-\infty.
\]
So, any point on the boundary $t_{1}=1$ can not be a maximizer. Similarly,
any point on the boundary $t_{2}=1$ can not be a maximizer as well. 

We next show that all points on the boundary $t_{1}+t_{2}=0$ with
$t_{1}\in(-1,1)$ can not be maximizers.  For this case, the inner
maximum in \eqref{eq:-4-1-2} reduces to 
\begin{equation}
	\max_{t\in(-1,1)}\frac{1-\rho}{2}\eta_{q}(-t,t,\beta),\label{eq:-4-1-3-1}
\end{equation}
%
where 
\begin{align*}
	\eta_{q}(-t,t,\beta) & =q_{1}\phi_{q}(\frac{1+t}{2})+q_{2}\phi_{q}(\frac{1-t}{2})\\
	& = \frac{1}{2-2\rho^{2}}\Bigl[\left(A-B\right)\phi_{q}(\frac{1+t}{2}) +\left(A+B\right)\phi_{q}(\frac{1-t}{2})\Bigr],
\end{align*}
with  $A=1+\rho-4\rho^{2}\omega(\beta)$ and $B=2\beta\rho^{2}$. Here the identity $(1+t)\phi_{q}(\frac{1+t}{2})=(1-t)\phi_{q}(\frac{1-t}{2})$ is used. 

We first consider the case $t\ge 0$.  By Lemma \ref{lem:phi-1}, 
\begin{align*}
	\phi_{q}(\frac{1+t}{2})+\phi_{q}(\frac{1-t}{2}) & \le2\conc\,\phi_{q}(\frac{1}{2})=2\phi_{q}(\frac{1}{2}),\\
	\phi_{q}(\frac{1+t}{2})-\phi_{q}(\frac{1-t}{2}) & \ge0,
\end{align*}
and by definition of $\phi_{q}$, 
\[
(1+t)\phi_{q}(\frac{1+t}{2})=(1-t)\phi_{q}(\frac{1-t}{2}).
\]
So,  the optimization in \eqref{eq:-4-1-3-1} is upper bounded by
\[
\frac{1-\rho}{2}\frac{1+\rho-4\rho^{2}\omega(\beta)}{1-\rho^{2}}\phi_{q}(\frac{1}{2}),
\]
which corresponds to the value of the objective function in \eqref{eq:-22-2}
at $t=0$, and hence, is no larger than the optimal value of the optimization
in \eqref{eq:-22-2}. 

We next consider the case $t<0$. 
For this case,  we rewrite 
\begin{align*}
	\eta_{q}(-t,t,\beta) & = \frac{A+Bt}{(1-\rho^{2})(1-t)} \phi_{q}(\frac{1+t}{2}) = \frac{A-Bs}{(1-\rho^{2})(1+s)} \phi_{q}(\frac{1-s}{2}),
\end{align*}
where $s=-t \in (0,1)$.
So,  the optimization in \eqref{eq:-4-1-3} is upper bounded by 
\[
\frac{1+\rho-4\rho^{2}\omega(\beta)-2\beta\rho^{2} s}{2(1+\rho)(1+s)}\phi_{q}(\frac{1-s}{2}).
\]

Summarizing all the above, the whole optimization in \eqref{eq:-4-1-2}
is equal to 
\begin{align*}
	& \max_{\beta\in[0,\frac{1}{2}-\epsilon^{*}(\rho)]}\max_{t\in[0,1)}\max \left\{\frac{\left(1-\rho\right)\left(1+\rho-4\rho^{2}\omega(\beta)\right)}{2\left(1+t-\rho^{2}\right)},\frac{1+\rho-4\rho^{2}\omega(\beta)-2\beta\rho^{2} t}{2(1+\rho)(1+t)}\right\} \phi_{q}(\frac{1-t}{2}) \nonumber \\
	& =\max_{t\in[0,1)}\max
	\left\{\eta_{q,1}(t),\eta_{q,2}(t)\right\}. 
\end{align*}

\section{\label{sec:Proof-of-Theorem-1-1-1}Proof of Theorem \ref{thm:CK-1}}

\begin{figure}
	\centering \includegraphics[scale=0.6]{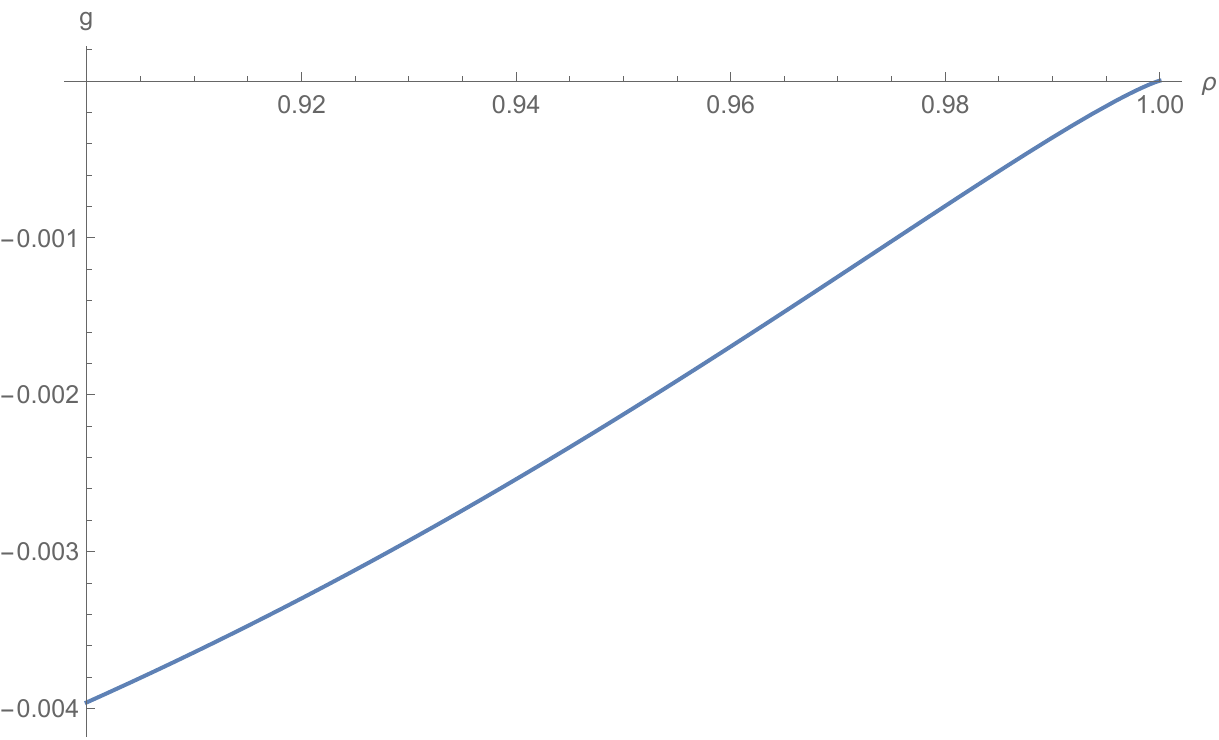}
	
	\caption{\label{fig:The-region-of-2}The function  $g$.}
\end{figure}

We need the following basic properties of $\epsilon_{q}^{*}(\rho)$,
$\omega_{q}^{*}(\rho)$, and $t_{q,\rho}$. 
\begin{lemma}
	\label{lem:The-following-hold.-2} Let $q=1.36$. Then, the following
	hold. 
	\begin{enumerate}
		\item For all $\rho\in[0.914,1]$, it holds that $\epsilon_{q}^{*}(\rho)\ge0.16$.
		\item For all $\rho\in[0.914,1]$, it holds that $\omega_{q}^{*}(\rho)\le\omega_{q}^{*}(1)\le\omega(\beta_{0})=0.193026...$ 
		\item For all $\rho\in[0.914,1]$, it holds that $t_{q,\rho}\in[0.41,0.91]$. 
	\end{enumerate}
\end{lemma}
\begin{proof}
	We first prove Statement 1. Consider the function 
	\begin{align*}
	f(\rho,\epsilon):&=\left(\epsilon+\left(\frac{1+\rho}{2}\right)^{p}(1-2\epsilon)\right)^{q/p}+\left(\epsilon+\left(\frac{1-\rho}{2}\right)^{p}(1-2\epsilon)\right)^{q/p} \\
	& \qquad -\left(\frac{1+\rho}{2}\right)^{q}-\left(\frac{1-\rho}{2}\right)^{q},
	\end{align*}
	with $p=1+(q-1)\rho^{2}$ and $q=1.36$. So, by Lemma \ref{lem:unique},
	given $\rho\in(0,1)$, $\epsilon_{q}^{*}(\rho)$ is the unique solution
	in $(0,1/2]$ to $f(\rho,\epsilon)=0$. Since $f(\rho,\epsilon)$
	is strictly convex in $\epsilon$ and $f(\rho,0)=f(\rho,\epsilon_{q}^{*}(\rho))=0$,
	it holds that $f(\rho,\epsilon)<0$ for all $\epsilon\in(0,\epsilon_{q}^{*}(\rho))$
	and $f(\rho,\epsilon)>0$ for all $\epsilon\in(\epsilon_{q}^{*}(\rho),1)$.
	So, to prove $\epsilon_{q}^{*}(\rho)\ge0.16$, it suffices to show
	that $f(\rho,0.16)\le0$. Note that $\rho\mapsto f(\rho,0.16)$ with
	$q=1.36$ is a one-variable function, whose graph is plotted in Fig.
	\ref{fig:The-region-of-2}. From the figure, it is easy to observe
	that $g(\rho):=f(\rho,0.16)\le0$ for all $\rho\in[0.914,1]$. We
	next prove this point rigorously. 
	
	We compute the derivative: 
	\begin{align*}
		g'(\rho) & =A_{+}+A_{-}+B_{+}+B_{-}+C,
	\end{align*}
	where 
	\begin{align*}
		A_{+} & =\frac{q}{p}(1-2\epsilon)\left((1-2\epsilon)(\frac{1+\rho}{2})^{p}+\epsilon\right)^{\frac{q}{p}-1}\left(2(q-1)\rho(\frac{1+\rho}{2})^{p}\ln\frac{1+\rho}{2}+\frac{p}{2}(\frac{1+\rho}{2})^{p-1}\right)\\
		A_{-} & =\frac{q}{p}(1-2\epsilon)\left((1-2\epsilon)(\frac{1-\rho}{2})^{p}+\epsilon\right)^{\frac{q}{p}-1}\left(2(q-1)\rho(\frac{1-\rho}{2})^{p}\ln\frac{1-\rho}{2}-\frac{p}{2}(\frac{1-\rho}{2})^{p-1}\right)\\
		B_{+} & =-\frac{2(q-1)q\rho}{p^{2}}\left((1-2\epsilon)(\frac{1+\rho}{2})^{p}+\epsilon\right)^{\frac{q}{p}}\ln\left((1-2\epsilon)(\frac{1+\rho}{2})^{p}+\epsilon\right)\\
		B_{-} & =-\frac{2(q-1)q\rho}{p^{2}}\left((1-2\epsilon)(\frac{1-\rho}{2})^{p}+\epsilon\right)^{\frac{q}{p}}\ln\left((1-2\epsilon)(\frac{1-\rho}{2})^{p}+\epsilon\right)\\
		C & =\frac{q}{2}(\frac{1-\rho}{2})^{q-1}-\frac{q}{2}(\frac{1+\rho}{2})^{q-1}.
	\end{align*}
	Since given $a>1$, $t\mapsto t^{a}\ln t$ is decreasing on $[0,e^{-1/a}]$
	and increasing on $[e^{-1/a},1]$, we have that 
	\begin{align*}
		B_{+} & \ge-\frac{2(q-1)q\rho}{p^{2}}\left(1-\epsilon\right)^{\frac{q}{p}}\ln\left(1-\epsilon\right),\\
		B_{-} & \ge-\frac{2(q-1)q\rho}{p^{2}}\epsilon^{\frac{q}{p}}\ln\epsilon,
	\end{align*}
	if 
	\begin{align}
		(1-2\epsilon)(\frac{1+\rho}{2})^{p}+\epsilon & \ge e^{-p/q},\label{eq:-39}\\
		(1-2\epsilon)(\frac{1-\rho}{2})^{p}+\epsilon & \le e^{-p/q}.\label{eq:-40}
	\end{align}
	We can also bound $A_{+}$ and $A_{-}$ as follows: 
	\begin{align*}
		A_{+} & =\frac{2(q-1)\rho q(1-2\epsilon)}{p}\left((1-2\epsilon)(\frac{1+\rho}{2})^{p}+\epsilon\right)^{\frac{q}{p}-1}(\frac{1+\rho}{2})^{p}\ln\frac{1+\rho}{2}\\
		& \qquad+\frac{q}{2}(1-2\epsilon)\left((1-2\epsilon)(\frac{1+\rho}{2})^{p}+\epsilon\right)^{\frac{q}{p}-1}(\frac{1+\rho}{2})^{p-1}\\
		& \ge\frac{2(q-1)\rho q(1-2\epsilon)}{p}(\frac{1+\rho}{2})^{p}\ln\frac{1+\rho}{2} +\frac{q}{2}(1-2\epsilon)\epsilon^{\frac{q}{p}-1}(\frac{1+\rho}{2})^{p-1},\\
		A_{-} & \ge\frac{q}{p}(1-2\epsilon)\left(2(q-1)\rho(\frac{1-\rho}{2})^{p}\ln\frac{1-\rho}{2}-\frac{p}{2}(\frac{1-\rho}{2})^{p-1}\right).
	\end{align*}
	Moreover, $C\ge-\frac{q}{2}.$ Hence, for all $\rho$ satisfying
	\eqref{eq:-39} and \eqref{eq:-40}, it holds that 
	\begin{align}
		g'(\rho)&\ge\frac{q}{p}(1-2\epsilon)\left(\epsilon^{\frac{q}{p}-1}\frac{p}{2}(\frac{1+\rho}{2})^{p-1}-\frac{p}{2}(\frac{1-\rho}{2})^{p-1}-2(q-1)\rho H_{p}(\frac{1-\rho}{2})\right) \nonumber\\
		&\qquad +\frac{2(q-1)q\rho}{p^{2}}H_{q/p}(\epsilon)-\frac{q}{2},\label{eq:-42}
	\end{align}
	where $H_{a}(t):=-t^{a}\ln t-(1-t)^{a}\ln(1-t)$. Let $\rho_{0}$
	be a number such that 
	\begin{align}
		(1-2\epsilon)(\frac{1+\rho_{0}}{2})^{q}+\epsilon & \ge e^{-p_{0}/q},\label{eq:-39-2}\\
		(1-2\epsilon)(\frac{1-\rho_{0}}{2})^{p_{0}}+\epsilon & \le e^{-1},\label{eq:-40-2}
	\end{align}
	and 
	\begin{align}
		&\frac{q}{2}(1-2\epsilon)\left(\epsilon^{\frac{q}{p_{0}}-1}(\frac{1+\rho_{0}}{2})^{q-1}-(\frac{1-\rho_{0}}{2})^{p_{0}-1}-\frac{4(q-1)}{p_{0}}H_{1}(\frac{1-\rho_{0}}{2})\right) \nonumber\\
		& \qquad +\frac{2(q-1)\rho_{0}}{q}H_{q/p_{0}}(\epsilon)-\frac{q}{2}\ge0,\label{eq:-41}
	\end{align}
	where $p_{0}=1+(q-1)\rho_{0}^{2}$. Then, the inequality in \eqref{eq:-42}
	implies that for all $\rho\in[\rho_{0},1]$, $g'(\rho)\ge0$. Using
	computers, it is easy to verify that $\rho_{0}=0.9999$ is a good
	choice which satisfies \eqref{eq:-39-2}--\eqref{eq:-41}. So, for
	all $\rho\in[0.9999,1]$, $g'(\rho)\ge0$. On the other hand, $g(1)=0$.
	So, for all $\rho\in[0.9999,1]$, $g(\rho)\le0$. 
	
	We next show $g(\rho)\le0$ for all $\rho\in[0.914,0.9999]$. Denote
	\begin{align*}
		F(\rho) & =\left(\epsilon+\left(\frac{1+\rho}{2}\right)^{p}(1-2\epsilon)\right)^{q/p}+\left(\epsilon+\left(\frac{1-\rho}{2}\right)^{p}(1-2\epsilon)\right)^{q/p},\\
		G(\rho) & =\left(\frac{1+\rho}{2}\right)^{q}+\left(\frac{1-\rho}{2}\right)^{q},
	\end{align*}
	with $p=1+(q-1)\rho^{2}$, $q=1.36$, and $\epsilon=0.16$. Then,
	$g(\rho)=F(\rho)-G(\rho)$. For any $\rho\in[\rho_{1},\rho_{2}]\subset[1/2,1]$,
	it holds that 
	\begin{align*}
		g(\rho) & \le\left(\epsilon+\left(\frac{1+\rho}{2}\right)^{p_{2}}(1-2\epsilon)\right)^{q/p_{2}}+\left(\epsilon+\left(\frac{1-\rho}{2}\right)^{p_{2}}(1-2\epsilon)\right)^{q/p_{2}}-G(\rho_{1})\\
		& \le\left(\epsilon+\left(\frac{1+\rho_{2}}{2}\right)^{p_{2}}(1-2\epsilon)\right)^{q/p_{2}}+\left(\epsilon+\left(\frac{1-\rho_{2}}{2}\right)^{p_{2}}(1-2\epsilon)\right)^{q/p_{2}}-G(\rho_{1})\\
		& =F(\rho_{2})-G(\rho_{1}),
	\end{align*}
	where $p_{2}=1+(q-1)\rho_{2}^{2}$. Here, the first inequality follows
	since $G$ is convex and $p$-norm in a probability space is increasing
	in $p$. The second inequality follows by the convexity of the expression
	at the first line in $\rho$.  We run $\rho_{1}$ and $\rho_{2}$
	over all two consecutive numbers in the finite set $[0.914:0.000002:0.9999]$.
	We numerically verify that $F(\rho_{2})-G(\rho_{1})<0$. So, $g(\rho)\le0$
	for all $\rho\in[0.914,0.9999]$. The Matlab code for this verification is included within the TeX source files of this arXiv version.
	
	Statement 2 can be verified by using Statement 1 and Lemma \ref{lem:omega}. 
	
	We now prove Statement 3. Observe that 
	\begin{align*}
		\phi_{q}'(\frac{1-t_{q,\rho}}{2})   =\frac{-2\phi_{q}(\frac{1-\rho}{2})}{1+\rho-4\rho^{2}\omega_{q}^{*}(\rho)}  \le-2\phi_{q}(\frac{1-\rho}{2})  \le-2\phi_{q}(0) =\frac{2q}{q-1},
	\end{align*}
	which, combined with the monotonicity of $\phi_{q}'$ (given in Lemma
	\ref{lem:phi-1}), implies $t_{q,\rho}\le0.91$ for all $\rho\in[0.914,1]$.
	Denoting $s=\frac{1-\rho}{2}$, we observe that 
	\begin{align*}
		\phi_{q}'(\frac{1-t_{q,\rho}}{2})   =\frac{-2\phi_{q}(\frac{1-\rho}{2})}{1+\rho-4\rho^{2}\omega_{q}^{*}(\rho)}  \ge\frac{-2\phi_{q}(\frac{1-\rho}{2})}{1+\rho} \ge\frac{-\phi_{q}(s)}{1-s}\Big|_{s=\frac{1-0.914}{2}}  =\frac{-\phi_{q}(0.043)}{0.957},
	\end{align*}
	which, combined with the monotonicity of $\phi_{q}'$ (given in Lemma
	\ref{lem:phi-1}), implies $t_{q,\rho}\ge0.41$ for all $\rho\in[0.914,1]$. 
\end{proof}

We next prove Theorem \ref{thm:CK-1}. We first swap two maximizations
in \eqref{eq:max} and then apply Lemma \ref{lem:The-following-hold.-2}
to obtain that 
\begin{equation}
	\max_{\rho\in[0.914,1]}\max_{t\in[0,1]}\chi_{q}(\rho,t)=\max_{t\in[0.41,0.91]}\max_{\rho\in[0.914,1]}\chi_{q}(\rho,t).\label{eq:-4-2-3-1}
\end{equation}
Observe that 
\begin{align*}
	\partial_{\rho}\chi_{q}(\rho,t) & =\left(1-8\rho\omega_{0}\right)\phi_{q}(\frac{1-t}{2})+2\rho\phi_{q}(\frac{1-\rho}{2})+\frac{1}{2}\left(1+t-\rho^{2}\right)\phi_{q}'(\frac{1-\rho}{2}).
\end{align*}
We bound this derivative as follows. For $q=1.36$,
\begin{align*}
	\partial_{\rho}\chi_{q}(\rho,t) & \ge-\left(|1-8\rho\omega_{0}|\left|\phi_{q}(\frac{1-t}{2})\right|+2\rho\left|\phi_{q}(\frac{1-\rho}{2})\right|\right)+\frac{1}{2}\left(1+t-\rho^{2}\right)\phi_{q}'(\frac{1-\rho}{2})\\
	& \ge2\phi_{q}(0)+\frac{0.41}{2}\phi_{q}'(\frac{1-\rho}{2})  =-\frac{2q}{q-1}+0.205\phi_{q}'(\frac{1-\rho}{2})\\
	& =-\frac{2.72}{0.36}+0.205\phi_{q}'(\frac{1-\rho}{2}).
\end{align*}
Since for all $\rho\ge0.993$, it holds that 
\[
\phi_{q}'(\frac{1-\rho}{2})\ge\frac{2.72}{0.36\times0.205},
\]
which implies that for all $\rho\ge0.993$, $\partial_{\rho}\chi_{q}(\rho,t)\ge0$.
That is, given any $t\in[0.41,0.91]$, $\chi_{q}(\rho,t)$ is increasing
in $\rho\in[0.993,1]$. So, for this case, the maximum of $\chi_{q}(\rho,t)$
is attained at $(\rho,t)=(1,t_{q,1})$.  It can be easily verified
that $\chi_{q}(1,t_{q,1})\le0$ for $q=1.36$. 

We next consider the case of $\rho\le0.993$. For all $\rho\le0.993$,
\begin{align*}
	|\partial_{\rho}\chi_{q}(\rho,t)| & \le|1-8\rho\omega_{0}|\left|\phi_{q}(\frac{1-t}{2})\right|+\max\left\{ 2\rho\left|\phi_{q}(\frac{1-\rho}{2})\right|,\frac{1}{2}\left(1+t-\rho^{2}\right)\left|\phi_{q}'(\frac{1-\rho}{2})\right|\right\} \\
	& \le\left|\phi_{q}(\frac{1-0.91}{2})\right|+\max\left\{ 2\left|\phi_{q}(0)\right|,\frac{1}{2}\left(1+0.91-0.914^{2}\right)\left|\phi_{q}'(\frac{1-0.993}{2})\right|\right\} \\
	& \le\left|\phi_{q}(0.045)\right|+\max\left\{ \frac{2.72}{0.36},0.537302\times\phi_{q}'(0.0035)\right\} \\
	& \le25=:M.
\end{align*}
We choose $\delta=0.2$. Then, $\Delta:=\frac{\delta}{M}=0.008$.
We finally use Matlab to verify $\chi_{q}(\rho,t)<-\delta$ numerically
for all $t\in[0.41,0.91]$ and $\rho\in\{0.914,0.922,0.930,...,0.986\}$,
which, by \eqref{eq:-24}, yields that $\chi_{q}(\hat{\rho},t)<0$
for all $t\in[0.41,0.91]$ and $\hat{\rho}\in[0.914,0.993]$. The Matlab code for this verification is included within the TeX source files of this arXiv version.

\begin{acks}[Acknowledgments]
The authors would like to thank the anonymous referees and the Editor Qi-Man Shao for their constructive comments that improved the
quality of this paper, and especially thank one of referees for pointing out that our proof cannot be extended to the case $q<0$ in the present form.
\end{acks}

\begin{funding}
This work was supported by the National Key Research and Development
Program of China under grant 2023YFA1009604 and the NSFC under grant
62101286.
\end{funding}

\bibliographystyle{imsart-number} 
\bibliography{ref}       

\begin{thebibliography}{37}

\bibitem{anantharam2017conjecture}
\begin{binproceedings}[author]
\bauthor{\bsnm{Anantharam},~\bfnm{V.}\binits{V.}},
  \bauthor{\bsnm{Bogdanov},~\bfnm{A.}\binits{A.}},
  \bauthor{\bsnm{Chakrabarti},~\bfnm{A.}\binits{A.}},
  \bauthor{\bsnm{Jayaram},~\bfnm{T.}\binits{T.}} \AND
  \bauthor{\bsnm{Nair},~\bfnm{C.}\binits{C.}}
(\byear{2017}).
\btitle{A conjecture regarding optimality of the dictator function under
  {Hellinger} distance}.
In \bbooktitle{Information Theory and Applications Workshop)}.
\bpublisher{Available at
  http://chandra.ie.cuhk.edu.hk/pub/papers/HC/hel-conj.pdf}.
\end{binproceedings}
\endbibitem

\bibitem{anantharam2013on}
\begin{binproceedings}[author]
\bauthor{\bsnm{Anantharam},~\bfnm{V.}\binits{V.}},
  \bauthor{\bsnm{Gohari},~\bfnm{A.~A.}\binits{A.~A.}},
  \bauthor{\bsnm{Kamath},~\bfnm{S.}\binits{S.}} \AND
  \bauthor{\bsnm{Nair},~\bfnm{C.}\binits{C.}}
(\byear{2013}).
\btitle{On hypercontractivity and the mutual information between {Boolean}
  functions}.
In \bbooktitle{Communication, Control, and Computing (Allerton), 2013 51th
  Annual Allerton Conference on}
\bpages{13--19}.
\bpublisher{IEEE}.
\end{binproceedings}
\endbibitem

\bibitem{barnes2020courtade}
\begin{binproceedings}[author]
\bauthor{\bsnm{Barnes},~\bfnm{L.~P.}\binits{L.~P.}} \AND
  \bauthor{\bsnm{{\"O}zg{\"u}r},~\bfnm{A.}\binits{A.}}
(\byear{2020}).
\btitle{The {Courtade-Kumar} Most Informative {Boolean} Function Conjecture and
  a Symmetrized {Li-M\'edard} Conjecture are Equivalent}.
In \bbooktitle{2020 IEEE International Symposium on Information Theory (ISIT)}
\bpages{2205--2209}.
\bpublisher{IEEE}.
\end{binproceedings}
\endbibitem

\bibitem{beltran2023sharp}
\begin{barticle}[author]
\bauthor{\bsnm{Beltran},~\bfnm{D.}\binits{D.}},
  \bauthor{\bsnm{Ivanisvili},~\bfnm{P.}\binits{P.}} \AND
  \bauthor{\bsnm{Madrid},~\bfnm{J.}\binits{J.}}
(\byear{2023}).
\btitle{On sharp isoperimetric inequalities on the hypercube}.
\bjournal{arXiv preprint arXiv:2303.06738}.
\end{barticle}
\endbibitem

\bibitem{bobkov1999discrete}
\begin{barticle}[author]
\bauthor{\bsnm{Bobkov},~\bfnm{S.~G.}\binits{S.~G.}} \AND
  \bauthor{\bsnm{G{\"o}tze},~\bfnm{F.}\binits{F.}}
(\byear{1999}).
\btitle{Discrete isoperimetric and {Poincar{\'e}-type} inequalities}.
\bjournal{Probability theory and related fields}
\bvolume{114}
\bpages{245--277}.
\end{barticle}
\endbibitem

\bibitem{borell1985geometric}
\begin{barticle}[author]
\bauthor{\bsnm{Borell},~\bfnm{C.}\binits{C.}}
(\byear{1985}).
\btitle{Geometric bounds on the {Ornstein--Uhlenbeck} velocity process}.
\bjournal{Probability Theory and Related Fields}
\bvolume{70}
\bpages{1--13}.
\end{barticle}
\endbibitem

\bibitem{chang2002polynomial}
\begin{barticle}[author]
\bauthor{\bsnm{Chang},~\bfnm{M.~C.}\binits{M.~C.}}
(\byear{2002}).
\btitle{A polynomial bound in {Freiman's} theorem}.
\bjournal{Duke mathematical journal}
\bvolume{113}
\bpages{399--419}.
\end{barticle}
\endbibitem

\bibitem{chen2025differential}
\begin{barticle}[author]
\bauthor{\bsnm{Chen},~\bfnm{Zijie}\binits{Z.}},
  \bauthor{\bsnm{Gohari},~\bfnm{Amin}\binits{A.}} \AND
  \bauthor{\bsnm{Nair},~\bfnm{Chandra}\binits{C.}}
(\byear{2025}).
\btitle{A Differential Equation Approach to the Most-Informative Boolean
  Function Conjecture}.
\bjournal{arXiv preprint arXiv:2502.10019}.
\end{barticle}
\endbibitem

\bibitem{chen2024optimality}
\begin{binproceedings}[author]
\bauthor{\bsnm{Chen},~\bfnm{Z.}\binits{Z.}} \AND
  \bauthor{\bsnm{Nair},~\bfnm{C.}\binits{C.}}
(\byear{2024}).
\btitle{On the Optimality of Dictator Functions and Isoperimetric Inequalities
  on Boolean Hypercubes}.
In \bbooktitle{2024 IEEE International Symposium on Information Theory (ISIT)}
\bpages{3380--3385}.
\bpublisher{IEEE}.
\end{binproceedings}
\endbibitem

\bibitem{courtade2014boolean}
\begin{barticle}[author]
\bauthor{\bsnm{Courtade},~\bfnm{T.~A.}\binits{T.~A.}} \AND
  \bauthor{\bsnm{Kumar},~\bfnm{G.~R.}\binits{G.~R.}}
(\byear{2014}).
\btitle{Which {Boolean} functions maximize mutual information on noisy inputs?}
\bjournal{IEEE Trans. Inf. Theory}
\bvolume{60}
\bpages{4515--4525}.
\end{barticle}
\endbibitem

\bibitem{durcik2024sharp}
\begin{barticle}[author]
\bauthor{\bsnm{Durcik},~\bfnm{P.}\binits{P.}},
  \bauthor{\bsnm{Ivanisvili},~\bfnm{P.}\binits{P.}} \AND
  \bauthor{\bsnm{Roos},~\bfnm{J.}\binits{J.}}
(\byear{2024}).
\btitle{Sharp isoperimetric inequalities on the {Hamming} cube near the
  critical exponent}.
\bjournal{arXiv preprint arXiv:2407.12674}.
\end{barticle}
\endbibitem

\bibitem{eldan2015two}
\begin{barticle}[author]
\bauthor{\bsnm{Eldan},~\bfnm{R.}\binits{R.}}
(\byear{2015}).
\btitle{A two-sided estimate for the {Gaussian} noise stability deficit}.
\bjournal{Inventiones Mathematicae}
\bvolume{201}
\bpages{561--624}.
\end{barticle}
\endbibitem

\bibitem{eldan2022noise}
\begin{barticle}[author]
\bauthor{\bsnm{Eldan},~\bfnm{R.}\binits{R.}},
  \bauthor{\bsnm{Mikulincer},~\bfnm{D.}\binits{D.}} \AND
  \bauthor{\bsnm{Raghavendra},~\bfnm{P.}\binits{P.}}
(\byear{2023}).
\btitle{Noise stability on the {Boolean} hypercube via a renormalized
  {Brownian} motion}.
\bpages{661--671}.
\end{barticle}
\endbibitem

\bibitem{folland_real_1999}
\begin{bbook}[author]
\bauthor{\bsnm{Folland},~\bfnm{Gerald~B.}\binits{G.~B.}}
(\byear{1999}).
\btitle{Real {Analysis}: {Modern} Techniques and Their Applications},
\bedition{2} ed.
\bpublisher{John Wiley \& Sons, Inc.}, \baddress{New York, NY, USA}.
\end{bbook}
\endbibitem

\bibitem{fu2001minimum}
\begin{barticle}[author]
\bauthor{\bsnm{Fu},~\bfnm{F.~W.}\binits{F.~W.}},
  \bauthor{\bsnm{Wei},~\bfnm{V.~K.}\binits{V.~K.}} \AND
  \bauthor{\bsnm{Yeung},~\bfnm{R.~W.}\binits{R.~W.}}
(\byear{2001}).
\btitle{On the minimum average distance of binary codes: {Linear} programming
  approach}.
\bjournal{Discrete Applied Mathematics}
\bvolume{111}
\bpages{263--281}.
\end{barticle}
\endbibitem

\bibitem{gacs1973common}
\begin{barticle}[author]
\bauthor{\bsnm{G{\'a}cs},~\bfnm{P.}\binits{P.}} \AND
  \bauthor{\bsnm{K{\"o}rner},~\bfnm{J.}\binits{J.}}
(\byear{1973}).
\btitle{Common information is far less than mutual information}.
\bjournal{Problems of Control and Information Theory}
\bvolume{2}
\bpages{149--162}.
\end{barticle}
\endbibitem

\bibitem{kahn2020isoperimetric}
\begin{barticle}[author]
\bauthor{\bsnm{Kahn},~\bfnm{J.}\binits{J.}} \AND
  \bauthor{\bsnm{Park},~\bfnm{J.}\binits{J.}}
(\byear{2020}).
\btitle{An isoperimetric inequality for the {Hamming} cube and some
  consequences}.
\bjournal{Proceedings of the American Mathematical Society}
\bvolume{148}
\bpages{4213--4224}.
\end{barticle}
\endbibitem

\bibitem{kindler2015remarks}
\begin{barticle}[author]
\bauthor{\bsnm{Kindler},~\bfnm{G.}\binits{G.}},
  \bauthor{\bsnm{O'Donnell},~\bfnm{R.}\binits{R.}} \AND
  \bauthor{\bsnm{Witmer},~\bfnm{D.}\binits{D.}}
(\byear{2015}).
\btitle{Remarks on the most informative function conjecture at fixed mean}.
\bjournal{arXiv preprint arXiv:1506.03167}.
\end{barticle}
\endbibitem

\bibitem{laguerre1883sur}
\begin{barticle}[author]
\bauthor{\bsnm{Laguerre},~\bfnm{E.~N.}\binits{E.~N.}}
(\byear{1883}).
\btitle{Sur la th\'eorie des \'equations num\'eriques}.
\bjournal{Journal de Math\'ematiques pures et appliqu\'ees. [Online].
  Available: http://sepwww.stanford.edu/oldsep/stew/laguerre.pdf}.
\end{barticle}
\endbibitem

\bibitem{ledoux1994semigroup}
\begin{barticle}[author]
\bauthor{\bsnm{Ledoux},~\bfnm{M.}\binits{M.}}
(\byear{1994}).
\btitle{Semigroup proofs of the isoperimetric inequality in {Euclidean} and
  {Gauss} space}.
\bjournal{Bulletin des sciences math{\'e}matiques}
\bvolume{118}
\bpages{485--510}.
\end{barticle}
\endbibitem

\bibitem{ledoux2014remarks}
\begin{binproceedings}[author]
\bauthor{\bsnm{Ledoux},~\bfnm{M.}\binits{M.}}
(\byear{2014}).
\btitle{Remarks on {Gaussian} noise stability, {Brascamp-Lieb} and {Slepian}
  inequalities}.
In \bbooktitle{Geometric Aspects of Functional Analysis: Israel Seminar (GAFA)
  2011-2013}
\bpages{309--333}.
\bpublisher{Springer}.
\end{binproceedings}
\endbibitem

\bibitem{li2020boolean}
\begin{barticle}[author]
\bauthor{\bsnm{Li},~\bfnm{J.}\binits{J.}} \AND
  \bauthor{\bsnm{M{\'e}dard},~\bfnm{M.}\binits{M.}}
(\byear{2020}).
\btitle{Boolean functions: noise stability, non-interactive correlation
  distillation, and mutual information}.
\bjournal{IEEE Trans. Inf. Theory}
\bvolume{67}
\bpages{778--789}.
\end{barticle}
\endbibitem

\bibitem{marshall2011inequalities}
\begin{bbook}[author]
\bauthor{\bsnm{Marshall},~\bfnm{A.~W.}\binits{A.~W.}},
  \bauthor{\bsnm{Olkin},~\bfnm{I.}\binits{I.}} \AND
  \bauthor{\bsnm{Arnold},~\bfnm{B.~C.}\binits{B.~C.}}
(\byear{2011}).
\btitle{Inequalities: theory of majorization and its applications}.
\bpublisher{Springer, New York}.
\end{bbook}
\endbibitem

\bibitem{mossel2005coin}
\begin{barticle}[author]
\bauthor{\bsnm{Mossel},~\bfnm{E.}\binits{E.}} \AND
  \bauthor{\bsnm{O'Donnell},~\bfnm{R.}\binits{R.}}
(\byear{2005}).
\btitle{Coin flipping from a cosmic source: {On} error correction of truly
  random bits}.
\bjournal{Random Structures \& Algorithms}
\bvolume{26}
\bpages{418--436}.
\end{barticle}
\endbibitem

\bibitem{nair2014equivalent}
\begin{binproceedings}[author]
\bauthor{\bsnm{Nair},~\bfnm{C.}\binits{C.}}
(\byear{2014}).
\btitle{Equivalent formulations of hypercontractivity using information
  measures}.
In \bbooktitle{International Zurich Seminar (IZS) Workshop}.
\end{binproceedings}
\endbibitem

\bibitem{O'Donnell14analysisof}
\begin{bbook}[author]
\bauthor{\bsnm{O'Donnell},~\bfnm{R.}\binits{R.}}
(\byear{2014}).
\btitle{Analysis of {Boolean} Functions}.
\bpublisher{Cambridge University Press}.
\end{bbook}
\endbibitem

\bibitem{ordentlich2016improved}
\begin{binproceedings}[author]
\bauthor{\bsnm{Ordentlich},~\bfnm{O.}\binits{O.}},
  \bauthor{\bsnm{Shayevitz},~\bfnm{O.}\binits{O.}} \AND
  \bauthor{\bsnm{Weinstein},~\bfnm{O.}\binits{O.}}
(\byear{2016}).
\btitle{An improved upper bound for the most informative {Boolean} function
  conjecture}.
In \bbooktitle{2016 IEEE International Symposium on Information Theory (ISIT)}
\bpages{500--504}.
\bpublisher{IEEE}.
\end{binproceedings}
\endbibitem

\bibitem{pichler2018dictator}
\begin{barticle}[author]
\bauthor{\bsnm{Pichler},~\bfnm{G.}\binits{G.}},
  \bauthor{\bsnm{Piantanida},~\bfnm{P.}\binits{P.}} \AND
  \bauthor{\bsnm{Matz},~\bfnm{G.}\binits{G.}}
(\byear{2018}).
\btitle{Dictator functions maximize mutual information}.
\bjournal{The Annals of Applied Probability}
\bvolume{28}
\bpages{3094--3101}.
\end{barticle}
\endbibitem

\bibitem{polyanskiy2010channel}
\begin{barticle}[author]
\bauthor{\bsnm{Polyanskiy},~\bfnm{Y.}\binits{Y.}},
  \bauthor{\bsnm{Poor},~\bfnm{H.~V.}\binits{H.~V.}} \AND
  \bauthor{\bsnm{Verd{\'u}},~\bfnm{S.}\binits{S.}}
(\byear{2010}).
\btitle{Channel coding rate in the finite blocklength regime}.
\bjournal{IEEE Trans. Inf. Theory}
\bvolume{56}
\bpages{2307--2359}.
\end{barticle}
\endbibitem

\bibitem{samorodnitsky2016entropy}
\begin{barticle}[author]
\bauthor{\bsnm{Samorodnitsky},~\bfnm{A.}\binits{A.}}
(\byear{2016}).
\btitle{On the entropy of a noisy function}.
\bjournal{IEEE Trans. Inf. Theory}
\bvolume{62}
\bpages{5446--5464}.
\end{barticle}
\endbibitem

\bibitem{talagrand1993isoperimetry}
\begin{barticle}[author]
\bauthor{\bsnm{Talagrand},~\bfnm{M.}\binits{M.}}
(\byear{1993}).
\btitle{Isoperimetry, logarithmic {Sobolev} inequalities on the discrete cube,
  and {Margulis'} graph connectivity theorem}.
\bjournal{Geometric \& Functional Analysis GAFA}
\bvolume{3}
\bpages{295--314}.
\end{barticle}
\endbibitem

\bibitem{tsallis1994numbers}
\begin{barticle}[author]
\bauthor{\bsnm{Tsallis},~\bfnm{C.}\binits{C.}}
(\byear{1994}).
\btitle{What are the numbers that experiments provide}.
\bjournal{Quimica Nova}
\bvolume{17}
\bpages{468--471}.
\end{barticle}
\endbibitem

\bibitem{witsenhausen1975sequences}
\begin{barticle}[author]
\bauthor{\bsnm{Witsenhausen},~\bfnm{H.~S.}\binits{H.~S.}}
(\byear{1975}).
\btitle{On sequences of pairs of dependent random variables}.
\bjournal{SIAM Journal on Applied Mathematics}
\bvolume{28}
\bpages{100--113}.
\end{barticle}
\endbibitem

\bibitem{yu2023phi}
\begin{barticle}[author]
\bauthor{\bsnm{Yu},~\bfnm{L.}\binits{L.}}
(\byear{2023}).
\btitle{On the {$\Phi$}-Stability and Related Conjectures}.
\bjournal{Probability Theory and Related Fields}
\bvolume{186}
\bpages{1045--1080}.
\end{barticle}
\endbibitem

\bibitem{yu2022theentro}
\begin{bbook}[author]
\bauthor{\bsnm{Yu},~\bfnm{L.}\binits{L.}}
(\byear{2023}).
\btitle{The Entropy Method}.
\bpublisher{DOI: 10.13140/RG.2.2.26552.11527/1}.
\end{bbook}
\endbibitem

\bibitem{yu2019improved}
\begin{barticle}[author]
\bauthor{\bsnm{Yu},~\bfnm{L.}\binits{L.}} \AND
  \bauthor{\bsnm{Tan},~\bfnm{V.~Y.~F.}\binits{V.~Y.~F.}}
(\byear{2019}).
\btitle{An Improved Linear Programming Bound on the Average Distance of a
  Binary Code}.
\bjournal{arXiv preprint arXiv:1910.09416}.
\end{barticle}
\endbibitem

\bibitem{yu2022common}
\begin{barticle}[author]
\bauthor{\bsnm{Yu},~\bfnm{L.}\binits{L.}} \AND
  \bauthor{\bsnm{Tan},~\bfnm{V.~Y.~F.}\binits{V.~Y.~F.}}
(\byear{2022}).
\btitle{Common information, noise stability, and their extensions}.
\bjournal{Foundations and Trends in Communications and Information Theory}
\bvolume{19}
\bpages{107--389}.
\end{barticle}
\endbibitem

\end{thebibliography}

%
%
%

\end{document}